\documentclass[opre,nonblindrev]{informs3} 

\DoubleSpacedXI 

\usepackage{xcolor}

\usepackage{enumitem}

\usepackage{natbib}
\bibpunct[, ]{(}{)}{,}{a}{}{,}%
\def\bibfont{\small}%

\usepackage{bbm}
\usepackage{endnotes}

\TheoremsNumberedThrough     
\ECRepeatTheorems

\EquationsNumberedThrough    

\usepackage[toc,page]{appendix}
\usepackage{subcaption}
\usepackage{mathtools}
\usepackage{bigstrut}
\usepackage{booktabs}
\usepackage{mathrsfs}

\usepackage[utf8x]{inputenc}
\usepackage{multicol}

\usepackage{algorithm}
\usepackage{algorithmicx}
\usepackage{algpseudocode}
\usepackage{adjustbox}
\usepackage{amsmath, cases, bm}
 \usepackage{relsize}
 \usepackage{hyperref}
 \usepackage{graphicx}
\usepackage{multirow}
\usepackage{comment}
\usepackage{dsfont}

\newcommand{\st}{{\mathrm{s.t.}}}

\newcommand\scalemath[2]{\scalebox{#1}{\mbox{\ensuremath{\displaystyle #2}}}}

\DeclareMathOperator{\expec}{\ensuremath{\mathbb{E}}}

\DeclareMathOperator{\tr}{\ensuremath{\mathrm{Tr}}}

\DeclareMathOperator{\realR}{\ensuremath{\mathbb{R}}}

\usepackage{footnote}
\makesavenoteenv{tabular}
\makesavenoteenv{table}
\usepackage{tablefootnote}
\usepackage{threeparttable}
\usepackage{color, colortbl}
\usepackage[hang,flushmargin,symbol]{footmisc} 
\usepackage{subcaption}
\DeclareMathAlphabet{\mathpzc}{OT1}{pzc}{m}{it}

\begin{document}
	\graphicspath{{figures/}}
	
	
	\RUNAUTHOR{Goyal, Zhang, and He}
	
	\RUNTITLE{Distributionally Robust Bilevel Programs}
	
\TITLE{Decision Rule Approaches for Pessimistic Bilevel Linear Programs under Moment Ambiguity with Facility Location Applications}

\ARTICLEAUTHORS{%
	\AUTHOR{Akshit Goyal}
	\AFF{Department of Industrial and Systems Engineering, University of Minnesota, \EMAIL{goyal080@umn.edu
	}}
	\AUTHOR{Yiling Zhang}
	\AFF{Department of Industrial and Systems Engineering, University of Minnesota, \EMAIL{yiling@umn.edu}}
	\AUTHOR{Chuan He}
	\AFF{Department of Industrial and Systems Engineering, University of Minnesota, \EMAIL{he000233@umn.edu}}
}
	
	\ABSTRACT{%
	We study a pessimistic stochastic bilevel program in the context of sequential two-player games, where the leader makes a binary here-and-now decision, and the follower responds a continuous wait-and-see decision after observing the leader’s action and revelation of uncertainty. Only the information of the mean, covariance, and support is known. We formulate the problem as a distributionally robust (DR) two-stage problem. The pessimistic DR bilevel program is shown to be equivalent to a generic two-stage distributionally robust stochastic (nonlinear) program with both a random objective and random constraints under proper conditions of ambiguity sets. Under continuous distributions, using linear decision rule approaches, we construct upper bounds on the pessimistic DR bilevel program based on (1) 0-1 semidefinite programming (SDP) approximation and (2) an exact 0-1 copositive programming reformulations. When the ambiguity set is restricted to discrete distributions, an exact 0-1 SDP reformulation is developed, and explicit construction of the worst-case distribution is derived.  To further improve the computation of the proposed 0-1 SDPs, a cutting-plane framework is developed. Moreover, based on a mixed-integer linear programming approximation, another cutting-plane algorithm is proposed. Extensive numerical studies are conducted to demonstrate the effectiveness of the proposed approaches on a facility location problem.
	 
	
				
	}%
	
	
	
	\maketitle

\normalsize


\section{Introduction}
Two-stage stochastic bilevel programming often arises in the context of sequential two-player games where players strive to optimize their individual objectives under uncertainty. In the game, the leader makes a decision first to optimize their utility function, and then the follower responds after observing the uncertain parameters and the leader's action. Specifically, the leader first chooses a here-and-now decision $x\in\mathcal{X}\subseteq\{0,1\}^{d}$, before the revelation of uncertain parameter $\xi\in\realR^k$. Next, after the revelation, the follower makes a wait-and-see decision $y\in\realR^{n}$ to minimize their utility function $c(\xi)^\top y$ subject to $Ay\le b_x(\xi)$ for some $A\in\realR^{m\times n}$ and $b_x(\xi)\in\realR^{m}$. 
We emphasize that the leader's decision affects the follower's feasible region as the right-hand side $b_x(\xi)$ depends on the leader's decision $x$. By taking the follower's potential choice $y$ into account, the leader minimizes their direct cost $w^\top x$ and an indirect expected  cost $\mathbb E_F[v(\xi)^\top y]$ via the impact on the follower's feasible region. That is, the leader minimizes $w^{\top}x+\expec_F[ v(\xi)^\top y]$ for some $w\in\realR^d$ and $v(\xi)\in\realR^n$, where $\expec_F[\cdot]$ denotes the expectation with respect to distribution $F$ of $\xi$. 
\par

One important concern in solving bilevel optimization problems regards the follower's decision when multiple alternative optimal solutions are present.  The leader's decision might have very different values depending on the choice of follower's optimal solutions. Denote the set of alternative optimal solutions for the follower $\Omega(x, \xi) := \argmin_{y}\left\{ c(\xi)^\top y: \ Ay \le b_x(\xi)\right\} \subset \mathbb R^n$.
%
First, assuming that the follower is cooperative, i.e., the follower chooses the solution in favor of the leader, the \emph{optimistic} stochastic bilevel program is formulated as
\begin{equation}\label{model:OPT-STO-BILVL}
	\min_{x\in\mathcal X} w^\top x + \mathbb E_F \left[\min_{ y\in {\Omega}(x,\xi)}  v(\xi)^\top y \right].
\end{equation}
%
%
	%
The solution of \eqref{model:OPT-STO-BILVL} can be risky for the leader assuming cooperation on the follower. 
The \emph{pessimistic} formulation considers the worst-case scenario in the follower's optimal solution set $\Omega(x,\xi)$ as follows.
\begin{equation}\label{model:PESS-STO-BILVL}
	\min_{x\in\mathcal X} w^\top x + \mathbb E_F \left[ \max_{y\in \Omega(x,\xi)} v(\xi)^\top y \right].
\end{equation}
%
	%
	%
	Both the optimistic and the pessimistic formulations have been recently studied in
	\citet{yanikoglu2018decision}. Using linear decision rules (LDRs), they propose an upper-bound mixed-integer linear programming (MILP) approximation of the pessimistic  problem and a lower-bound MILP approximation of  the optimistic problem.  The approximations are robust with respect to distributions matching the first- and second-order (empirical) moments.

	However, there may exist estimation errors in the moments when, for example, the historical data is inadequate. As a consequence, the solution obtained from the stochastic bilevel program may yield poor out-of-sample performance (which is demonstrated in Section \ref{sec:out-of-sample}).
	In this paper, we employ a set of plausible probability distributions, termed as $\mathcal D$, to take into account the ``moment ambiguity''. 
	%
	From a robust perspective, we seek for a solution to the bilevel programs that hedges against all probability distributions belonging to the ambiguity set $\mathcal D$, as follows.
	\begin{eqnarray}
		\label{model:DR-OPTBI}(\mathcal O): \	&&\min_{x\in\mathcal X} w^\top x + \sup_{F\in\mathcal D}\mathbb E_F \left[\min_{y\in \Omega(x,\xi)} v(\xi)^\top y \right],\\
		\label{model:DR-PESBI}(\mathcal P): \		&& \min_{x\in\mathcal X} w^\top x + \sup_{F\in\mathcal D}\mathbb E_F \left[\max_{y\in \Omega(x,\xi)} v(\xi)^\top y \right].
	\end{eqnarray}
	We call  $(\mathcal O)$ and $(\mathcal P)$ the optimistic and pessimistic distributionally robust bilevel programs (DRBPs).
	%
	They minimize the worst-case expected indirect costs of the follower's decision regarding the ambiguity set $\mathcal D$. The two DRBPs can be viewed as special two-stage distributionally robust stochastic programs with $\argmin$ operators in $\Omega(x,\xi)$ which returns the set of optimizers of the follower's utility minimization programs.
	
	\subsection{Assumptions and Ambiguity Sets}
	
	In this paper, we make the following assumptions on the DRBPs. 
	
	\begin{assumption}(Square integrability)
		We assume that $c(\xi),v(\xi)\in\mathcal{L}^2_n(F)$ and $b_x(\xi)\in\mathcal{L}^2_m(F)$ for all $x\in\mathcal{X}$, where $\mathcal{L}^2_r(F)$ denotes the set of $r$-dimensional square-integrable functions of $\xi$ with respect to its probability distribution $F$. That is, all Borel-measurable functions $g$  with $\expec_F[\|g(\xi)\|^2]<\infty$.
	\end{assumption}
	
	\begin{assumption}(Linearity)\label{asp:Linear} We assume that $c(\xi)=C\xi+{c_0}$ for $C\in\realR^{n\times k}$, {$c_0\in\realR^n$}, $v(\xi)=V\xi+{v_0}$ for $V\in\realR^{n\times k}$, {$v_0\in\realR^n$}, $b_x(\xi)=B_x\xi+{b_{x0}}$ for $B_x\in\realR^{m\times k}$, {$b_{x0}\in\realR^m$}, where $B_x=\sum_{i=1}^d B_ix_i+{B_0}$ for $B_i\in\realR^{m\times k},\;i=0,1,\ldots,d$ and {$b_{x0} = \sum_{i=1}^d b_ix_i+b_0$ for $b_i\in\realR^m,\;i=0,1,\ldots,d$}. \end{assumption}
	
	\begin{assumption}(Relatively complete recourse)\label{asp:RECOURSE}
		We assume that $\{y\in\realR^n|Ay\le b_x(\xi)\}$ is non-empty and bounded  almost surely for every $x\in\mathcal{X}$ and for every $\xi$ of any $F\in\mathcal D$.
	\end{assumption}
	We note that the boundedness assumption is equivalent to either statements: (i) $\{y\in\realR^n\ |\ Ay\le 0\}=\{0\}$ or (ii) $\{A^\top p\ |\ p\ge0,\ p\in\realR^m\}=\realR^n$. 
	Suppose that a series of independent data samples $\{\xi^n\}_{n=1}^N$ are drawn from the true probability distribution of $\xi$. We calculate the sample mean $\mu_0$ and covariance matrix $\Sigma_0$ as
	\begin{align*}
		\mu_0:=\frac{1}{N}\sum_{n=1}^N\xi^n\quad \text{and }\quad \Sigma_0:=\frac{1}{N-1}\sum_{n=1}^N(\xi^n-\mu_0)(\xi^n-\mu_0)^\top.
	\end{align*}
	Two moment-based ambiguity sets are considered and defined as follows. 
	\begin{definition}Moment-based ambiguity set  with continuous distributions \citep{delage2010distributionally}. 
		\begin{equation*}\label{def:ambiguity-set}
			\mathcal{D}_M\left(\mathcal{S}, {\mu}_{0}, {\Sigma}_{0}, \gamma_{1}, \gamma_{2}\right):=\left\{F \in \mathcal{M} \left| \begin{array}{l}{\mathbb{P}_F({\xi} \in \mathcal{S})=1} \\ {\left(\mathbb{E}_F[{\xi}]-{\mu}_{0}\right)^{\top} {\Sigma}_{0}^{-1}\left(\mathbb{E}_F[{\xi}]-{\mu}_{0}\right) \leq \gamma_{1}} \\ {\mathbb{E}_F\left[\left({\xi}-{\mu}_{0}\right)\left({\xi}-{\mu}_{0}\right)^{\top}\right] \preceq \gamma_{2} {\Sigma}_{0}}\end{array}\right.\right\},
		\end{equation*}
		where $\gamma_1\ge0$, $\gamma_2\ge 1$, $\mathcal{M}$ is the set of all probability measures, 
		and the support set of $\xi$ is 
		\begin{align*}\label{def:support}
			\mathcal{S}:=\{\xi\in\realR^k\ |\ W\xi\ge h\},
		\end{align*}
		where $W = [w_1,\ldots, w_l]^\top \in \mathbb R^{l \times k}$ and $h = [h_1,\ldots, h_l]^\top \in \mathbb{R}^l$.
	\end{definition}
	
	\begin{assumption}{(Compact support)} The support set $\mathcal S$ is compact.\label{assump:compact}
	\end{assumption}
	\begin{definition} Moment-based ambiguity
		set  with discrete distributions. Suppose that the distribution of uncertainty $\xi$ is supported on a finite set $\mathcal N := \{\xi^1,\ldots, \xi^N\}$ with probability $p_1,\ldots p_N$. The ambiguity set is given by   
		\begin{equation*}\label{def:ambiguity-set-discrete}
			\mathcal{D}_\text{dis}\left(\mathcal{N}, {\mu}_{0}, {\Sigma}_{0}, \gamma_{1}, \gamma_{2}\right):=\left\{p \in \mathbb R_+^{N} \left| \begin{array}{l}{\sum_{s=1}^N p_s = 1} \\ {\left(\sum_{s=1}^N p_s \xi^s-{\mu}_{0}\right)^{\top} {\Sigma}_{0}^{-1}\left(\sum_{s=1}^N p_s \xi^s-{\mu}_{0}\right) \leq \gamma_{1}} \\ {\sum_{s=1}^N p_s \left[\left({\xi^s}-{\mu}_{0}\right)\left({\xi^s}-{\mu}_{0}\right)^{\top}\right] \preceq \gamma_{2} {\Sigma}_{0}}\end{array}\right.\right\},
		\end{equation*}
		where $\gamma_1\ge0$, $\gamma_2\ge 1$.
	\end{definition}
	
	The constraints in the sets $\mathcal D_M$ and $\mathcal D_\text{dis}$ ensure that (i) the true first-order moment $\expec[\xi]$ lies in an ellipsoid centered at $\mu_0$ and (ii) the true second-order moment lies in a cone bounded above by $\gamma_2\Sigma_0$. Parameters $\gamma_1$ and $\gamma_2$ control the distance between the true moments and the moment estimates. Theoretically, \citet{delage2010distributionally}provide guidance on the choice of  $\gamma$'s to guarantee out-of-sample performance. In practice, the parameters can be decided using cross-validation. 
	\subsection{Contributions}
	In this paper, we focus on the pessimistic DRBP where the follower does not act in favor of the leader. The leader therefore  makes decisions against the worst case of the follower's decision. 
	The primary motivation of this work is to take the first step towards modeling and solving distributionally robust pessimistic bilevel programs. 
	Specifically, we show that, under proper conditions of ambiguity sets (not only restricted to the two moment-based ambiguity sets), the pessimistic DRBP is equivalent to a generic two-stage DR  stochastic (nonlinear)  program (TSDR) with random objective and right-hand side.
	For moment-based ambiguity sets of continuous distributions, even a two-stage distributionally robust linear program (TSDRLP), a special case of the resulting TSDR, is NP-hard  with both random objective and right-hand side \citep[e.g., ][]{bertsimas2010models}. Therefore, we approach the problem by linear decision rules (LDRs) to tractably approximate the TSDR.
	%
	When the ambiguity set is restricted to discrete distributions, an exact reformulation is derived for the resulting TSDR. 
	Our main contributions are summarized as follows.
	\begin{enumerate}
		\item We derive an equivalent TSDR reformulation (Section \ref{sec:ts-equivalent}) for the pessimistic DRBP under ambiguity sets satisfying proper conditions. The binary restriction of the leader's decision is not necessarily required for the reformulation. {The derivation is based on finding the optimal value function of the follower's problem and establishing the equivalence between the least favorable outcome  and   the robust outcome with respect to the ambiguity set.}
		\item 
		We consider the resulting TSDR under moment-based ambiguity sets \citep{delage2010distributionally}. 
		\begin{itemize}
			\item  For the  ambiguity set $\mathcal D_M$ of continuous distributions (Section \ref{sec:ldr}), we use LDR techniques which leads to a non-convex program. (1) We  derive a 0-1 semidefinte programming (SDP) approximation via the standard convex duality theory. (2) Alternatively,  leveraging the modern conic programming techniques, we develop an exact 0-1 copositive program (COP) reformulation, which, although intractable, {admits a  0-1 SDP inner approximation}. Both 0-1 SDP approximations provide upper bounds to the DRBP under the ambiguity set of continuous distributions.
			\item  A lower bound to the DRBP with $\mathcal D_M$  can be provided by solving the DRBP under a constructed ambiguity set of discrete distributions.  
			For the ambiguity set $\mathcal D_\text{dis}$ of discrete distributions (Section \ref{sec:discrete}), the resulting TSDR is exactly reformulated as a 0-1 SDP and explicit constructions of the worst-case distributions are provided. 
		\end{itemize}
		\item  The resulting 0-1 SDPs still remain computationally challenging. In Section \ref{sec:benders},	we propose a cutting-plane framework to solve 0-1 SDPs. Moreover, based on the MILP approximation proposed in \citet{yanikoglu2018decision}, another cutting-plane algorithm is proposed for the continuous ambiguity set $\mathcal D_M$.  We computationally evaluate efficiency and effectiveness of the proposed 0-1 SDP approximations in solving a facility location problem under various settings (Section \ref{sec:comp}).
	\end{enumerate}
	%
	
	\section{Related Literature}
	\label{sec:litreview}
	
	\subsection{Bilevel Programs}
	The bilevel program can be interpreted as an optimization problem of the leader who searches for the global minimum. The feasible solutions lie among the optimal solutions of the follower's problem, of which the objective function and (or) the feasible region depend on the leader's decision. In the situation where the follower's problem have more than one optimal solution, the follower returns solutions either in favor (\emph{optimistic version}) of  or adverse (\emph{pessimistic version}) to the  leader's objective. 
	Bilevel program has been widely applied to various problems including revenue management \citep{cote2003bilevel}, supply chain management \citep{ryu2004bilevel}, production planning \citep{iyer1998bilevel},  security \citep{scaparra2008bilevel}, transportation \citep{migdalas1995bilevel}, energy market \citep{carrion2009bilevel}, and many other fields. 
	
	The study of deterministic bilevel programs can be dated back to 1934 from the investigation of market equilibria by Stackelberg  \citep{fudenberg1991game}. However, comprehensive studies began only over the last few decades.
	Most existing work focuses on optimistic bilevel programs. The reason is that, in the optimistic version, the minimization over the follower's and leader's decisions can be combined into a joint minimization in the leader's objective function. Thus, the three interrelated optimizations (i.e., (1) the leader minimizing over his decision, (2) the leader minimizing over the follower's optimal decisions, and (3) the follower optimizing her decisions) can be transformed into a two-level problem.  Moreover, if the follower's problem is convex, the two-level problem  can be further reformulated as a single-level program by replacing the constraints, which enforce the solution being optimal for the follower's problem,  with the corresponding Karush-Kuhn-Tucker optimality conditions \citep[see e.g.,][]{dempe2020bilevel,beck2021gentle}. Global optimization algorithms have been developed to solve the resulting single-level (non-convex) problems \citep[e.g., see,][]{faisca2007parametric,tuy2007novel}. \citet{mitsos2008global} and \citet{tsoukalas2009global} further study the case where the follower's problem is nonconvex. 
	We refer interested readers to \citet{dempe2002foundations} and \citet{dempe2013bilevel} for a comprehensive review to deterministic bilevel program research (up to 2002) and more details on the optimality conditions.

	Due to inherit nature of successive decision making under uncertainty, many applications can be modeled as stochastic bilevel programs. However, stochastic bilevel programs attract less attention than their deterministic counterparts. \citet{ivanov2014bilevel} consider optimistic stochastic bilevel linear programs with right-hand side uncertainty in the follower's problem using quantile criterion. 
	They develop a two-stage stochastic programming equivalence with equilibrium constraints and a mixed-integer linear programming reformulation under a discrete distribution. More coherent or convex risk measures are considered in 
	\citet{burtscheidt2019bilevel} and \citet{burtscheidt2020risk} for  stochastic bilevel linear programs when the right-hand side of the follower's problem is stochastic. They outline Lipschitzian properties, conditions for existence and optimality, and stability results. 
	\citet{zhang2021bilevel} studies stochastic  bilevel programs with follower's integer recourse and stochastic right-hand side using an exact value function-based approach.
	All the work mentioned above studies optimistic stochastic bilevel programs, which can  be viewed as special cases of stochastic mathematical programming with equilibrium constraints (SMPEC) problem \citep[see, e.g.,][]{patriksson1999stochastic,xu2005mpcc, shapiro2006stochastic,shapiro2008stochastic}.
	
	Recently, \citet{yanikoglu2018decision} study pessimistic stochastic bilevel programs with binary leader's decision and continuous follower's decision. They apply decision rules to construct upper and lower bounds on the leader's problem. \citet{tavasliouglu2019solving} extend the problem to mixed-integer follower's decision and develop a generalized value function-based approach. 

	\subsection{Distributionally Robust Optimization}
	In contrast to the previous stochastic bilevel programming work requiring full distributional information of uncertainties, we study the DRO variant, which assumes that only limited distributional information is known when making decisions.
	A set (termed as ambiguity set) of distributions are considered based on the partially known information. Two types of information are widely considered in literature, including the moment information such as mean and covariance \citep[see, e.g.,][]{bertsimas2010models, delage2010distributionally, wiesemann2014distributionally,zhang2018ambiguous}, 
	and the distance to a reference distribution, for example, Wasserstein metrics \citep[see, e.g.,][]{gao2016distributionally,esfahani2018data,zhao2018data}.

	
	In this paper, we  show that the pessimistic DRBPs   admit equivalent TSDR nonlinear reformulations with uncertainty in both  the objective and right-hand side. A special case of the TSDR is 
	TSDRLPs, which have attracted many interests since a decade ago.
	\citet{bertsimas2010models} consider risk averse TSDRLPs when the first- and second-order moments  are known. When the uncertainty only appears in the objective function, they obtain a tight SDP reformulation. When the uncertainty is only in the right-hand side, they show that the problem is NP-hard in general. 
	\citet{xie2018distributionally_int} study TSDRLPs with simple integer round-up recourse using mean and support to form the ambiguity set. They obtain  a mixed-integer second-order cone programming reformulation. Recently, \citet{fan2021decision} consider TSDRLP
	with random recourse matrix using piecewise decision rules based on a partitioning scheme and accordingly a conditional ambiguity set is constructed for marginal probabilities on the partitions. 
	%
	%
	Using Wasserstein metrics, the TSDRLPs admit tractable reformulations under various conditions \citep[see, e.g.,][]{hanasusanto2018conic,xie2019tractable,wang2020second}. For ambiguity sets of distributions supported on a finite set, decomposition methods  have been studied for various settings: \citet{love2015phi} based on $\Phi$-divergence ambiguity set, \citet{bansal2018decomposition,luo2021decomposition} for problems with binary variables or/and conic constraints. Using $L_1$-norm ball ambiguity sets, \citet{jiang2018risk} present a sample average approximation algorithm to solve a TSDRLP with binary first-stage decisions.
	

	\section{Distributionally Robust Two-stage Stochastic Programming Equivalence}
	\label{sec:ts-equivalent}
	
	In this section, we will show that the pessimistic DRBPs can be reformulated to equivalent  TSDRs with a random objective and right-hand side when the ambiguity set is weakly compact.   We first consider a special case when the objective coefficients of the follower's decision  in the leader's problem and in the follower's problem coincide, the DRBP admits a TSDRLP reformulation in Section \ref{sec:speical-case}. Then we focus on the more general case when they are not equal in Section \ref{sec:general-case}. 
	
	\subsection{A Special Case: $c(\xi) = v(\xi)$}
	\label{sec:speical-case}
	
	When the objective coefficients of the leader's and follower's are equal (or when $v(\xi) = Kc(\xi)$, $K$ is a positive constant), the pessimistic version ($\mathcal P$) and the optimistic version ($\mathcal O$) coincide as the following TSDRLP with uncertainties in the objective and the right-hand side.
	\begin{equation}\label{eq:speical-tslp}
		\min_{x\in\mathcal X} w^\top x + \sup_{F\in\mathcal D}\mathbb E_F \left[\min_{Ay \le b_x(\xi)} c(\xi)^\top y \right].
	\end{equation}
	\subsection{General Case: $c(\xi) \neq v(\xi)$}
	\label{sec:general-case}
	Now, we consider the more general case when the follower and the leader concern  different objective coefficients for the pessimistic version ($\mathcal P$). 
	We denote $\Psi(x,\xi) := \max_{y\in \Omega(x,\xi)} v(\xi)^\top y$. Then the pessimistic DRBP $(\mathcal P)$ is rewritten as $\min_{x\in\mathcal X} w^\top x + \sup_{F\in\mathcal D} \mathbb E_F[\Psi(x,\xi)]$.
	%
	%
	%
	%
	Recall that the second-stage feasible region $\Omega(x, \xi) := \argmin_{y}\left\{ c(\xi)^\top y: \ Ay \le b_x(\xi)\right\} \subset \mathbb R^n$ involves  an $\argmin$ operator. 
	The most frequently used solution approach for bilevel problems is to reformulate the bilevel problem into a single-level problem by eliminating the $\argmin$ operator. To eliminate the $\argmin$ operator, for a given distribution $F\in\mathcal D$, we first rewrite the second-stage problem in a square-integrable functional space in Lemma \ref{lem:SECSTAG-PESSBI-CLP}.
	
	%

	\begin{lemma}\label{lem:SECSTAG-PESSBI-CLP}
		For a given leader's solution $x$ and distribution $F\in\mathcal D$, the second-stage  expected cost  of the pessimistic DRBP ($\mathcal P$) with respect to distribution $F$, i.e., $\expec_F[\Psi(x,\xi)]$, is equivalent to the following stochastic program. 
		\begin{subequations}\label{model:SECSTAG-PESSBI-CLP}
			\begin{eqnarray}
				Q_F(x):=\max_{y(\xi)\in\mathcal{L}^2_n(F)}&& {\expec_F[v(\xi)^\top y(\xi)]} \\ 
				\mbox{s.t.} &&y(\xi)\in\underset{y^{\prime}(\xi)\in\mathcal{L}^2_n(F)}{\arg\min}\ \left\{c(\xi)^\top y^{\prime}(\xi):Ay^{\prime}(\xi)\le b_x(\xi)\right\}.
			\end{eqnarray}
		\end{subequations}
	\end{lemma}
	See the detailed proof in Section \ref{sec:proof-to-function_space} of the Appendices. Next, following \cite{yanikoglu2018decision}, an optimal value function $\bar{Q}_F(x)$ of the follower's problem is derived in the following Proposition.
	\begin{proposition}\label{prop:yanikoglu}(Adapted from Proposition 2.1 in \cite{yanikoglu2018decision}) 
		For a given leader's solution $x$ and a distribution $F\in\mathcal D$, the  second-stage expected cost $\mathbb E_F[\Psi(x,\xi)]$ in pessimistic DRBP ($\mathcal P$)	
		is equivalent to
		\begin{subequations}\label{model:SECSTAG-PESSBI-CLP-ELIARG}
			\begin{eqnarray}
				{\underset{y(\xi)\in\mathcal{L}^2_n(F)}{\sup}\ }&& {\expec_F[v(\xi)^\top y(\xi)]} \\ 
				\label{model:SECSTAG-PESSBI-CLP-ELIARG-constr1}\mbox{s.t.}& & {\expec_F[c(\xi)^\top y(\xi)]\le \bar{Q}_F(x)}\\
				&&Ay(\xi)\le b_x(\xi),
			\end{eqnarray}
		\end{subequations}
		where $\bar{Q}_F(x):=\underset{y(\xi)\in\mathcal{L}^2_n(F)}{\min}\ \left\{\expec_F[c(\xi)^\top y(\xi)] :\ Ay(\xi)\le b_x(\xi)\right\}$. 
	\end{proposition}
	In constraint \eqref{model:SECSTAG-PESSBI-CLP-ELIARG-constr1}, the right-hand side $\bar Q_F(x)$ involves a two-stage stochastic program (with only wait-and-see decisions).  Lemma \ref{lem:SECSTAG-PESSBI-CLP-DUAL} incorporates the two-stage stochastic problem into the objective according to the duality of convex stochastic program using risk functions \citep{ruszczynski2006optimization}. The proof is presented in Section \ref{sec:proof-lemma-to-real} of the Appendices. \begin{lemma}\label{lem:SECSTAG-PESSBI-CLP-DUAL}
		For a given leader's solution $x$ and distribution $F\in\mathcal D$, problem \eqref{model:SECSTAG-PESSBI-CLP-ELIARG}, or equivalently $\mathbb E_F[\Psi(x,\xi)]$, is equivalent to the two-stage stochastic program as follows.
		\begin{equation}\label{model:SECSTAG-PESSBI-DUAL}
			\min_{\lambda \ge 0} \mathbb E_F \left[\Phi_\lambda (x,\xi) \right],
		\end{equation}	
		where
		\begin{subequations}
			\label{eq:secstag-pessbi-dual}\begin{eqnarray}\Phi_\lambda (x,\xi) := &\underset{p\in\realR^m,y\in\realR^n}{\min}&b_x(\xi)^\top p+c(\xi)^\top y \\ &{\text { s.t. }} & {A^\top p+\lambda c(\xi)=v(\xi),\ p\ge0}\label{cnstr:SECSTAG-PESSBI-DUAL1}\\
				&{}&{\ Ay\le\lambda b_x(\xi)}. 
				\label{cnstr:SECSTAG-PESSBI-DUAL2}
		\end{eqnarray}\end{subequations}
		%
		%
		%
	\end{lemma}
	\begin{remark}\label{rem:BIL-TWOSTG-FSB}
		The feasibility of problem \eqref{eq:secstag-pessbi-dual} is guaranteed by 
		Assumption \ref{asp:RECOURSE}. Specifically, constraint \eqref{cnstr:SECSTAG-PESSBI-DUAL1} is feasible for any given $\xi$ as $\left\{ A^\top p: \ p\ge 0\right\} = \mathbb R^n$. For constraint \eqref{cnstr:SECSTAG-PESSBI-DUAL2}: if $\lambda>0$, then $\left\{y: \ Ay\le b_x(\xi)\right\} \neq \emptyset$; if $\lambda = 0$, $\left\{y: \ Ay\le 0\right\} = \{0\}$.
		\Halmos
	\end{remark}

	
	Combining the minimization problem in Lemma \ref{lem:SECSTAG-PESSBI-CLP-DUAL} with the outer maximization over the ambiguity set, the worst-case outcome $\sup_{F\in\mathcal D} \min_{\lambda \ge 0} \mathbb E_F[\Phi_\lambda (x,\xi)]$ is upper bounded by the robust outcome $ \min_{\lambda \ge 0} \sup_{F\in\mathcal D} \mathbb E_F[\Phi_\lambda (x,\xi)]$ with respect to the ambiguity set $\mathcal D$. 
	The following theorem further establishes the equivalence between the worst-case outcome and the robust outcome, whose proof is provided in Section \ref{sec:proof-interchangeability} of the Appendices.
\begin{theorem}\label{thm:sup-min-interchangeability}
	Given a leader's solution $x\in\mathcal X$ and  a  convex and weakly compact  ambiguity set $\mathcal D$ of probability measures on $(\mathcal S, \mathcal F)$ with a  compact metric space $\mathcal S$ and its Borel $\sigma$-algebra, the worst-case second-stage expected cost $\sup_{F\in\mathcal D} \mathbb E_F[\Psi(x,\xi)] = $
	\begin{equation}\label{eq:min-max-lambda}
		\min_{\lambda \ge 0} \sup_{F\in\mathcal D} \mathbb E_F[\Phi_\lambda (x,\xi)].
	\end{equation}
\end{theorem}
We are now ready to show that the pessimistic DRBP ($\mathcal P$) is equivalent to a generic TSDR, which is formalized in Corollary \ref{coro:minimax_DRBP}.

\begin{corollary}\label{coro:minimax_DRBP}
	Given  a convex and weakly compact ambiguity set $\mathcal D$ of probability measures on $(\mathcal S, \mathcal F)$ with a  compact metric space $\mathcal S$ and its Borel $\sigma$-algebra, the pessimistic DRBP ($\mathcal P$) under $\mathcal D$ is equivalent to the following TSDR
	\begin{equation}\label{eq:general-tsnlp}
		\min_{x\in\mathcal X, \lambda \ge 0} \left\{ w^\top x + \sup_{F\in\mathcal D} \mathbb E_F \left[\Phi_\lambda (x,\xi) \right] \right\}.
	\end{equation}
	%
	%
\end{corollary}
Note that the two-stage distributionally robust problem \eqref{eq:general-tsnlp} is nonconvex, inherent from the nonconvexity of bilevel program, due to the bilinear term $b_x(\xi)^\top p$ in the objective function of $\Phi_\lambda (x,\xi)$ and $\lambda b_x(\xi)$ in the constraint of $\Phi_\lambda (x,\xi)$, where $b_x(\xi)$ is an affine function of the leader's decision $x$ and $p$ is a recourse decision. 
\begin{remark}
	The binary constraints in the leader's feasible region $\mathcal X$ are not necessarily required
	in Theorem \ref{thm:sup-min-interchangeability} and Corollary \ref{coro:minimax_DRBP}. The results can be easily extended (1) to the case where  linking constraints of the here-and-now leader's decision $x$ and the follower's optimal wait-and-see decision $y$ by properly assuming boundedness and relatively complete recourse for the second-stage problem, and (2) to the case where the objective coefficient of the follower's problem is an affine function of the leader's decision $x$.
	\Halmos
\end{remark}
\begin{remark}
	The moment-based ambiguity set $\mathcal D_M$, which will be discussed in Section \ref{sec:ldr}, satisfies the weak compactness \citep[Proposition 7,][]{sun2016convergence} required by
	Theorem \ref{coro:minimax_DRBP} and Corollary \ref{coro:minimax_DRBP}. 
	Other common choices of ambiguity sets are also compact with respect to the weak topology under  some proper conditions. 
	For example, a $p$-Wasserstein ball is weakly compact when the reference distribution has a finite $p$th moment \citep[Theorem 1,][]{yue2021linear}. For an $f$-divergence (of which examples include Kullback-Leibler divergence, Hellinger divergence, J-divergence, etc.) ball, if the metric space is a compact Polish space, then it is weakly compact  \citep[Lemma 3.2,][]{birghila2021distributionally}. We refer interested readers to \citet{sun2016convergence} for the compactness results of other types of ambiguity sets constructed through moments and mixture distributions. Furthermore,  
	given that the finite support set $\mathcal N$ is compact and $\mathcal D_\text{dis}$ is closed \citep[see Proposition 7][]{sun2016convergence}, by Prokhorov's theorem, the ambiguity set $\mathcal D_\text{dis}$  of discrete distributions (further discussed in Section \ref{sec:discrete})  is as well weakly compact  \citep[see, e.g., ][]{shapiro2002minimax,prokhorov1956convergence}.
	\Halmos
\end{remark}

Corollary \ref{coro:minimax_DRBP} extends the results in \citet{yanikoglu2018decision} from a stochastic bilevel program to a distributionally robust variant. Given that the TSDRLP is a special case of the resulting TSDR \eqref{eq:general-tsnlp} in Theorem \ref{thm:sup-min-interchangeability}, in the rest of the paper, we will focus on the solution approaches to TSDR \eqref{eq:general-tsnlp} with binary leader's decisions under the two types of moment-based ambiguity sets $\mathcal D_M$ and $\mathcal D_\text{dis}$. 
%
%
%

\section{Decision Rule Approximations under Ambiguity Set of Continuous Distributions}
\label{sec:ldr}

In this section, we consider the pessimistic problem ($\mathcal P$) with binary leader's decisions under the moment-based ambiguity set of continuous distributions $\mathcal D = \mathcal D_M$. As pointed out in \citet{bertsimas2010models}, even the TSDRLP \eqref{eq:speical-tslp}  is NP-hard even under ambiguity sets that match the moment information.  We  employ LDRs to conservatively approximate the resulting TSDRs under moment-based ambiguity sets $\mathcal D_M$.
Specifically, we consider the recourse decisions $p$ and $y$ as affine functions of the uncertainty $\xi$: $y(\xi):=Y\xi+{y_0}$ and $p(\xi):=P\xi+{p_0}$, where $P\in\mathbb R^{m\times k},\ p_0\in\mathbb R^m$, $Y\in\mathbb R^{n\times k}$, and $y_0 \in \mathbb R^{n}$.
In this section, we  focus on a conservative approximation of \eqref{eq:general-tsnlp} in the following form:
\small
\begin{equation}\label{eq:tsnlp-ldr}
	\min_{x\in\mathcal X, \lambda \ge 0, P, p_0, Y, y_0} \left\{ w^\top x + \sup_{F\in\mathcal D_M} \mathbb E_F [b_x(\xi)^\top p(\xi) + c(\xi)^\top y(\xi)]: \ A^\top p(\xi) + \lambda c(\xi) = v(\xi), \ p(\xi)\ge 0, \ Ay(\xi)\le \lambda b_x(\xi) \right\}.	
\end{equation}
\normalsize
By deriving the dual of the inner maximization problem \citep[see Lemma 1 of][]{delage2010distributionally}, we obtain an equivalent formulation:
\small
\begin{subequations}\label{model:DR-LDR}
	\begin{eqnarray}
		\underset{Q,q,r,t,\lambda,P,Y,p_0,y_0,x\in\mathcal X}{\min}&&w^\top x + r+t\\
		\label{eq:DR-LDR-constr1}\st&&r\ge \xi^\top \left(B_x^\top P + C^\top Y - Q\right)\xi+  \left(p_0^\top B_x + b_{x0}^\top P +  y_0^\top C + c_0^\top Y -q\right)\xi+
		{b_{x0}^\top p_0}+{c_0^\top y_0},\nonumber \\ 
		&&\forall \xi\in\mathcal S \label{cnstr:DR-LDR-SECSTAG-PESSBI1-SEMIINF1}\\
		\label{eq:DR-LDR-constr3}	&&A (Y\xi+{y_0})\le \lambda (B_x\xi+{b_{x0}}),\ \forall \xi\in\mathcal S\\
		\label{eq:DR-LDR-constr4}	&&A^\top (P\xi+{p_0})+\lambda (C\xi+{c_0})=V\xi+{v_0},\ P\xi+{p_0}\ge0,\ \forall \xi\in\mathcal S \\
		\label{cnstr:DR-PESSBI-POLY-2}		&&t\ge \left(\gamma_{2} \Sigma_{0}+\mu_{0} \mu_{0}^{\top}\right) \cdot Q+\mu_{0}^{\top} q+\sqrt{\gamma_{1}}\left\|\Sigma_{0}^{1 / 2}\left(q+2 Q \mu_{0}\right)\right\|\\
		\label{eq:DR-LDR-sign}		&&Q\succeq 0,\ \lambda\ge0.
	\end{eqnarray}
\end{subequations}
\normalsize
The problem \eqref{model:DR-LDR} is a semi-infinite problem with bilinear terms: $B_x^\top P, p_0^\top B_x, b_{x0}^\top P$ in constraint \eqref{eq:DR-LDR-constr1} and $\lambda B_x, \lambda b_{x0}$ 
in \eqref{eq:DR-LDR-constr3}. As the leader's decision $x$ is binary and  $B_x = \sum_{i=1}^d B_ix_i + B_0, \ b_{x0} = \sum_{i=1}^d b_i x_i + b_0$ for some $B_i\in\mathbb R^{m\times k}, \ b_{i}\in\mathbb R^m, \ i=0,1,\ldots,d$, one can linearize the bilinear terms using McCormick inequalities \citep[e.g.,][]{mccormick1976computability}. Introduce auxiliary variables $\Gamma_i =  B_i^\top P x_i,\ \theta_i = \lambda x_i, \ \omega_i = p_0 x_i,\ \rho_i = P^\top b_i x_i, \ i=1,\ldots,d$ by the following McCormick inequalities:
\small
\begin{eqnarray}
	\label{cnstr:DR-PESSBI-POLY-7}	&&B_i^\top TW-(1-x_i)M \mathbf{1}_{k\times k} \le \Gamma_i\leq B_i^\top TW+(1-x_i)M\mathbf{1}_{k\times k},\ 
	-x_iM\mathbf{1}_{k\times k}\le \Gamma_i\le x_iM\mathbf{1}_{k\times k},\ i=1,\ldots,d\hspace{8mm}\\
	&&	{\lambda-M(1-x_i) \le \theta_i\le Mx_i,\ 0\le\theta_i\le\lambda, \ i=1,\ldots,d}\label{cnstr:DR-PESSBI-POLY-10}\\
	\label{cnstr:omeg1}	&&p_0-(1-x_i)M\mathbf{1}_m\le \omega_i\le p_0+(1-x_i)M\mathbf{1}_m,\ {-x_iM\mathbf{1}_m\le \omega_i \le x_iM\mathbf{1}_m, i=1,\ldots,d}\\
	\label{cnstr:rho2}	\label{cnstr:DR-PESSBI-POLY-77}&&(TW)^\top b_i-(1-x_i)M\mathbf{1}_k \le \rho_i \le (TW)^\top b_i+(1-x_i)M\mathbf{1}_k \  -x_iM\mathbf{1}_k\le \rho_i \le x_iM\mathbf{1}_k,\ i=1,\ldots,d,
\end{eqnarray}
\normalsize
where $M > 0$ is a sufficiently large big-M constant, $\mathbf{1}_r\in\mathbb R^r$ is a $r$-dimensional vector of all ones, and $\mathbf{1}_{k\times k}\in\mathbb R^{k\times k}$ is a $k$-by-$k$ all-ones matrix. 

Recall that	$\mathcal{S}:=\{\xi\in\realR^k\ |\ W\xi\ge h\}$. The  two semi-infinite constraints \eqref{eq:DR-LDR-constr3}-\eqref{eq:DR-LDR-constr4}  can be transformed into finite constraints as follows, by applying standard duality techniques for robust optimization.
\begin{eqnarray}
	&& {AY+\Lambda W=\sum_{i=1}^{d}\theta_iB_i+{\lambda B_0},\ \Lambda h-{Ay_0+\sum_{i=1}^d\theta_i b_i+\lambda b_{0}}\ge0,\ \Lambda\ge0}\label{cnstr:DR-PESSBI-POLY-6}			\\
	&&{A^\top T W+\lambda C=V,\ A^\top p_0+\lambda c_0 = v_0, \ Th+{p_0}\ge0,\ T\ge0}\label{cnstr:DR-PESSBI-POLY-5b}
\end{eqnarray}
%
The remaining semi-infinite constraint \eqref{eq:DR-LDR-constr1}, given any fixed $r,Q,q,Y,y_0,P,p_0$,  involves a nonconvex quadratic function of $\xi$.
In Section \ref{sec:0-1SDP}, we derive a conservative approximation of \eqref{eq:DR-LDR-constr1}  and present a 0-1 SDP approximation of \eqref{eq:tsnlp-ldr}.  In Section \ref{sec:0-1COP}, we utilize the copositve programming scheme to derive an exact reformulation of \eqref{eq:DR-LDR-constr1} and thus provide an exact 0-1 COP reformulation of \eqref{eq:tsnlp-ldr} which, although intractable, also admits a 0-1 SDP approximation. 

\subsection{0-1 Semidefinite Programming Approximation}
\label{sec:0-1SDP}

\begin{theorem}\label{thm:0-1SDP}
	The following 0-1 SDP is a conservative approximation of Problem \eqref{eq:tsnlp-ldr}.
	\small
	\begin{subequations}\label{model:DR-PESSBI-POLY}
		\begin{eqnarray}
			\label{eq:DR-LDR-obj-sdp}		 {\displaystyle \min_{\substack{Q, q, r, t,Y,\lambda,\tau,\theta_i,\Gamma_i,T,\Lambda,\\ \omega_i,\rho_i,y_0,p_0, x\in\mathcal X}}} && w^\top x + r+t \\ 
			\label{eq:DR-LDR-constr1-sdp}		{\st}  && \scalemath{0.9}{\begin{bmatrix}
					\begin{aligned}
						&Q-\Big[\frac{1}{2}\Big(\sum_{i=1}^d\Gamma_i+{B_0^\top TW}+C^\top Y\Big)\\&+\frac{1}{2}\left(\sum_{i=1}^d\Gamma_i+{B_0^\top TW}+C^\top Y\right)^\top\Big]
					\end{aligned} & {\begin{aligned}
							\frac{1}{2}\Big(&q-{\sum_{i=1}^d B_i^\top \omega_i -B_0^\top p_0 -\sum_{i=1}^d \rho_i}\\
							& {-(TW)^\top b_0-C^\top y_0-Y^\top c_0}-W^\top\tau\Big) \end{aligned}}
					\\ 
					{\begin{aligned}
							\frac{1}{2}\Big(&q-{\sum_{i=1}^d B_i^\top \omega_i -B_0^\top p_0 -\sum_{i=1}^d \rho_i}\\
							& {-(TW)^\top b_0-C^\top y_0-Y^\top c_0}-W^\top\tau\Big)^\top \end{aligned}} & 
					\quad r-{\sum_{i=1}^d b_i^\top \omega_i-b_0^\top p_0-c_0^\top y_0}+\tau^\top h \end{bmatrix} \succeq0 
			}\hspace{10mm}			\label{cnstr:DR-PESSBI-POLY-1} \\
			\label{eq:DR-LDR-constr1-sdp2}		& & { \tau\ge0}\label{cnstr:DR-PESSBI-POLY-3}\\
			&& \nonumber \eqref{cnstr:DR-PESSBI-POLY-2}-\eqref{eq:DR-LDR-sign}, \ \eqref{cnstr:DR-PESSBI-POLY-7}-\eqref{cnstr:DR-PESSBI-POLY-5b}.		
		\end{eqnarray}
	\end{subequations} 
	\normalsize
\end{theorem}
Section \ref{sec:proof_thm_sdp} in the Appendices presents the detailed proof.
\subsection{0-1 Copositive Programming Exact Reformulation}
\label{sec:0-1COP}

Constraint \eqref{eq:DR-LDR-constr1} is equivalent to the following constraint involving a nonconvex quadratic maximization problem:
\begin{equation}\label{eq:DR-LDR-constr1-max}
	r \ge \max_{\xi\in\mathcal S}\  \xi^\top \left(B_x^\top P + C^\top Y - Q\right)\xi+  \left(p_0^\top B_x + b_{x0}^\top P +  y_0^\top C + c_0^\top Y -q\right)\xi+
	{b_{x0}^\top p_0}+{c_0^\top y_0}.
\end{equation}
In this section, we reformulate the constraint above using the generalized copositive programming techniques \citep[e.g.,][]{burer2012representing}. First, we review the definitions of \emph{copositive cone} and \emph{completely positive cone}. 

\begin{definition}
	Let $\mathcal K\subset \mathbb R^n$ be a closed convex cone. The \textit{copositive cone} with respect to $\mathcal{K}$ is 
	\begin{align*}
		\mathcal{COP}(\mathcal{K}):=\bigg\{M\in\mathbb S^n:\ z^\top Mz\geq0\ \forall z\in\mathcal{K}\bigg\},
	\end{align*}
	and its dual cone, the \textit{completely positive cone}, with respect to $\mathcal K$, is 
	\begin{align*}
		\mathcal{CP}(\mathcal{K}): = \bigg\{X\in\mathbb S^n:\ X=\sum_{i}
		z^i(z^i)^\top,\  z^i\in\mathcal{K}\bigg\}, 
	\end{align*}
	where the summation is finite and the cardinality is unspecified.
\end{definition}  
We define $z:=\begin{bsmallmatrix} \xi\\ \\ 1  \end{bsmallmatrix}$. The polyhedral support $\mathcal S$ of $\xi$ then induces a polyhedral cone of $z$,  $\hat{\Xi} := \left\{z\in\mathbb R^{k+1}: \ \mathcal{H}z\geq0\right\}$, where  $\mathcal{H}=\left[W\ |\ -h \right]\in\realR^{l\times(k+1)}$. In addition to the compactness assumption \ref{assump:compact} on the support $\mathcal S$, the polyhedral cone is assumed to satisfy the following  condition.
\begin{assumption}\label{assump:proper}
	$\hat{\Xi}$ is full dimensional. That is, there is no implicit equalities in $\hat{\Xi}$.
\end{assumption}
Under Assumption \ref{assump:proper}, the polyhedral cone $\hat\Xi$ is proper \citep[see Lemma 2 of][]{mittal2021finding}. 
Now, consider
$$\mathcal{L}:=\left\{\begin{bmatrix} 1 \\ z\end{bmatrix} \begin{bmatrix} 1 \\ z\end{bmatrix}^\top\in\mathbb S^{k+2}:\ e_{k+1}^\top z=1,\ z\in\hat \Xi\right\},$$
where $e_{k+1}\in \realR^{k+1}$ denotes the vector with one in the $k+1$th coordinate and zeros elsewhere.
Then constraint \eqref{eq:DR-LDR-constr1-max} is rewritten as $r\geq v^*$, where 
\begin{equation}
	v^* : = \max\left\{ \tr( {\mathcal{Q}} Z): \ \begin{bmatrix}
		1 & z^\top \\ z & Z
	\end{bmatrix} \in \text{clconv}\left(\mathcal L\right) \right\}.
	\label{eq:RHS_r}
\end{equation}
	with 
	$$\mathcal{Q}:=
	\begin{bmatrix}
		\begin{aligned}
			&\Big[\frac{1}{2}\Big(\sum_{i=1}^d\Gamma_i+{B_0^\top TW}+C^\top Y\Big)\\&+\frac{1}{2}\left(\sum_{i=1}^d\Gamma_i+{B_0^\top TW}+C^\top Y\right)^\top\Big]-Q
		\end{aligned}
		& &   {\begin{aligned}
				\frac{1}{2}\Big(&{\sum_{i=1}^d B_i^\top \omega_i +B_0^\top p_0 +\sum_{i=1}^d \rho_i}\\
				& {+(TW)^\top b_0+C^\top y_0+Y^\top c_0}-q\Big) \end{aligned}}\\
		{\begin{aligned}
				\frac{1}{2}\Big(&{\sum_{i=1}^d B_i^\top \omega_i +B_0^\top p_0 +\sum_{i=1}^d \rho_i}\\
				& {+(TW)^\top b_0+C^\top y_0+Y^\top c_0}-q\Big)^\top \end{aligned}} & & b_{x0}^\top p_0+c_0^\top y_0
	\end{bmatrix}\in\mathbb{S}^{k+1}.$$ 
	
	\begin{proposition}[Adapted from Theorem 8.1 of \citet{burer2012copositive}]  $\text{clconv}(\mathcal L) = \mathcal R$, where $$\mathcal R := \left\{\begin{bmatrix} 1 & z^\top \\ z & Z \end{bmatrix}\in\mathcal{CP}(\realR_+\times\hat \Xi):\ \begin{matrix} e_{k+1}^\top z=1,\\ e_{k+1}^\top Ze_{k+1}=1 \end{matrix}\right\}.$$
	\end{proposition}
	Then the optimal value $v^*$ of \eqref{eq:RHS_r} equals to that of
	the following completely positive program (CPP).
	\begin{equation}
		\max \left\{ \tr({\mathcal{Q}}Z): \ \begin{bmatrix}
			1 & z^\top \\ z & Z
		\end{bmatrix} \in \mathcal R\right\}.
		\label{model:CPP1}
	\end{equation}
		It is possible to further eliminate the last variable $z_{k+1}$.
		\begin{proposition}\label{prop:extraneous} The CPP \eqref{model:CPP1} is equivalent to
			\begin{equation}\label{model:CCP2}
				v^* =  \max \left\{  \tr\left(\mathcal{Q} Z\right): \
				Z\in\mathcal{CP}(\hat \Xi),\ e_{k+1}^\top Ze_{k+1}=1 \right\}.\end{equation}
		\end{proposition}
		The proof is provided in Section \ref{sec:proof-prop-cop} of the Appendices.
		Besides a more compact feasible region, an additional benefit of \eqref{model:CCP2} over \eqref{model:CPP1} is that \eqref{model:CCP2} has  nonempty interior due to Assumption \ref{assump:proper} under which the polyhedral cone $\hat{\Xi}$ is a proper cone \citep[see Lemma 4 of][]{mittal2021finding}. 
		Thus, strong duality holds between \eqref{model:CCP2} and its copositive programming (COP) dual problem:
		$$\min\left\{u: \ ue_{k+1}e_{k+1}^\top-\mathcal{Q}\in\mathcal{COP}(\hat \Xi) \right\}.$$
		We are now ready to present the 0-1 COP reformulation.
		\begin{theorem}\label{thm:0-1COP}
			The following 0-1 COP is an exact reformulation of Problem \eqref{eq:tsnlp-ldr}.
			\begin{subequations}\label{model:COP_reform}
				\begin{eqnarray}
					{\displaystyle \min_{\substack{Q, q, r, t,Y,\lambda,\theta_i,\Gamma_i,T,\Lambda,\\ \omega_i,\rho_i,y_0,p_0, x\in\mathcal X}}} && w^\top x + r+t \\ 
					{\st}  && re_{k+1}e_{k+1}^\top-\mathcal{Q}\in\mathcal{COP}(\hat \Xi)\label{constr:COP1}\\
					&& \nonumber \eqref{cnstr:DR-PESSBI-POLY-2}-\eqref{eq:DR-LDR-sign}, \ \eqref{cnstr:DR-PESSBI-POLY-7}-\eqref{cnstr:DR-PESSBI-POLY-5b}.	
				\end{eqnarray}
			\end{subequations} 
		\end{theorem}
		
		The detailed proof is in Section \ref{sec:proof-thm-cop} of the Appendices.	The resulting generalized copositve program is, however, in general intractable \citep{burer2012copositive} and one alternative is to replace the copositive cone $\mathcal {COP}(\hat\Xi)$ with a semidefinite-based inner approximation \citep{xu2018data}: $\textrm{IA}(\hat\Xi) = \left\{\mathcal{\mathcal{H}}^\top U\mathcal{H}\ :\ U\geq 0,\ U\in\mathbb S^l\right\}\subset \mathcal{COP}(\hat{\Xi})$.
		\begin{corollary}
			\label{thm:cop-ia}
			A conservative 0-1 SDP approximation of \eqref{eq:tsnlp-ldr} is given by
			\begin{subequations}\label{model:DR-PESSBI-POLY-IA_COP}
				\begin{eqnarray}
					{\displaystyle \min_{\substack{Q, q, r, t,Y,\lambda,\theta_i,\Gamma_i,T,\Lambda,\\ \omega_i,\rho_i,y_0,p_0, x\in\mathcal X}}} && w^\top x + r+t \\ 
					{\st}  && re_{k+1}e_{k+1}^\top-\mathcal{Q}= \mathcal H^\top U \mathcal H, \ U\ge 0, \ U\in\mathbb S^l
					\label{constr:COP1-ia}\\
					&& \nonumber \eqref{cnstr:DR-PESSBI-POLY-2}-\eqref{eq:DR-LDR-sign}, \ \eqref{cnstr:DR-PESSBI-POLY-7}-\eqref{cnstr:DR-PESSBI-POLY-5b}.	
				\end{eqnarray}
			\end{subequations} 
		\end{corollary}

		\section{Exact Reformulation under Ambiguity Set of Discrete Distributions}
		\label{sec:discrete}
		Under an ambiguity set of discrete distributions $\mathcal D =\mathcal D_\text{dis}$, the pessimistic bilevel program ($\mathcal P$)  admits exact 0-1 SDP reformulation, which is formally stated in Theorem \ref{thm:discrete-reformulation} and constitutes a simple corollary of Theorem \ref{thm:0-1SDP}. Thus, the proof is omitted for brevity.
		\begin{theorem}\label{thm:discrete-reformulation} Under the ambiguity set of discrete distributions $\mathcal D_\text{dis}$, the pessimistic bilevel program ($\mathcal P$) is equivalent to the 0-1 SDP:
			\begin{subequations}\label{eq:sdp-discrete}
				\begin{eqnarray}
					&\underset{x\in\mathcal{X},Q,q,r,t,\lambda,p^s,y^s}{\min}&w^\top x +r+t\\
					&\st&r\ge \sum_{i=1}^d(B_i\xi^s+b_i)^\top\omega_i^s+(B_0\xi^s+b_0)^\top p^s+c(\xi^s)^\top y^s-{\xi^s}^\top Q\xi^s- {\xi^s}^\top q,\nonumber\\
					&&A^\top p^s+\lambda c(\xi^s)=v(\xi^s),\ A y^s\le \sum_{i=1}^d \theta_i(B_i\xi^s+b_i)+\lambda(B_0\xi^s+b_0),\ s=1,\ldots,N\hspace{10mm}\\
					\label{eq:sdp-discrete-t}	&&t\ge\left(\gamma_{2} \Sigma_{0}+\mu_{0} \mu_{0}^{\top}\right) \cdot Q+\mu_{0}^{\top} q+\sqrt{\gamma_{1}}\left\|\Sigma_{0}^{1 / 2}\left(q+2 Q \mu_{0}\right)\right\|\\
					&&\omega_i^s\le p^s+(1-x_i)M\mathbf{1}_m,\ \omega_i^s\ge p^s-(1-x_i)M\mathbf{1}_m,\ i=1,\ldots,d,\ s=1,\ldots,N\\
					&&-x_iM\mathbf{1}_m\le \omega^s_i \le x_iM\mathbf{1}_m,\ i=1,\ldots,d,\ s=1,\cdots,N\\
					&&{\theta_i\le Mx_i,\ \theta_i\ge\lambda-M(1-x_i),\ 0\le\theta_i\le\lambda,\ i=1,\ldots,d}\\
					\label{eq:sdp-discrete-sign}	&&Q\succeq 0,\ \lambda\ge0,\ p^s\geq0,\ s=1,\ldots,N,
				\end{eqnarray}
			\end{subequations}
			where $M>0$ is a sufficiently large big-M constant.
		\end{theorem}
		
		Given a solution of the leader's decision $\hat x$, the worst-case distribution can be derived  based on the dual problem of the following SDP equivalent of the second-stage problem.
		\begin{subequations}\label{eq:discrete-sub-primal}
			\begin{eqnarray}
				\min_{Q,q,r,t,\lambda,p^s,y^s} && r+t\\
				\mbox{s.t.}
				\label{eq:discrete-sub-constr1}		&&r\ge b_{\hat x}(\xi^s)^\top p^s+c(\xi^s)^\top y^s-{\xi^s}^\top Q\xi^s- {\xi^s}^\top q,\ \forall s=1,\ldots,N\\
				&&A^\top p^s+\lambda c(\xi^s)=v(\xi^s),\ p^s\ge0,\ A y^s\le \lambda b_{\hat x}(\xi^s),\ \forall s=1,\ldots,N\\ &&\eqref{eq:sdp-discrete-t}, \ \eqref{eq:sdp-discrete-sign}.\nonumber
			\end{eqnarray}
		\end{subequations}

		\begin{theorem}\label{thm:discrete-worst-case}
			Given a leader's decision $\hat x$, assume the dual optimal solution associated with constraints \eqref{eq:discrete-sub-constr1},  denoted as $\gamma^*\in\mathbb R^N$. The worst-case distribution is characterized as $\mathbb P\{\xi=\xi^s\} = \gamma^{s*}, \ s=1,\ldots,N$.
		\end{theorem}
		The proof is given in  Section \ref{sec:discrete-distribution} of the Appendices.
		For any given continuous ambiguity set $\mathcal D_M$, one can generate samples in the support set of $\mathcal D_M$ and construct a discrete ambiguity set as an approximation. Naturally, the DRBP model under the discrete ambiguity set provides a lower bound to that under the continuous ambiguity set. Given that the two 0-1 SDP formulations proposed in Section \ref{sec:ldr} provide upper bounds to the DRBP, together with the lower bounds using discrete ambiguity sets, the approximate gaps of the LDR method can be quantified. Computational studies on the gaps are presented in Section \ref{sec:time_comp}. 
		
		\section{Solution Algorithms}
		\label{sec:benders}
		
		In this section, we first propose a cutting-plan framework to efficiently solve the proposed 0-1 SDP formulations. Then another cutting-plane procedure is developed based on the MILP upper-bounding approximation proposed in \citet{yanikoglu2018decision}.
		
		\subsection{Cutting-Plane Framework for the Three 0-1 SDP}
		The three 0-1 SDP formulations: \eqref{model:DR-PESSBI-POLY} for conservatively approximating the pessimistic problem $(\mathcal P)$ under the ambiguity set of continuous distributions, \eqref{model:DR-PESSBI-POLY-IA_COP} for an inner approximation of the exact 0-1 COP reformulation under the ambiguity set of continuous distributions, and \eqref{eq:sdp-discrete} for $(\mathcal P)$ under the ambiguity set of discrete distributions, are computationally challenging 
		as the problem size increases. 
		In this section, we develop a Benders-type cutting-plane algorithm to solve the three 0-1 SDP problems.
		The cutting-plane algorithm iteratively solves a relaxed master problem and a subproblem to generate optimality cuts (due to the relatively complete recourse assumption \ref{asp:RECOURSE},  no feasibility cuts are needed.) We define variable $\nu$ as an underestimator of the worst-case second-stage expected cost $\sup_{F\in\mathcal D}\mathbb E_F\left[\Psi(x,\xi)\right]$. 
		The relaxed master problem  is given by 
		\begin{subequations}\label{eq:RLX-MASTERPB}
			\begin{eqnarray}
				\textbf{MP}:\quad &\underset{x\in\mathcal{X},\nu}{\min}& w^\top x+\nu  \\
				&\st&\nu\ge u_l^\top x+a_l,\ l=1,\ldots,L, 
			\end{eqnarray}
		\end{subequations}
		where the parameter $u_l\in\realR^d$ is the coefficient of $x$ and $a_l\in\realR$ is the scalar parameter of the valid inequality generated from the $l$-th subproblem. After obtaining the optimal solution $(\hat x, \hat \nu)$, we solve a subproblem. If the optimal value of the subproblem is greater than the underestimator $\hat\nu$, an optimality cut is generated and added to the  MP. The details of the algorithm are presented in Algorithm \ref{alg:BENDER-LAG}. The algorithm guarantees finite termination for any given $\epsilon \ge 0$  \citep{geoffrion1972generalized}.
		%
		%
		\begin{algorithm}[h]
			\caption{A cutting-plane framework for solving 0-1 SDP formulations (\eqref{model:DR-PESSBI-POLY}, \eqref{model:DR-PESSBI-POLY-IA_COP}, and \eqref{eq:sdp-discrete})}
			\label{alg:BENDER-LAG}
			\begin{algorithmic}[1]
				\State Set $LB\leftarrow -\infty,\ UB\leftarrow +\infty$ and $\epsilon > 0$.
				\For {$\ell=0,1,2,\ldots$}
				\State Solve the relaxed master problem MP. If the problem is infeasible, claim the infeasibility and quit the loop. Otherwise let $(x^\ell,\nu^\ell)$ be the optimal solution and $z^\ell$ be the optimal value. (When $\ell=0$, let $\nu^{\ell}=-\infty$, and $x^\ell=\underset{x\in\mathcal{X}}{\arg\min}\ \left\{w^\top x\right\}$.) 
				\State Update the lower bound: $LB\leftarrow z^\ell$.
				\State Solve the subproblem $\text{SP}(x^{\ell})$. Obtain the optimal solution and optimal value $obj^\ell$.\label{algstep:SP} 
				\If {$w^\top x^{\ell}+obj^\ell<UB$}
				\State Let $x^*\leftarrow x^\ell$ be the incumbent solution and update the upper bound: $UB\leftarrow w^\top x^{\ell}+obj^\ell$.
				\EndIf
				\If{$UB-LB>\varepsilon$}
				\State Add an optimality cut $\nu\ge u_l^\top x + a_l$.\label{algstep:OPTCUT}
				\Else 
				\State \Return $x^*$ as the optimal solution to 0-1 SDP problem.
				\EndIf
				\EndFor
			\end{algorithmic}
		\end{algorithm}
		The subproblem solved in 
		Line \ref{algstep:SP} and the optimality cut in Line \eqref{algstep:OPTCUT} are specified in Section \ref{sec:benders-app} of the Appendices for the three  0-1 SDPs, respectively.

		
		%
		\subsection{MILP-based Cutting-Plane Algorithm}
		\label{sec:dcg}
		When the two moments $(\mu, \Omega)$ of the worst-case distribution are known, the second-stage cost of the conservative decision-rule based approximation \eqref{eq:tsnlp-ldr} is equivalent to an MILP \citep{yanikoglu2018decision}.
		\small
		\begin{subequations}
			\begin{eqnarray}
				\min_{\substack{Y,\lambda,\theta_i,\Gamma_i,T,\Lambda,\\ \omega_i,\rho_i,y_0,p_0, }} && \tr{\left[\Omega\left(\sum_{i=1}^d \Gamma_i + B_0^\top TW+ C^\top Y \right)\right]} + \mu^\top \left(\sum_{i=1}^d B_i^\top \omega_i + B_0^\top p_0 + \sum_{i=1}^d \rho_i + (TW)^\top b_0 + C^\top y_0 + Y^\top c_0\right) \nonumber\\
				&& + \sum_{i=1}^d b_i^\top \omega_i + b_0^\top p_0 + c_0^\top y_0\hspace{9mm}\nonumber\\
				\mbox{s.t.}	&& \eqref{cnstr:DR-PESSBI-POLY-7}-\eqref{cnstr:DR-PESSBI-POLY-5b}. \nonumber 
			\end{eqnarray}
		\end{subequations}
		\normalsize
		It is easy to see that the decision-rule based approximation \eqref{eq:tsnlp-ldr} is upper bounded by 
		\small
		\begin{subequations}
			\label{eq:proj_inf}		\begin{eqnarray}
				\min_{\substack{v,Y,\lambda,\theta_i,\Gamma_i,T,\Lambda,\\ \omega_i,\rho_i,y_0,p_0, x\in\mathcal X}} &&w^\top x + v \\
				\mbox{s.t.}			&& v\ge  \tr{\left[\Omega\left(\sum_{i=1}^d \Gamma_i + B_0^\top TW+ C^\top Y \right)\right]} + \mu^\top \left(\sum_{i=1}^d B_i^\top \omega_i + B_0^\top p_0 + \sum_{i=1}^d \rho_i + (TW)^\top b_0 + C^\top y_0 + Y^\top c_0\right) \nonumber \\ 
				\label{eq:proj_inf_constr}			&&+ \sum_{i=1}^d b_i^\top \omega_i + b_0^\top p_0 + c_0^\top y_0, \ \forall (\mu,\Omega)\in \mathcal D_\text{moment} \\
				&&  \lambda \ge 0, \ \eqref{cnstr:DR-PESSBI-POLY-7}-\eqref{cnstr:DR-PESSBI-POLY-5b},\nonumber 
			\end{eqnarray}
		\end{subequations}
		\normalsize
		where $\mathcal D_\text{moment} = \left\{(\mu,\Omega): \ (\mu-\mu_0)^\top \Sigma_0^{-1}(\mu - \mu_0) \le \gamma_1, \ \Omega - \mu\mu_0^\top - \mu_0\mu^\top + \mu_0\mu_0^\top  \preceq \gamma_2 \Sigma_0, \ \mu\mu^\top \preceq \Omega \right\}$, which contains all the moments of any distributions in the ambiguity set $\mathcal D_M$. Although \eqref{eq:proj_inf_constr} incorporates an infinite number of constraints, we can relax constraints \eqref{eq:proj_inf_constr} and iteratively add them back if needed. Specifically, in each iteration, we obtain an incumbent solution $(\hat x, \hat v, \hat \Gamma_i, \hat T,\hat Y, \hat y_0,\hat \rho_i,\hat \omega_i, \hat p_0)$ from the relaxed formulation. Then, we solve the following separation problem to decide if any solution violates constraints \eqref{eq:proj_inf_constr}.
		\footnotesize
			\begin{eqnarray}\label{eq:proj_separation}
				\max_{(\mu,\Omega)\in\mathcal D_\text{moment}}&& \tr{\left[\Omega\left(\sum_{i=1}^d \hat\Gamma_i + B_0^\top \hat TW+ C^\top \hat Y \right)\right]} + \mu^\top \left(\sum_{i=1}^d B_i^\top \hat \omega_i + B_0^\top p_0 + \sum_{i=1}^d \hat \rho_i + (\hat TW)^\top b_0 + C^\top \hat y_0 + \hat Y^\top c_0\right)\nonumber\\
				&& + \sum_{i=1}^d b_i^\top \hat \omega_i + b_0^\top \hat p_0 + c_0^\top \hat y_0
			\end{eqnarray}
		\normalsize
		If not, then $(\hat x, \hat v, \hat \Gamma_i, \hat T,\hat Y, \hat y_0,\hat \rho_i,\hat \omega_i, \hat p_0)$ is an optimal solution to \eqref{eq:proj_inf}; otherwise, we obtain the optimal solution $(\hat\mu, \hat\Omega)$ to \eqref{eq:proj_separation} which violates constraints \eqref{eq:proj_inf_constr}. We add this violated constraint back into the relaxed formulation of \eqref{eq:proj_inf} to cut off the incumbent solution. Since $\mathcal X$ is finite, the iterative cutting-plane algorithm  terminates within a finite number of iterations. Note that the proposed cuts can be added in	a delayed constraint generation fashion in a branch-and-cut framework.

		\section{Computational Studies}
		\label{sec:comp}
		
		We evaluate the performance of the proposed 0-1 SDP approximations for the pessimistic DRBP $(\mathcal P)$ on a facility location problem following \citet{yanikoglu2018decision}.
		The facility location problem is described in Section \ref{sec:fl_prob}. 
		The experimental setup is described in Section \ref{sec:exp_setup1}. 
		The comparison of CPU time and approximate optimality gap of solving the pessimistic DRBP using different formulations are presented in Section \ref{sec:time_comp}. We then consider two cases : (i) when the uncertainty is only in the constraints ($C=V=0$), and (ii) when the uncertainty is in both the 
		constraints and objective ($C,V\neq 0$). The solution details of the two cases are presented in Section \ref{sec:sol-details} and the out-of-sample performance in Section \ref{sec:out-of-sample}.  {All the computational tests are performed on a 64-bit Windows 10 machine with Intel Core i7-4770 CPU 3.40 GHz and 16 GB memory. All the models have been implemented using YALMIP in MATLAB.} Specifically, all the  0-1 SDP formulations are solved using the cutting-plane algorithm proposed in Section \ref{sec:benders}, of which the master problems are solved using Gurobi 9.1.2 and the subproblems are solved using MOSEK 9.3.

		\subsection{Facility Location Problem of a Market Entrant} \label{sec:fl_prob}
		Consider two market companies A and B selling homogeneous products to $d$ demand locations. Among the $d$ locations, a subset of them denoted by vector $l^S\in\{0,1\}^d$ are eligible for accommodating retail stores where $l^S_i = 1$ if a retail store can be accommodated at location $i$, 0 otherwise. Company A already operates retail stores at a subset of these eligible locations. We define $l^A\in\{0,1\}^d$ for Company A where $l^A_i = 1$ if Company A owns a store at location $i$, 0 otherwise. Company B wants to enter the market to build at most $N_B$ number of new stores where Company A does not have stores yet. As a leader in the bilevel program, Company B decides where to build stores denoted by a binary vector $x\in\{0,1\}^d$ where $x_i = 1$ if Company B opens a store at location $i$, 0 otherwise. We assume that at each eligible location, at most one new store can be built either by Company A or by Company B. The feasible region of $x$ is given by $\mathcal{X}=\big\{x\in\{0,1\}^d\;|\;x+l^A\leq l^S,\;\|x\|_1\leq N_B\big\}$. 
		
		Let $\xi_i$ be the uncertain demand at location $i$ and $\xi = (\xi_i,\ i=1,\dots,d)^\top$ be the vector of the demand. We assume that the customers select the nearest store for purchase regardless of the owner. An aggregate customer is considered as the follower whose decision $y = (y_{ij}, \ i,j=1,\ldots,d)^\top\in\mathbb R^{d^2}$ denotes the amount of the products which are supplied from location $i$ to location $j$. Let $b_i$ be the store capacity at location $i$. Given the leader's decision $x$ and demand $\xi$, the follower solves a transportation problem as follows.
		\vspace{-0.05cm}
		\begin{subequations}\label{eq:follower}
			\begin{eqnarray}
				&\underset{y}{\min} & \sum_{i=1}^d\sum_{j=1}^d c_{ij}y_{ij} \label{follwr1}\\
				&\st & \sum_{i=1}^d y_{ij} \ge \xi_j, \ j=1,\ldots,d \label{follwr2}\\
				&& \sum_{j=1}^d y_{ij} \le b_i(x_i + l^A_i), \ i=1,\ldots,d \label{follwr3}\\
				&& y_{ij}\ge 0, \ i,j=1,\ldots,d \label{follwr4}.
			\end{eqnarray}
		\end{subequations} 
		The objective \eqref{follwr1} minimizes the total shipping cost where $c_{ij} $ is the unit transportation cost from location $i$ to location $j$. Constraints \eqref{follwr2} ensure that all the demands are satisfied and constraints \eqref{follwr3} enforce the store capacity if there is a store at location $i$. Constraints \eqref{follwr4} ensure the nonnegativity. Let $Y^*(x,\xi)$ be the set of optimal solutions of the follower's problem \eqref{eq:follower}. We assume that $\max_{\xi \in \Xi} \|\xi\|_1\le \sum_{i=1}^d b_i l^A_i$ and thus the follower's problem is always feasible (i.e., $Y^*(x,\xi) \neq \emptyset$) for every $x$ and $\xi$. 
		
		For the leader (i.e., Company B), let $w_i$ be the cost of opening a store at location $i$ and $v_{ij}$ be the negative of revenue per unit product sold to location $j$ from location $i$. Denote vector $w = (w_i, \ i=1,\ldots, d)^\top$ and vector $v = (v_{ij}, \ i,j = 1,\ldots,d)^\top$. The pessimistic DRBP is formulated as follows.
		\vspace{-0.05cm}
		\begin{align}
			\begin{array}{cl}{\underset{x\in\mathcal{X}}{\inf}\ }& {\left\{w^\top x+\underset{F\in\mathcal{D}}{\sup}\ \expec_F\left[ \underset{y \in\realR^n}{\sup}v(\xi)^\top y\right] : \ {y\in Y^*(x,\xi)  }\right\}}
			\end{array}
		\end{align}
		The model decides how many stores to open and where to build them for Company B under the uncertain demand and the aggregate customer's optimal decisions. 
		\subsection{Computational Setup}
		\label{sec:exp_setup1}
		
		We consider first $d=8$ locations of the SGB128 dataset\footnote{The dataset is available from John Burkardt's website: \url{https://people.sc.fsu.edu/~jburkardt/datasets/cities/cities.html}}  of North America cities, among which $\|l^S\|_1 = 5$ are eligible for building a retail store, where $l^S_i = 1$ for $i=1,2,3,4,6$, and $l^S_i = 0$ for other $i$'s.  Company A already operates at location $i=6$, i.e., $\|l^A\|_1 = 1$, where $l^A_6 = 1$ with store capacity $b_6 = 240 \cdot (d-\| l^S \|_1 ) / \| l^A \|_1 = 720$\footnote{The capacity is chosen such that the follower's problem is always feasible irrespective of whether Company B opens its stores or not for the relatively complete recourse assumption.}. Company B can build at most $N_B = 4$ stores among the remaining eligible locations. {The store capacities at the eligible locations are  $b_i = 360$ for $i=1,2,3,4$. } The fixed cost of opening a new store is $w_i = 305$ and the fixed revenue of selling one product is $\$ 5$ at eligible locations for company B, i.e., $v_{ij} = -5$ if $l_i^S = 1 \text{ and } l_i^A = 0$, otherwise $0$. The fixed unit transportation cost $c_{ij}$ is set to be the Euclidean distance between location $i$ and $j$. We generate $N$ realizations of the demands $\xi_1^n,\ldots,\xi_d^n, \ n=1,\ldots, N$ sampled using the following procedure. We assume that all the demands are independent, identical and generate $\xi_i$, at each demand location $i$ with $l_i^S = 0$ (not eligible to build stores), following the uniform distribution on $[30,240]$. There is no demand, $\xi_i = 0$, at eligible location $i$ where $l_i^S = 1$. We solve the 0-1 SDP approximations  using  the empirical mean $\mu_0$ and covariance $\Sigma_0$ calculated from $N= 10$ samples and support set $\mathcal S = \{ 30 \le \xi_i \le 240,\ \text{for } i \text{ such that } l_i^S = 0\}$. All big-M constants are set to $10^6$. Given the optimal solutions obtained by solving different models, we generate $N^\prime = 5000$ data samples for out-of-sample test with details given in Section \ref{sec:out-of-sample}.
		
		\subsection{CPU Time and Computational Details} \label{sec:time_comp}
		To assess the accuracy and runtime of the proposed approximations, we randomly generate instances of the facility location problem with $d\in\{15,20,25\}$ locations, $\|l^S\|_1 \in\{5,10\}$ eligible locations for building stores, $\|l^A\|_1\in\{2,3,5\}$ locations occupied by Company A. The coordinates of the $d$ locations are chosen independently and uniformly from $[0,1]^2$. The random demand  $\xi_i$ at location $i$ follows the uniform distribution on $[50, 150]$.
		{The store capacity is $b_i=150 \cdot (d-\| l^S \|_1 ) / \| l^A \|_1$ for $i$ such that $l_i^S=1$. }
		There are $N_B \in \{2,3,5\}$ candidate locations for Company B to establish new stores at a fixed cost of $w_i\in\{750,1000\}$ if the uncertainty only appears in the constraints ($C=V=0$). When the uncertainty is in both the objective and the constraints ($C,V\neq 0$), we consider a double  cost $w_i\in\{1500,2000\}$  {so that the ratio of the fixed store opening cost ($w_i$) and the product selling price ($v_{ij}$) remains compatible}. 
		The entries of $C\in\realR^{d^2\times k}$ are independently and randomly generated so that each entry of vector $C\xi$ falls into 
		$[0.09,0.092]$.
		The entries of the unit revenue $V\in\mathbb R^{d^2\times k}$ are independently and randomly generated so that each entry is on the interval $[4v_{ij},0]$ (note that $v_{ij} = -5$). 
		In particular, we randomly generate 10 instances for each setting as shown in Table \ref{tab:setting_gaps_10}. For each instance, we solve it using the 0-1 SDP approximation (SDP) proposed in Section \ref{sec:0-1SDP}, the inner approximation of the 0-1 COP formulation (IA-COP), and the MILP-based cutting-plane (MILP-Cut) method in Section \ref{sec:dcg},  by constructing the ambiguity set $\mathcal D_M$ with $N=10$ samples of the demand. {By generating 10 more additional samples of the demand,  we solve the same instance under the discrete ambiguity set $\mathcal D_\text{dis}$ with $N=20$ samples. } Given that the optimal value of the discrete ambiguity set serves as a valid lower bound of the problems under the continuous ambiguity set, an optimality gap of the LDR approximations can be computed based on the lower bound.  We denote the optimal value of SDP  (IA-COP, or MILP-Cut) as $V_\text{approx}$ and that of the discrete ambiguity set as $V_\text{dis}$. The optimality gap for SDP  (IA-COP, or MILP-Cut) is calculated as
		$$Gap = \frac{V_\text{approx}-V_\text{dis}}{|V_\text{dis}|}\times 100\%.$$
		All the models are solved with $(\gamma_1,\gamma_2) = (0,1)$. The results for $(\gamma_1,\gamma_2) = (1,1)$ are similar and
		are presented in Section   \ref{sec:cpu_11} of the Appendices.
		
		Table \ref{tab:setting_gaps_10}  reports the 25\%, 50\% and 75\% quantiles of the optimality gaps for the two approximations using linear decision rules.  When the uncertainty only appears in the constraints ($C=V=0$), the three approaches provide the same solutions and the same optimal values. Considering both random objective and constraints ($C,V\neq 0$), SDP and MILP-Cut provide better approximations than IA-COP. Table \ref{tab:cpu_10} summarizes the computational performance, across the same test instances as those reported in Table \ref{tab:setting_gaps_10}, of  using the cutting-plane algorithms and the MILP-Cut approach proposed in Section \ref{sec:benders}. The columns $t_\text{tot}$, \# It., and Gap report the average total run time, the average number of iterations (of the cutting-plane algorithm), and the average optimality gap, over the 10 instances of each setting. The numbers in bold are either the fastest run times or the smallest  gaps among the SDP, the IA-COP, and the MILP-Cut approaches. When $C=V=0$, the solutions obtained by solving the three approaches are almost the same and thus yield the similar gaps except that the MILP-Cut performs slightly better at for the last three settings. All three approaches have better performance (except for settings 1, 3, 8, and 10) in terms of faster computation time and smaller optimality gaps when $C=V=0$ compared to $C,V\neq 0$. In contrast, the computational performance of the discrete ambiguity set is similar in both cases. For either case ($C=V=0$ or $C,V\neq 0$), most of the time, MILP-Cut is faster than the other two approaches and
		it takes relatively longer time to solve SDP  than IA-COP.
		The optimality gaps of the three approaches are almost the same when $C=V=0$. Whereas when $C,V\neq 0$, SDP and MILP-Cut yield better gaps than IA-COP.
		Tables \ref{tab:setting_gaps_10} and \ref{tab:cpu_10} suggest that when the uncertainty appears only in the constraints ($C=V=0$), the decision maker may prefer MILP-Cut for better computational performance. Otherwise,  SDP and MILP-Cut are preferred for smaller gaps or IA-COP and MILP-Cut for faster computational time. 
		\begin{table}[htbp]
			\centering
			\caption{Quantiles of optimality gaps with $(\gamma_1,\gamma_2) = (0,1)$ }
			\resizebox{.95\textwidth}{!}{%
				\begin{tabular}{crrrrrr|rrr|rrr|rrr|rrr}
					\toprule
					\multirow{3}[1]{*}{Setting} & \multicolumn{1}{c}{\multirow{3}[1]{*}{$d$}} & \multicolumn{1}{c}{\multirow{3}[1]{*}{$l^S$}} & \multicolumn{1}{c}{\multirow{3}[1]{*}{$l^A$}} & \multicolumn{1}{c}{\multirow{3}[1]{*}{$N_B$}} & \multicolumn{2}{c|}{\multirow{2}[1]{*}{$w_i$}} & \multicolumn{3}{c|}{$C=V=0$} & \multicolumn{9}{c}{$C, V \neq 0$} \\
					&       &       &       &       & \multicolumn{2}{c|}{} & \multicolumn{3}{c|}{SDP/IA-COP/MILP-Cut} & \multicolumn{3}{c|}{SDP} & \multicolumn{3}{c|}{IA-COP} & \multicolumn{3}{c}{MILP-Cut} \\
					&       &       &       &       & \multicolumn{1}{c}{$C=V = 0$} & \multicolumn{1}{c|}{$C,V \neq 0$} & \multicolumn{1}{c}{25\%-Q} & \multicolumn{1}{c}{50\%-Q} & \multicolumn{1}{c|}{75\%-Q} & \multicolumn{1}{c}{25\%-Q} & \multicolumn{1}{c}{50\%-Q} & \multicolumn{1}{c|}{75\%-Q} & \multicolumn{1}{c}{25\%-Q} & \multicolumn{1}{c}{50\%-Q} & \multicolumn{1}{c|}{75\%-Q} & \multicolumn{1}{c}{25\%-Q} & \multicolumn{1}{c}{50\%-Q} & \multicolumn{1}{c}{75\%-Q} \\
					\hline
					1     & 15    & 5     & 2     & 3     & 750   & 1500  & 0.00  & 11.81 & 18.00 & 2.41  & 12.14 & 15.12 & 20.11 & 25.21 & 27.32 & 2.41  & 12.14 & 15.12 \\
					2     & 15    & 5     & 3     & 2     & 750   & 1500  & 0.00  & 18.68 & 26.65 & 20.88 & 24.75 & 66.21 & 30.64 & 35.20 & 72.50 & 20.65 & 24.75 & 66.21 \\
					3     & 15    & 5     & 2     & 3     & 1000  & 2000  & 7.03  & 15.29 & 21.83 & 3.56  & 17.73 & 21.70 & 26.45 & 32.42 & 35.69 & 3.50  & 17.73 & 21.70 \\
					4     & 15    & 5     & 3     & 2     & 1000  & 2000  & 0.00  & 23.22 & 33.80 & 0.00  & 28.96 & 36.39 & 0.00  & 37.59 & 43.08 & 0.00  & 25.65 & 30.59 \\
					5     & 20    & 5     & 2     & 3     & 750   & 1500  & 0.00  & 0.00  & 0.01  & 0.83  & 2.65  & 9.18  & 13.83 & 14.01 & 14.85 & 0.65  & 0.70  & 0.77 \\
					6     & 20    & 5     & 3     & 2     & 750   & 1500  & 0.00  & 11.48 & 45.23 & 2.02  & 17.20 & 46.25 & 16.30 & 20.65 & 54.11 & 0.81  & 6.49  & 46.22 \\
					7     & 20    & 5     & 2     & 3     & 1000  & 2000  & 0.00  & 0.00  & 0.01  & 0.97  & 3.51  & 9.94  & 14.40 & 15.26 & 15.78 & 0.70  & 0.74  & 0.78 \\
					8     & 20    & 5     & 3     & 2     & 1000  & 2000  & 0.00  & 25.37 & 54.88 & 3.99  & 21.38 & 53.63 & 19.35 & 30.79 & 62.74 & 1.02  & 14.42 & 53.60 \\
					9     & 25    & 10    & 5     & 5     & 750   & 1500  & 0.00  & 29.96 & 61.67 & 1.08  & 25.46 & 60.81 & 23.46 & 36.00 & 67.17 & 1.04  & 25.43 & 60.81 \\
					10    & 25    & 10    & 5     & 5     & 1000  & 2000  & 14.77 & 51.58 & 70.41 & 1.75  & 35.76 & 74.31 & 37.24 & 51.56 & 82.08 & 1.69  & 35.61 & 74.31 \\
					\bottomrule
				\end{tabular}%
			}
			\label{tab:setting_gaps_10}%
			
		\end{table}%

			\begin{table}[htbp]
				\centering
				\caption{Computational comparison of the three approximation approaches with $(\gamma_1,\gamma_2)=(0,1)$}
				\resizebox{\textwidth}{!}{%
					\begin{tabular}{c|rcr|rcr|rcr|rr|rrr|rrr|rcr|rr}
						\toprule
						\multicolumn{1}{r|}{\multirow{3}[2]{*}{Setting}} & \multicolumn{11}{c|}{$C=V=0$}                                                           & \multicolumn{11}{c}{$C, V \neq 0$} \\
						& \multicolumn{3}{c|}{SDP} & \multicolumn{3}{c|}{IA-COP} & \multicolumn{3}{c|}{MILP-Cut} & \multicolumn{2}{c|}{Discrete} & \multicolumn{3}{c|}{SDP} & \multicolumn{3}{c|}{IA-COP} & \multicolumn{3}{c|}{MILP-Cut} & \multicolumn{2}{c}{Discrete} \\
						& \multicolumn{1}{c}{$t_\text{tot}$ (s)} & \# It. & \multicolumn{1}{c|}{Gap (\%)} & \multicolumn{1}{c}{$t_\text{tot}$ (s)} & \# It. & \multicolumn{1}{c|}{Gap (\%)} & \multicolumn{1}{c}{$t_\text{tot}$ (s)} & \# It. & \multicolumn{1}{c|}{Gap (\%)} & \multicolumn{1}{c}{$t_\text{tot}$ (s)} & \multicolumn{1}{c|}{\# It.} & \multicolumn{1}{c}{$t_\text{tot}$ (s)} & \multicolumn{1}{c}{\# It.} & \multicolumn{1}{c|}{Gap (\%)} & \multicolumn{1}{c}{$t_\text{tot}$ (s)} & \multicolumn{1}{c}{\# It.} & \multicolumn{1}{c|}{Gap (\%)} & \multicolumn{1}{c}{$t_\text{tot}$ (s)} & \# It. & \multicolumn{1}{c|}{Gap (\%)} & \multicolumn{1}{c}{$t_\text{tot}$ (s)} & \multicolumn{1}{c}{\# It.} \\
						\hline    1     & 5.00  & 9.0   & \textbf{10.79} & 3.74  & 8.8   & \textbf{10.79} & \textbf{0.40} & 1.0   & \textbf{10.79} & 1.20  & 8.8   & 6.00  & 9.0   & \textbf{10.62} & 4.71  & 8.6   & 24.23 & \textbf{1.75} & 3.0   & \textbf{10.62} & 1.32  & 8.8 \\
						2     & 2.48  & 5.0   & \textbf{26.81} & 1.69  & 4.5   & \textbf{26.81} & \textbf{0.62} & 1.7   & \textbf{26.81} & 0.50  & 4.9   & 3.15  & 5.0   & 38.58 & \textbf{1.88} & 4.5   & 44.11 & 2.69  & 6.8   & \textbf{38.32} & 0.59  & 4.9 \\
						3     & 4.83  & 9.0   & \textbf{14.99} & 3.49  & 8.9   & \textbf{14.99} & \textbf{0.35} & 1.0   & \textbf{14.99} & 0.86  & 8.9   & 5.97  & 9.0   & \textbf{14.58} & 4.42  & 8.7   & 30.38 & \textbf{1.71} & 3.0   & \textbf{14.58} & 0.98  & 8.8 \\
						4     & 2.64  & 5.0   & \textbf{30.65} & 1.71  & 4.5   & \textbf{30.65} & \textbf{0.71} & 1.9   & \textbf{30.65} & 0.48  & 4.9   & 3.23  & 5.0   & 32.55 & \textbf{1.84} & 4.5   & 38    & 2.74  & 6.8   & \textbf{31.28} & 0.60  & 5.0 \\
						5     & 8.17  & 9.0   & \textbf{1.79} & 4.54  & 9.0   & \textbf{1.79} & \textbf{0.43} & 1.0   & \textbf{1.79} & 2.02  & 8.8   & 25.19 & 9.0   & 5.67  & 29.09 & 8.9   & 15.87 & \textbf{1.63} & 2.2   & \textbf{2.33} & 3.15  & 8.7 \\
						6     & 4.47  & 5.0   & \textbf{23.36} & 2.50  & 4.6   & \textbf{23.36} & \textbf{0.79} & 1.6   & \textbf{23.26} & 1.12  & 4.6   & \textbf{9.61} & 5.0   & 26.57 & 12.67 & 4.7   & 33.35 & 13.77 & 15.9  & \textbf{21.75} & 1.62  & 4.7 \\
						7     & 8.14  & 9.0   & \textbf{1.86} & 4.48  & 9.0   & \textbf{1.86} & \textbf{0.44} & 1.0   & \textbf{1.86} & 1.85  & 8.8   & 22.87 & 9.0   & 6.74  & 26.89 & 8.9   & 16.7  & \textbf{2.03} & 2.7   & \textbf{2.71} & 2.77  & 8.8 \\
						8     & 4.47  & 5.0   & 31.31 & 2.28  & 4.4   & 31.31 & \textbf{0.77} & 1.5   & \textbf{31.19} & 1.09  & 4.6   & \textbf{9.64} & 5.0   & 30.85 & 12.66 & 4.7   & 40.76 & 13.97 & 16.3  & \textbf{27.50} & 1.41  & 4.5 \\
						9     & 82.01 & 33.0  & 33.97 & 34.22 & 33.0  & 33.97 & \textbf{3.97} & 1.2   & \textbf{33.94} & 15.71 & 33.0  & 478.34 & 33.0  & \textbf{31.07} & \textbf{55.52} & 32.9  & 42.5  & 176.59 & 15.5  & 33.55 & 24.47 & 33.0 \\
						10    & 83.19 & 33.0  & 49.11 & 33.82 & 33.0  & 49.11 & \textbf{4.18} & 1.2   & \textbf{49.05} & 16.65 & 33.0  & 458.24 & 33.0  & 41.14 & \textbf{57.02} & 32.9  & 55.3  & 208.78 & 19.2  & \textbf{41.10} & 22.79 & 33.0 \\
						\bottomrule
					\end{tabular}%
				}
				\label{tab:cpu_10}%
			\end{table}%

			\subsection{In-Sample Results}
			\label{sec:sol-details}
			
			Given that in general the 0-1 SDP approximations in Section \ref{sec:0-1SDP} yield better gaps than IA-COP, we focus on the solution details of the 0-1 SDP approximations.
			Specifically,  Section \ref{sec:sol_profit} shows the solutions  and profits obtained using different values of $\gamma_1$ and $\gamma_2$. In Section \ref{sec:impact_of-gamma}, we further investigate how the profits change with $\gamma$'s and the support size. 
			
			\subsubsection{Solution details.}
			\label{sec:sol_profit}
			
			Table \ref{tab:sol_prof} presents the solution details when the uncertainty is (i) only in the constraints ($C=V=0$), and (ii) in both the objectives and constraints ($C,V\neq 0$) under various combination of $\gamma_1$ and $\gamma_2$. The  profit is the negative of the objective value \eqref{eq:DR-LDR-obj-sdp} associated with the optimal solution. For the optimal solution, only the first four components of $x$ are reported as they correspond to the four candidate locations for Company B to operate stores.
			
			In both cases ($C=V=0$ and $C,V\neq 0$), fewer stores are decided to open and the profits are non-increasing with larger $\gamma$'s. The solutions are invariant for a fixed $\gamma_1$ with the values of $\gamma_2$ considered here. In the next section, we show more results on how the profits change under various values of $\gamma$'s.

			We
			note that  the MILP model proposed in \citet{yanikoglu2018decision}  can be viewed as a DRBP problem under an distributional ambiguity set of distributions requiring to match a given mean and a given covariance. Let the mean and covariance be the same empirical mean and covariance used to solve the 0-1 SDP approximations for the MILP model. We solve for the optimal solution and the corresponding (estimated) profit using the MILP. Specifically, when $C=V=0$, the solutions and profits are the same as those of the 0-1 SDP  with $(\gamma_1,\gamma_2) = (0,1)$; when $C,V\neq 0$, the solution is the same as that of the 0-1 SDP with $(\gamma_1,\gamma_2) = (0,1)$, whereas, the profit is  1519.05, higher than that of the 0-1 SDP.

			\begin{table}[htbp]
				\centering
				\caption{Profits and solutions of the 0-1 SDP approximations under various $\gamma_1$ and $\gamma_2$}
				\resizebox{.95\textwidth}{!}{%

					\begin{tabular}{ccc|cc|cc|ccp{4.045em}|cp{4.045em}|cp{4.045em}}
						\toprule
						\multicolumn{7}{c|}{$C=V=0$}                            & \multicolumn{7}{c}{$C,V \neq 0$} \\
						\hline
						\multicolumn{1}{l}{$\gamma_1$} & \multicolumn{2}{c|}{0} & \multicolumn{2}{c|}{1} & \multicolumn{2}{c|}{1.5} & \multicolumn{1}{l}{$\gamma_1$} & \multicolumn{2}{c|}{0} & \multicolumn{2}{c|}{0.2} & \multicolumn{2}{c}{0.5} \\
						\hline
						\multicolumn{1}{l}{$\gamma_2$} & Profit  & Sol.  & Profit  & Sol.  & Profit  & Sol.  & \multicolumn{1}{l}{$\gamma_2$} & Profit  & \multicolumn{1}{c|}{Sol.} & Profit  & \multicolumn{1}{c|}{Sol.} & Profit  & \multicolumn{1}{c}{Sol.} \\
						1     & 808.97 & (0,0,1,1) & 344.07 & (0,1,0,0) & 344.07 & (0,1,0,0) & 1     & 1329.57 & (0,0,1,1) & 587.80 & (0,0,1,1) & 340.26 & (0,1,0,0) \\
						2     & 808.97 & (0,0,1,1) & 344.07 & (0,1,0,0) & 260.85 & (0,1,0,0) & 3     & 1255.79 & (0,0,1,1) & 513.98 & (0,0,1,1) & 270.06 & (0,1,0,0) \\
						\bottomrule
					\end{tabular}%

				}
				\label{tab:sol_prof}%
			\end{table}%
			\subsubsection{Impact of $\gamma$'s and support.}
			\label{sec:impact_of-gamma}
			
\begin{figure}[h]
				\centering
				\begin{subfigure}{0.45\textwidth}
					\centering
					\includegraphics[ width=0.85\linewidth]{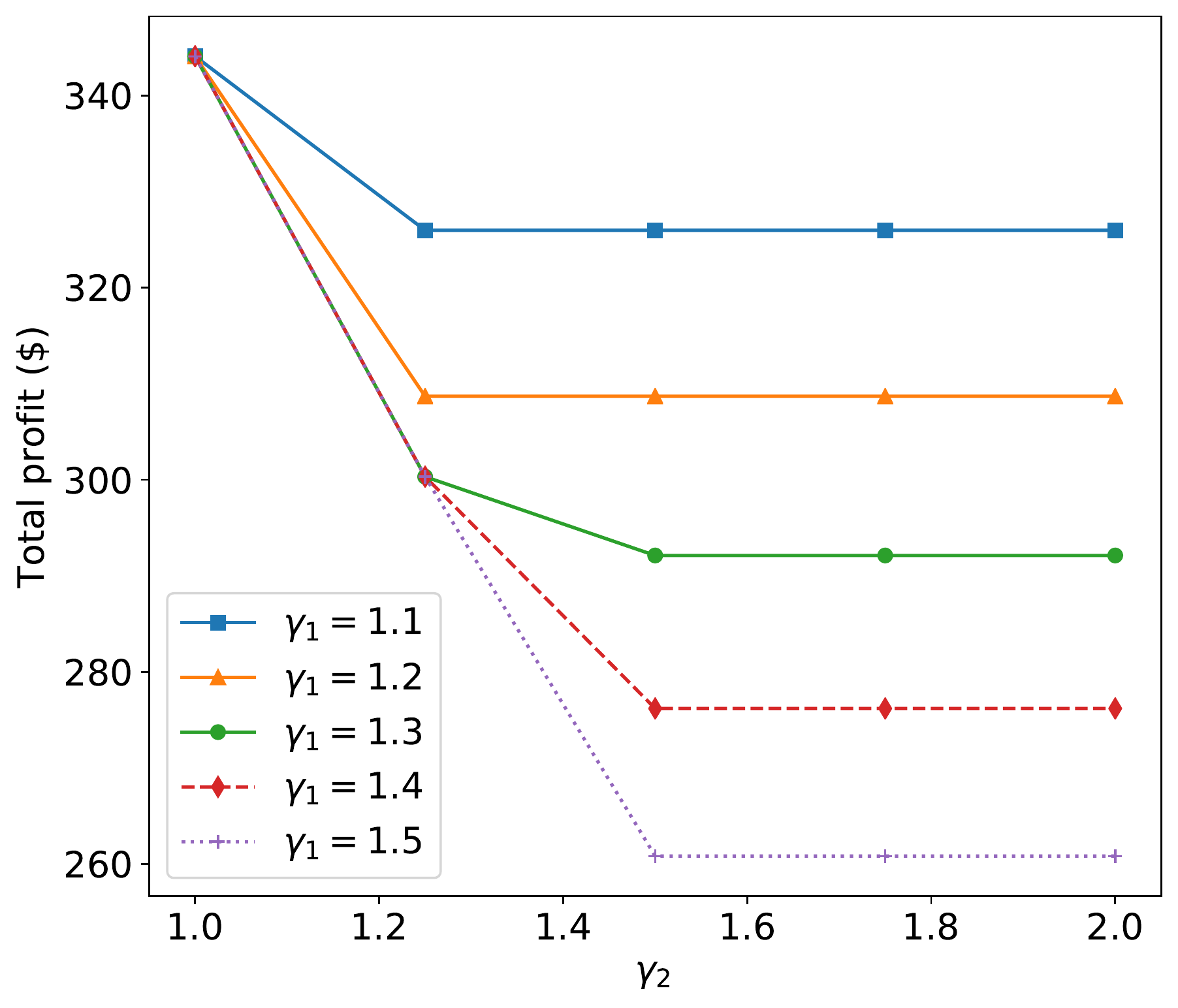}
					\caption{$C=V=0$}
					\label{fig:gamma_C_V_0}
				\end{subfigure}
				\begin{subfigure}{0.45\textwidth} 
					\centering
					\includegraphics[ width=.85\linewidth]{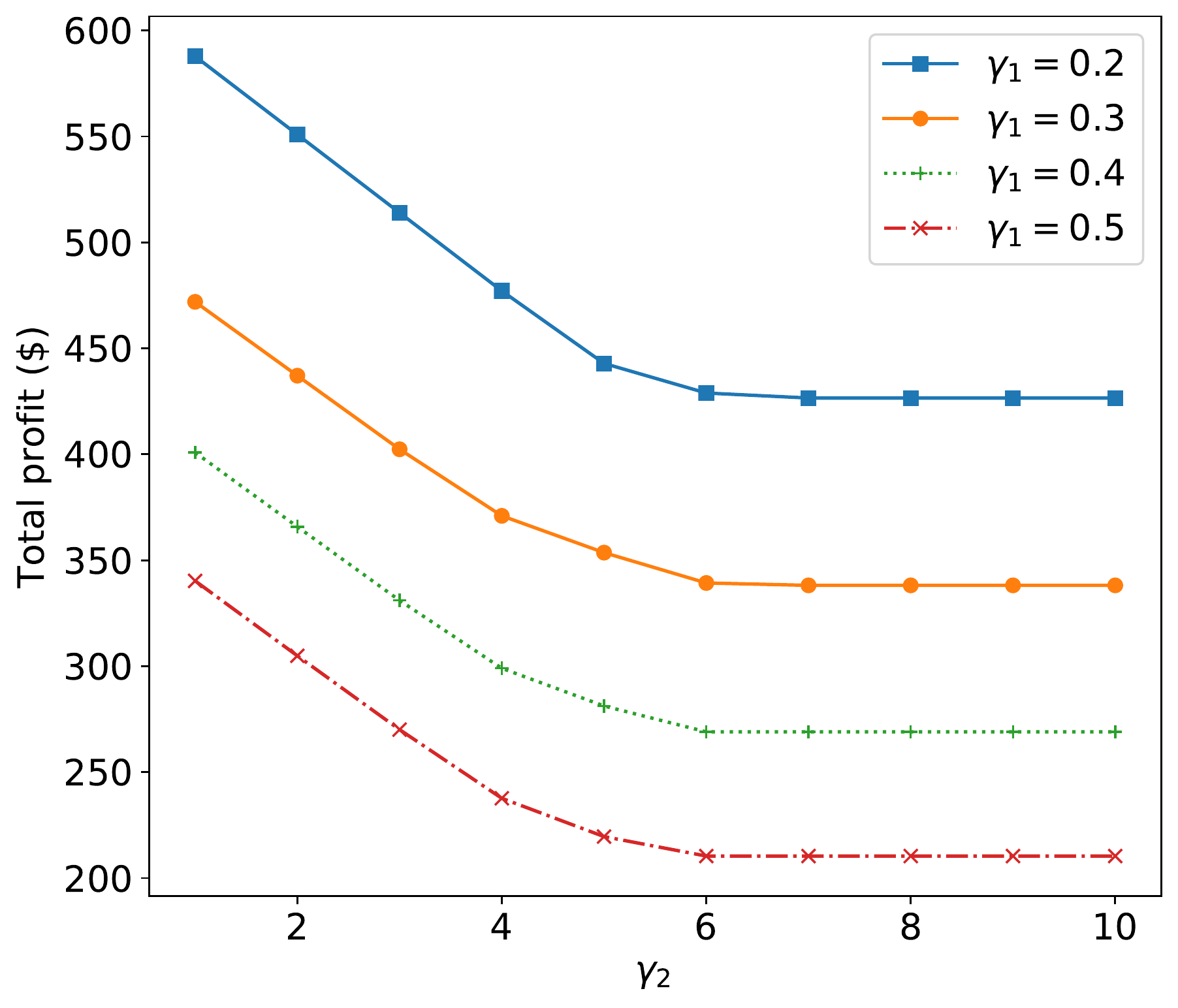}
					\caption{$C,V \neq 0$}
					\label{fig:gamma_C_V_N0}
				\end{subfigure}
				\centering
				\caption{Impact of $\gamma_1, \ \gamma_2$}
				\label{fig:gamma}
			\end{figure}

			In Figure \ref{fig:gamma}, each curve corresponds to the profits of a fixed $\gamma_1$ and a varying $\gamma_2$. The total profits are less when $\gamma$'s are larger as Company B needs to hedge against more (ambiguous) uncertainties.
			In both cases (Figures \ref{fig:gamma_C_V_0} and \ref{fig:gamma_C_V_N0}), for a fixed $\gamma_1$, the profits decrease at the beginning and remain unchanged. Comparing across different $\gamma_1$, the turning point of $\gamma_2$, beyond which the profits become unchanged, are non-decreasing as $\gamma_1$ increases. 
			The ambiguity set's parameter $\gamma_2$ impacts the profits by implicitly impacting  the uncertainties' dispersion via controlling the size of the ambiguity set. When $\gamma_2$ is large enough (at and beyond the turning point), the support size becomes the major factor on the uncertainties dispersion. 
			
\begin{figure}[h]
				\centering
				\begin{subfigure}{0.45\textwidth}
					\centering
					\includegraphics[height = .75\linewidth,  width=0.85\linewidth]{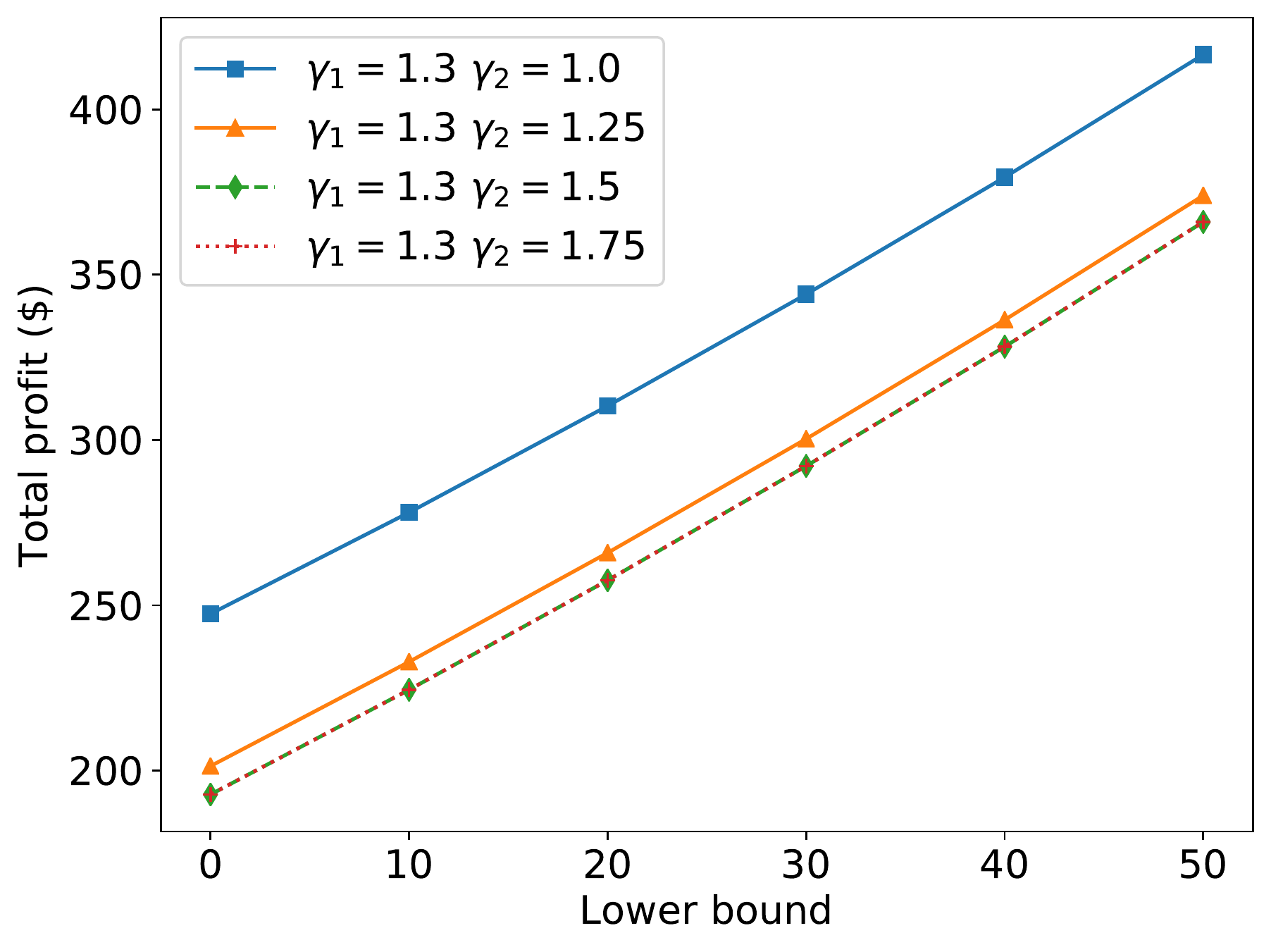}
					\caption{$C=V=0$}
					\label{fig:lb_C_V_0}
				\end{subfigure}
				\begin{subfigure}{0.45\textwidth} 
					\centering
					\includegraphics[height = .75\linewidth, width=.85\linewidth]{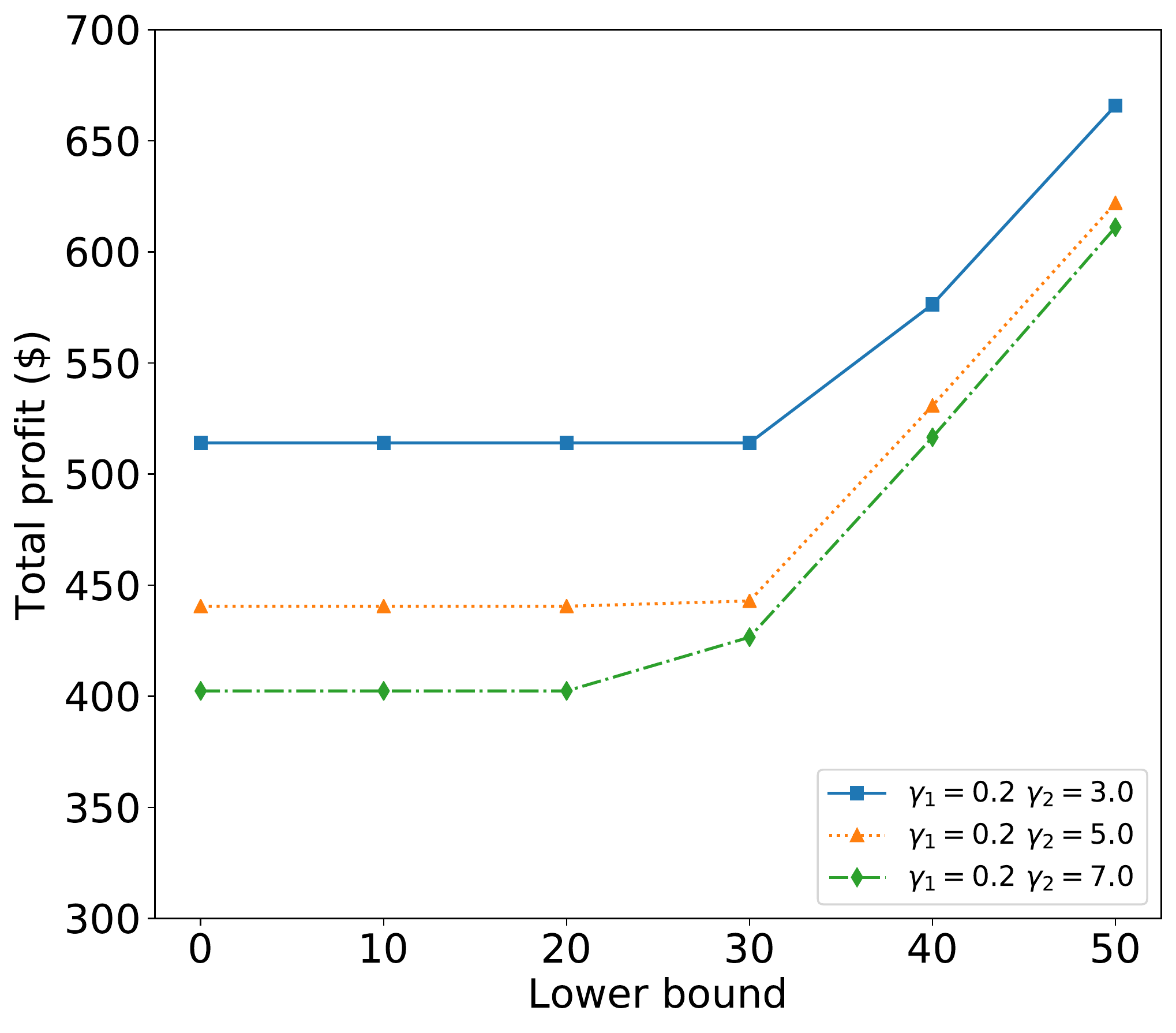}
					\caption{$C,V \neq 0$}
					\label{fig:lb_C_V_N0}
				\end{subfigure}
				\centering
				\caption{Impact of the support's lower bound}
				\label{fig:lb_C_V}
			\end{figure}

			Next, we investigate how the support size impacts the total profits by varying the lower bound of the support set. In Figure \ref{fig:lb_C_V}, each curve represents the total profits at a combination of $\gamma$'s using various lower bounds. When $C=V=0$, larger lower bounds yield higher profits as smaller uncertainties' dispersion is considered. When $C,V\neq 0$, for a fixed $(\gamma_1,\gamma_2)$, the profits remain the same when the lower bound is relatively small  and the profits start to climb up when continuing improving the lower bound.
			
			\subsection{Out-of-Sample Performance}
			\label{sec:out-of-sample}
			
			We perform out-of-sample simulations on the optimal solutions obtained by the pessimistic DRBP models under the same distribution used for the in-sample computation, i.e., uniform distribution on $[30,240]$. Specifically, we generate ten out-of sample sets, each containing $N^\prime = 5,000$ i.i.d. out-of-sample data $\xi_1^m,\ldots,\xi_d^m, \ m=1,\ldots,N^\prime$ of the demand vector $\xi$ following the uniform distribution. In Section \ref{sec:oos-misspecified} of the Appendices, we perform simulations on the out-of-sample data generated from distributions different than the in-sample uniform distribution to demonstrate the robustness when the distribution is misspecified. 
			
\begin{figure}[h]
				\centering
				\begin{subfigure}{0.45\textwidth}
					\centering
					\includegraphics[ width=1.\linewidth]{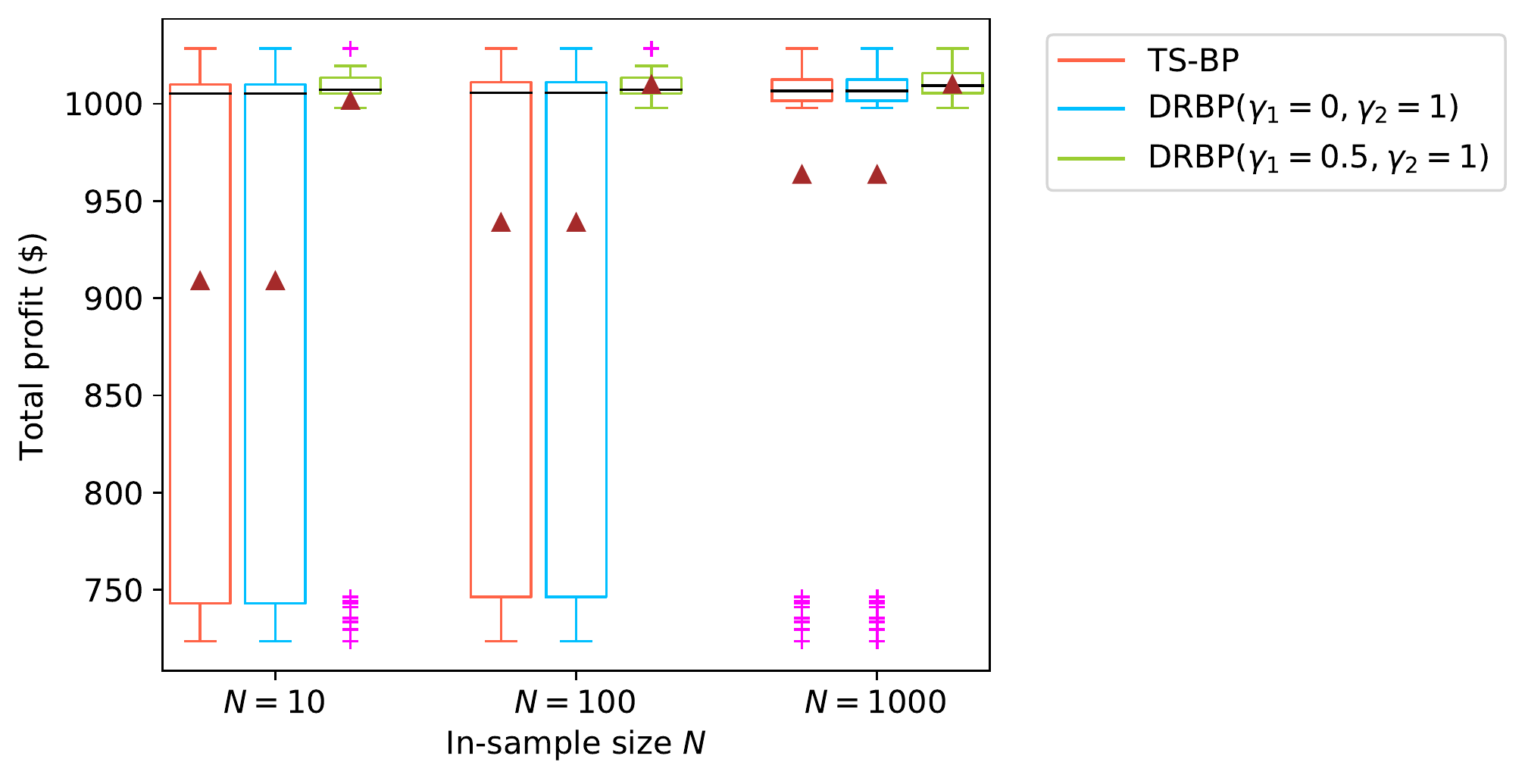}
					\caption{$C=V=0$}
					\label{fig:ofs-0}
				\end{subfigure}
				\hfill
				\begin{subfigure}{0.45\textwidth} 
					\centering
					\includegraphics[ width=1.\linewidth]{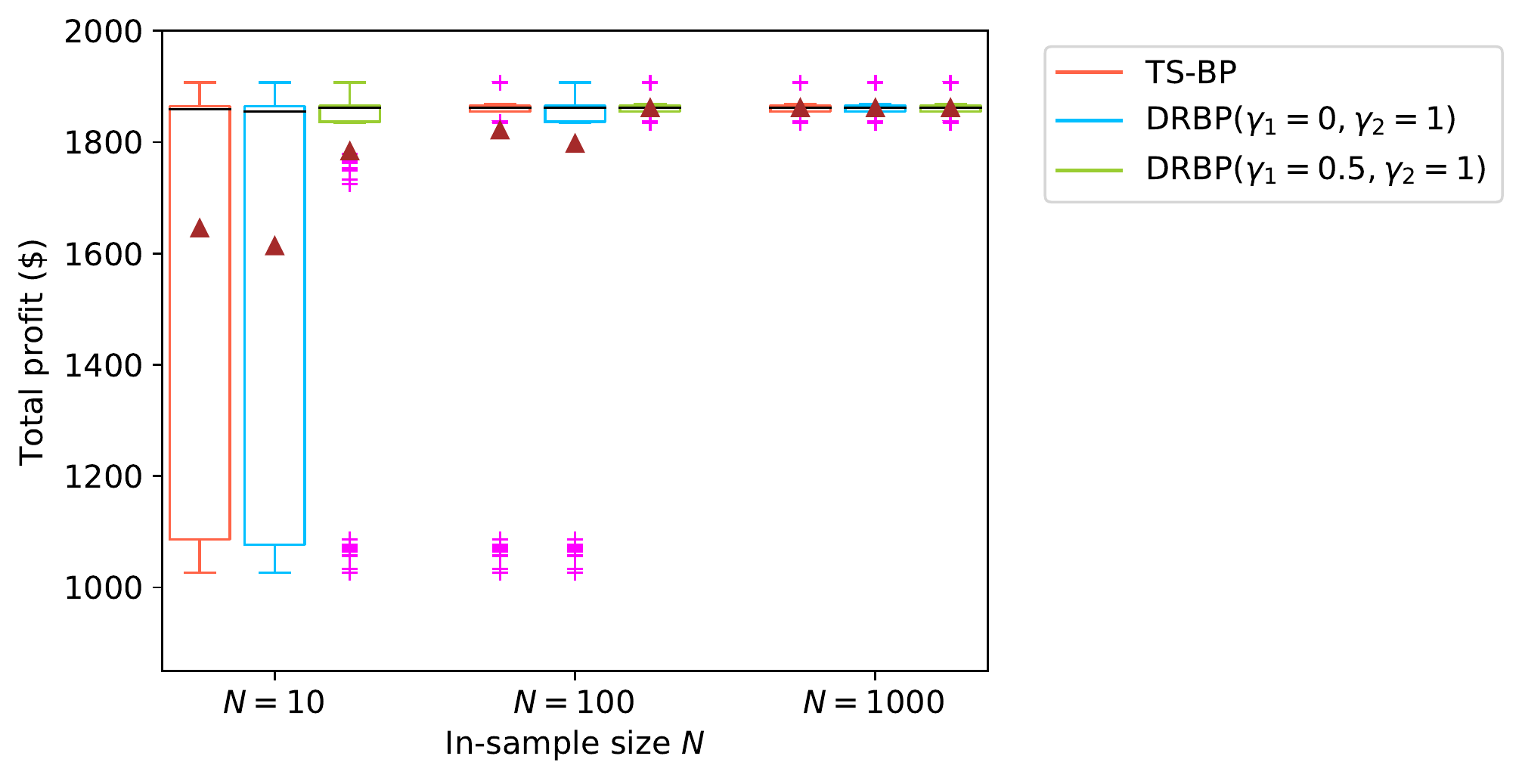}
					\caption{$C,V \neq 0$}
					\label{fig:ofs-n0}
				\end{subfigure}
				\centering
				\caption{Out-of-sample performance: expected profits}
				\label{fig:out-of-sample}
			\end{figure}

			Figure \ref{fig:out-of-sample} shows the total profits of the pessimistic DRBP models under two $\gamma$ combinations: $(\gamma_1,\gamma_2)= (0,1)$ and $(0.5,1)$. We benchmark our models with the MILP (TS-BP) proposed in
			\citet{yanikoglu2018decision}.
			The boxplots are generated by simulating the expected profits using the 100 solutions, each of an independent in-sample set with sizes $N = 10,100,1000$, on the ten out-of-sample sets. In both cases $C=V=0$ and $C,V\neq 0$, the expected profits vary in a wide range for TS-BP and DRBP with $(\gamma_1,\gamma_2)=(0,1)$. For larger $\gamma$'s, $(\gamma_1,\gamma_2)=(0.5,1)$, the expected profits show less variation. As the in-sample size $N$ increases, the variability of the expected profits reduces for all three models.

			\section{Conclusions}
			\label{sec:conclusion}
			
			This paper develops a general distributionally robust two-stage programming equivalence for the pessimistic bilevel program under proper conditions of  the ambiguity set. Typical choices of the ambiguity sets such as the moment-based and Wasserstein ball ambiguity sets satisfy the conditions under mild assumptions. The resulting two-stage formulation involves both uncertain objective and right-hand side. 
			The paper then focuses on the pessimistic DRBP when the leader chooses a binary decision under  moment-based ambiguity sets. Under the ambiguity set of continuous distributions, using LDRs,  a 0-1 SDP approximation and an exact 0-1 COP reformulation are provided. When the ambiguity set is restricted to only discrete distributions, an exact 0-1 SDP reformulation is developed and the worst-case distribution is explicitly constructed.
			To solve the resulting three 0-1 SDPs, we developed a Bender's type cutting-plane framework. Furthermore, another cutting-plane method is proposed based on the MILP approximation developed in \citet{yanikoglu2018decision}. Via numerical studies, we showed the efficiency and effectiveness of the proposed approximations. When the uncertainty  only presents in the constraints,  the approximations are solved faster and yield smaller gaps compared to the case when the uncertainty is both in the objectives and constraints. The  $\gamma$ parameters and the support size of the ambiguity sets play an important role in in-sample objectives. A reasonable choice of the $\gamma$ parameters has great impact on the variation of the out-of-sample objectives.
			
			Future research can be extended to distributionally robust bilevel programs with (1) integer or (and) nonlinear follower's problem, (2) continuous leader's decisions, (3) recourse leader's decisions. It would also be interesting to investigate the impact of the choice of ambiguity sets on the bilevel program's performance.
			
			\newpage
			\bibliographystyle{apalike}
			\bibliography{RO,YilingRef}
			
			
			
			
			\newpage
			\begin{APPENDICES}
				
				\section{Proof of Lemma \ref{lem:SECSTAG-PESSBI-CLP}}
				\label{sec:proof-to-function_space}
				
				\begin{proof}{Proof of Lemma \ref{lem:SECSTAG-PESSBI-CLP}:}
					First, we have $\expec_F[\Psi(x,\xi)]\ge {Q}_F(x)$ as the recourse variables $y(\xi)$ are restricted to a smaller class of function map. It remains to show that there exists a $y^*(\xi)\in\underset{y^{\prime}(\xi)\in\mathcal{L}^2_n(F)}{\arg\min}\ \left\{c(\xi)^\top y^{\prime}(\xi):Ay^{\prime}(\xi)\le b_x(\xi)\right\}$ such that $\expec_F[v(\xi)^\top y^*(\xi)]\ge\expec_F[\Psi(x,\xi)]$.\par
					According to Assumption \ref{asp:RECOURSE}, for a given leader's decision $x\in\mathcal X$ and any given $\xi\in\mathcal S$, the set of optimal solutions of the follower's problem $\Omega(x,\xi)$
					is a polyhedral set spanned over a subset of the extreme points of the follower's feasible region $\mathcal{Y}(\xi): = \left\{Ay \le b_x(\xi)\right\}$. Hence, one optimal solution of $\Psi(x,\xi)$  must be obtained at some extreme point of the  feasible region $\mathcal{Y}(\xi)$. There exist finitely many matrices $L_j^y\in\mathbb R^{m\times n}, \ j = 1,\ldots, J_1$, for every basis of the constraint matrix $A$, such that each extreme point of $\mathcal{Y}(\xi)$ is represented as $L_j^y b_x(\xi)$. For every $x\in\mathcal{X}$ and $\xi\in\mathcal{S}$, we define the following index set.
					\begin{align*}
						\begin{array}{rl}
							\mathcal{J}^y(\xi,x);=\Big\{j: \ L_j^yb_x(\xi)\in\mathcal{Y}(\xi),\ c(\xi)^\top\left(  L_{j}^yb_x(\xi)-L_{j^\prime}^yb_x(\xi)\right)\le0 ,\ \forall 1\le j^{\prime}\le J_1 \text{ with } L_{j^\prime}^yb_x(\xi)\in\mathcal{Y}(\xi)\Big\},
						\end{array}
					\end{align*}
					where the set $\mathcal{J}^y(\xi,x)$ consists of all the indices $j$ such that $L_j^yb_x(\xi)$ is an optimal solution in the set $\Omega(x,\xi)$ of optimal solutions to the follower's problem. We then consider a partition of the support set $\mathcal S$ of $\xi$ using  $\mathcal{J}^y(\xi,x)$. 
					\begin{align*}
						\begin{array}{rl}
							\Xi_j^{y}(x)=&\Big\{\xi\in\realR^k: \ v(\xi)^\top\left(L_j^{y}b_x(\xi)-L^{y}_{j^\prime}b_x(\xi)\right)\ge 0,\ \forall j^\prime\in\{1,\ldots,j-1\}\cap \mathcal{J}^y(\xi,x),\\
							&v(\xi)^\top\left(L_j^{y}b_x(\xi)-L^{y}_{j^\prime}b_x(\xi)\right)> 0,\ \forall j^\prime\in\{j+1,\ldots,J_1\}\cap \mathcal{J}^y(\xi,x)\Big\}.
						\end{array}
					\end{align*}
					The partition $\Xi_j^y(x)$ consists of all the $\xi$ such that $j$ is the largest index for which $j\in\mathcal{J}^y(\xi,x)$ and $y^*(\xi):=L_j^yb_x(\xi)$ maximizes $v(\xi)^\top y$. 
					Moreover,
					\begin{subequations}\label{eqn:SECSTAG-PESSBI-BNDYS}
						\begin{eqnarray*}
							\expec_F[\|y^*(\xi)\|^2]&\le&\expec_F\left[\underset{1\le j\le {J_1}}{\max}\|L_j^{y}b_x(\xi)\|^2\right]=\expec_F\left[\underset{1\le j\le {J_1}}{\max} \left\{b_x(\xi)^\top (L_j^{y})^\top L_j^{y} b_x(\xi)\right\}\right] \label{eqn:SECSTAG-PESSBI-BNDYS1}\\
							&\le&\underset{1\le j\le {J_1}}{\max}\|L_j^y\|^2\expec_F[\|b_x(\xi)\|^2]<\infty\label{eqn:SECSTAG-PESSBI-BNDYS2}
						\end{eqnarray*}
					\end{subequations}
					The first inequality holds since $\|y^*(\xi)\|\le \underset{1\le j\le {J_1}}{\max}\|L_j^{y}b_x(\xi)\|$  for any given $\xi\in\Xi_j^y(x)$. The second inequality holds because $b_x(\xi)^\top (L_j^{y})^\top L_j^{y} b_x(\xi)\le\|L_j^y\|^2\|b_x(\xi)\|^2$ for any $\xi$, where the norm of a matrix refers to the Frobenius norm, i.e., $\|X\| = \sqrt{\sum_{i=1}^m \sum_{j=1}^n X_{ij}^2}$, for $X\in\mathbb R^{m\times n}$.
					Thus, $y^*(\xi) \in \mathcal L_n^2(F)$ and is feasible to the maximization problem of $Q_F(x)$. Therefore, $Q_F(x)$ is no less than $\mathbb E_F[\Psi(x,\xi)]$ and  we conclude the proof. \Halmos
				\end{proof}
				
				\section{Proof of Lemma \ref{lem:SECSTAG-PESSBI-CLP-DUAL}}
				\label{sec:proof-lemma-to-real}
				\begin{proof}{Proof of Lemma \ref{lem:SECSTAG-PESSBI-CLP-DUAL}: }
					Define%
					\begin{align*}
						\varphi_1(x,\lambda):=&\underset{p(\xi)\in\mathcal{L}^2_m(F)}{\inf}\ \left\{ \expec_F[b_x(\xi)^\top p(\xi)]: \ A^\top p(\xi)+\lambda c(\xi)=v(\xi),\ p(\xi)\ge0\right\},\\
						\varphi_2(x,\lambda):=&\underset{y(\xi)\in\mathcal{L}^2_n(F)}{\inf}\ \left\{ \expec_F[c(\xi)^\top y(\xi)]: \ Ay(\xi)\le \lambda b_x(\xi)\right\}.
					\end{align*}
					To	establish the equivalence between problem \eqref{model:SECSTAG-PESSBI-CLP-ELIARG}   and \eqref{model:SECSTAG-PESSBI-DUAL}, we first will show that (i) problem \eqref{model:SECSTAG-PESSBI-CLP-ELIARG}   is equivalent to the following minimization problem over a scaler variable $\lambda\ge0$, $p(\xi)\in\mathcal L_m^2(F)$, and $y(\xi)\in\mathcal L_n^2(F)$.
					\begin{align}\label{eq:functional-lambda}
						\underset{\lambda\ge0}{\min}\ \varphi_1(x,\lambda)+\varphi_2(x,\lambda),
					\end{align}
					Then it remains to (ii) show the equivalence between \eqref{eq:functional-lambda} and \eqref{model:SECSTAG-PESSBI-DUAL}.
					
					\begin{enumerate}
						\item[(i)] Equivalence between \eqref{model:SECSTAG-PESSBI-CLP-ELIARG} and \eqref{eq:functional-lambda}. For this part, the proof is inspired by the proof of Theorem 4.1 in \citet{yanikoglu2018decision}. Dualizing problem \eqref{model:SECSTAG-PESSBI-CLP-ELIARG} yields the following dual problem.
						\begin{equation}\label{eq:pess-dual-intermediant}
							\inf_{\lambda,y(\xi)\in\mathcal{L}^2_n(F)}\varphi_1(x,\lambda)+\lambda\mathbb{E}_F[\bar Q_F(x)]
						\end{equation}
						The strong duality holds for the primal and dual pair \citep{ruszczynski2006optimization}. Given that $\lambda\ge 0$ and the definition of $\bar Q_F(x)$ in Proposition \ref{prop:yanikoglu}, problem \eqref{eq:pess-dual-intermediant} is rewritten as
						\begin{subequations}\label{eq:pess-dual-intermediant2}
							\begin{eqnarray}
								\inf_{\lambda,y(\xi)\in\mathcal{L}^2_n(F)}&&	\varphi_1(x,\lambda)+\lambda\mathbb{E}_F[c(\xi)^\top y(\xi)]\\
								\mbox{s.t.}
								&& \ Ay(\xi)\le\ b_x(\xi).
							\end{eqnarray}
						\end{subequations}
							Problem \eqref{eq:pess-dual-intermediant2} is nonlinear because of the product of $\lambda$ and $y(\xi)$ in the objective. We linearize it by letting $y(\xi)$ be $\lambda y(\xi)$ and obtain the reformulation \eqref{eq:functional-lambda}.
								The reformulation is exact. First, given a feasible solution $(\lambda, y, p)$ to \eqref{eq:pess-dual-intermediant2}, then $(\lambda, \lambda y, p)$ is feasible to \eqref{eq:functional-lambda} with the same objective value. Now, consider a feasible solution $(\lambda, y, p)$ to \eqref{eq:functional-lambda}. If $\lambda > 0$, then $(\lambda, y/\lambda, p)$ is feasible to \eqref{eq:pess-dual-intermediant2} with the same objective value. If $\lambda = 0$, the constraint of $\varphi_2(x,\lambda)$   becomes $Ay(\xi) \le 0$ and thus implies that $y(\xi) = 0$ almost surely due to Assumption \ref{asp:RECOURSE}. Hence, $(0,y^\prime, p)$, with any $y^\prime \in \mathcal L_n^2$ such that $Ay^\prime(\xi)\le b_x(\xi)$ holds almost surely, is feasible to \eqref{eq:pess-dual-intermediant2} and achieves the same objective value as $(\lambda, y,p)$. Therefore,  \eqref{eq:functional-lambda} is an exact reformulation of \eqref{eq:pess-dual-intermediant2} or, equivalently, \eqref{model:SECSTAG-PESSBI-CLP-ELIARG}. 
								
								\item[(ii)] Equivalence between \eqref{eq:functional-lambda} and \eqref{model:SECSTAG-PESSBI-DUAL}.
								To establish the equivalence, we need to show the following equivalence.
								\begin{align}\label{model:SECSTAG-PESSBI-DUAL-PHIVAR1}
									\varphi_1(x,\lambda)&=\phi_1(x,\lambda):=\expec_F\left[\underset{p\in\realR^m}{\min}\left\{b_x(\xi)^\top p\ |\ A^\top p+\lambda c(\xi)=v(\xi),p\ge0\right\}\right],\\
									\label{model:SECSTAG-PESSBI-DUAL-PHIVAR2}
									\varphi_2(x,\lambda)&=\phi_2(x,\lambda):=\expec_F\left[\underset{y\in\realR^n}{\min}\left\{c(\xi)^\top y\ |\ Ay\le \lambda b_x(\xi)\right\}\right].
								\end{align}
								Note that $\varphi_1(x,\lambda)$ and $\varphi_2(x,\lambda)$ are bounded below as the second-stage expected cost is bounded below with respect to distribution $F\in\mathcal D$ by virtue of Assumption \ref{asp:RECOURSE}.
								\begin{enumerate}
									\item[\eqref{model:SECSTAG-PESSBI-DUAL-PHIVAR1}:] 
									It is clear that the right-hand side $\phi_1(x,\lambda)$ is no greater than $\varphi_1(x,\lambda)$ as $\varphi_1(x,\lambda)$  is more restrictive on the function map $p(\xi)$. It remains to show that $\varphi_1(x,\lambda)$ is no greater than $\phi_1(x,\lambda)$. 
									
									We first show that $\phi_1(x,\lambda)$ is bounded below. 
									Denote $u(x, \lambda,\xi) : = \min_{p\in\mathbb R^m}\{b_x(\xi)^\top p : \ A^\top p + \lambda c(\xi) = v(\xi), \ p\ge 0\}$.
									Rewrite $u(x, \lambda,\xi)  = \max_{\eta\in\mathbb R^n} \left\{(v(\xi) - \lambda c(\xi))^\top \eta: \ A\eta \le b_x(\xi)\right\}$ with the dual variable $\eta$. Due to Assumption \ref{asp:RECOURSE}, $u(x, \lambda,\xi)$ is bounded for every $\xi$ and $x$.

											Next, we establish the equivalence of \eqref{model:SECSTAG-PESSBI-DUAL-PHIVAR1}.
											We form the second partition $\left\{\bar\Xi_j^{p}(x,\lambda)\right\}_{j=1}^{J_2}$ of $\mathcal{S}$ for $x\in\mathcal{X}$ and $\lambda\ge 0$, such that $j$ is the largest index so that $L_j^{p}d(\lambda,\xi)$ is an optimal solution of $u(x,\lambda,\xi)$.
											\begin{align*}
												\begin{array}{rl}
													\bar\Xi_j^{p}(x,\lambda):=&\Big\{\xi\in\realR^k: \ A^\top L_j^{p}d(\lambda,\xi)=d(\lambda, \xi),\ L_j^{p}d(\lambda, \xi)\ge0\\
													& b_x(\xi)^\top\left(L_j^{p}-L^{p}_{j^\prime}\right)d(\lambda,\xi)\le 0\ \forall j^\prime=1,\ldots,j-1\ \text{ with } A^\top L_{j^\prime}^{p}d(\lambda,\xi)=d(\lambda,\xi),\ L_{j^\prime}^{p}d(\lambda,\xi)\ge0,\\
													& b_x(\xi)^\top\left(L_j^{p}-L^{p}_{j^\prime}\right)d(\lambda,\xi)\le 0\ \forall j^\prime=j+1,\ldots,{J_2}\ \text{ with } A^\top L_{j^\prime}^{p}d(\lambda,\xi)=d(\lambda,\xi),\ L_{j^\prime}^{p}d(\lambda,\xi)\ge0\Big\}.
												\end{array}
											\end{align*}
											Denote $p^*(x,\lambda,\xi): = L_j^p d(\lambda, \xi)$ for $\xi \in \bar{\Xi}(x,\lambda), \ j=1,\ldots, J_2$. We can show that $p^*(x,\lambda,\xi)\in\mathcal L_m^2(F)$ and the equivalence between $\phi_1(x,\lambda)$ and $\varphi_1(x,\lambda)$ is established. The proof is similar to that of Lemma \ref{lem:SECSTAG-PESSBI-CLP} and thus is omitted. 
											\item[\eqref{model:SECSTAG-PESSBI-DUAL-PHIVAR2}:] According to Assumption \ref{asp:RECOURSE}, $\{y\in\realR^n | Ay\le\lambda b_x(\xi)\}$ is a compact polyhedron and thus $\underset{y\in\realR^n}{\min}\left\{c(\xi)^\top y|Ay\le \lambda b_x(\xi)\right\}$  is bounded below almost surely for $\xi\in \mathcal S$. Then following the similar proof of Lemma \ref{lem:SECSTAG-PESSBI-CLP}, we establish the equivalence as desired. 
										\end{enumerate}
										
									\end{enumerate}
									Therefore, the proof is complete. 
									\Halmos
								\end{proof}

								\section{Proof of Theorem \ref{thm:sup-min-interchangeability}: Interchangeability between the Minimization and Maximization}
								\label{sec:proof-interchangeability}
								
								\begin{proof}{Proof of Theorem \ref{thm:sup-min-interchangeability}:} As Lemma \ref{lem:SECSTAG-PESSBI-CLP-DUAL} implies that $\sup_{F\in\mathcal D} \mathbb E_F [\Psi(x,\xi)] = \sup_{F\in\mathcal D} \min_{\lambda\ge 0} \mathbb E_F[\Phi_\lambda (x,\xi)]$, to prove Theorem \ref{thm:sup-min-interchangeability}, it suffices to show that 
									\begin{equation*}
										\sup_{F\in\mathcal D} \min_{\lambda\ge 0} \mathbb E_F[\Phi_\lambda (x,\xi)] =  \min_{\lambda\ge 0} \sup_{F\in\mathcal D} \mathbb E_F[\Phi_\lambda (x,\xi)].
									\end{equation*}
									Consider a weakly compact ambiguity set $\mathcal D$ of probability measures on $(\mathcal S,\mathcal F)$ and a convex neighborhood $N_\delta : = \{\lambda: \ \lambda \ge -\delta\}$ of the feasible region of $\lambda$, where $\delta>0$.  According to the minimax analysis of stochastic problems (see Theorem 2.1 of \citet{shapiro2002minimax}, Theorem 49 of \citet{ruszczynski2003optimality}), 
									we need to prove the following statements to show the interchangeability of the maximization and minimization.
									\begin{enumerate}
										\item For $\lambda \in N_\delta$, $\sup_{F\in\mathcal D}\mathbb E_F [\Phi_\lambda (x,\xi)] < +\infty$.
										\item 
										For every $\xi\in\mathcal S$, function $\Phi_\lambda (x,\xi)$ is convex for $\lambda\in N_\delta$.
										\item Given an optimization solution $\bar\lambda$  of \eqref{eq:min-max-lambda}, for every $\lambda$ in a neighborhood of $\bar\lambda$, the function $\Phi_{\lambda}(x,\xi)$ is bounded and upper semicontinuous on $\xi\in\mathcal S$ and function $\Phi_{\bar\lambda}({x},\xi)$ is bounded and continuous on $\xi\in \mathcal S$.
									\end{enumerate}
									
									\noindent \underline{Statement 1.} 
									Denote $\Phi_\lambda^1(x,\xi) := \min_{y}\left\{ c(\xi)^\top y: \ Ay\le \lambda b_x(\xi)\right\}$ and $\Phi_\lambda^2(x,\xi): = \min_{p}\left\{b_x(\xi)^\top p: \ A^\top p+ \lambda c(\xi) = v(\xi), \ p\ge 0\right\}$. To show that $\sup_{F\in\mathcal D}\mathbb E_F [\Phi_\lambda (x,\xi)] < +\infty$, it suffices to show that  $\Phi_\lambda (x,\xi) = \Phi_\lambda^1 (x,\xi) + \Phi^2_\lambda (x,\xi)$. Due to  Assumption \ref{asp:RECOURSE}, $\Phi_\lambda^1(x,\xi)$ is bounded for every $\xi$ and $x$. It remains to prove that $\Phi_\lambda^2 (x,\xi)$ is bounded. Rewrite $\Phi_\lambda^2 (x,\xi) = \max_\eta \left\{(v(\xi) - \lambda c(\xi))^\top \eta: \ A\eta \le b_x(\xi)\right\}$ with the dual variable $\eta$. Due to Assumption \ref{asp:RECOURSE}, $\Phi_\lambda^2(x,\xi)$ is bounded for every $\xi$ and $x$. Thus $\sup_{F\in\mathcal D}\mathbb E_F[\Phi_\lambda(x,\xi)]<+\infty$.
									
									
									\noindent \underline{Statement 2.}  Given any $\xi\in\mathcal S$, we derive the dual problem of $\Phi_\lambda(x,\xi)$.
									\begin{align*}
										\Phi_\lambda(x,\xi)= & \max_{\varrho\in\realR^m,\eta\in\realR^n} v(\xi)^\top\eta-\lambda c(\xi)^\top\eta-\lambda b_x(\xi)^\top\varrho\\
										& \quad \st\ A^\top\varrho+c(\xi)=0,\ A\eta\leq b_x(\xi),\ \varrho\geq0,
									\end{align*}
									where $\eta\in\mathbb R^n$ and $\varrho \in\mathbb R^m$ are the dual variables. As $\Phi_\lambda(x,\xi)$ is bounded below, the optimal value $\Phi_\lambda(x,\xi)$ can be view as a pointwise maximum of (finitely many) linear functions. Hence $\Phi_\lambda(x,\xi)$ is a convex function of $\lambda\in N_\delta$.
									
									\noindent \underline{Statement 3.} As Statement 1 shows that $\Phi_\lambda(x,\xi)$ is bounded below for any $\xi\in\mathcal S$,  it remains to show that, for any $\lambda\in N_\delta$, $\Phi_\lambda (x,\xi)$ is Lipschitz continuous on $\xi\in\mathcal S$. 
									%
									%
									Given a $\lambda\ge 0$, consider any $\xi_1\neq \xi_2$,  $\xi_1,\xi_2\in\mathcal S$. Without loss of generality, assume that $\Phi_\lambda(x,\xi_1) \le \Phi_\lambda(x,\xi_2)$. Let $(y^*,p^*)$ denote the optimal solution of $\Phi_\lambda(x,\xi_1)$. We then rewrite $\Phi_\lambda (x,\xi_2)$ as
									\begin{subequations}
										\begin{eqnarray}
											\min_{\Delta y\in\realR^n,\Delta p\in\realR^m} && b(\xi_2)^\top (p^* + \Delta p) + c(\xi_2)^\top (y^* + \Delta y)\\
											\label{eq:reformulate-xi2_p-constr}		\mbox{s.t.} && A^\top (p^* + \Delta p) + \lambda c(\xi_2) = v(\xi_2), \ p^*+\Delta p \ge 0\\
											\label{eq:reformulate-xi2_y-constr}	&& A(y^*+\Delta y) \le \lambda b_x(\xi_2).
										\end{eqnarray}
									\end{subequations}
									Note that $(y^*,p^*)$ is feasible  to $\Phi_\lambda (x,\xi_1)$ and thus $A^\top p^* + \lambda c(\xi_1) = v(\xi_1), \ p^* \ge 0 , \ Ay^*\le \lambda b_x(\xi_1)$. Define
									\begin{eqnarray}
										\label{eq:reformulation_xi2_y}\phi_\lambda^y(x,\xi_1,\xi_2) :=	 &&\min_{\Delta y \in\mathbb R^n} \left\{c(\xi_2)^\top \Delta y: \ A\Delta y\leq \lambda (b_x(\xi_2) - b_x(\xi_1))	\right\},\\
										\label{eq:reformulation_xi2_p}\phi_\lambda^p (x,\xi_1,\xi_2): = &&\min_{\Delta p\in\mathbb R^m}  \left\{b(\xi_2)^\top \Delta p: \ \Delta p\geq0,\ A^\top \Delta p=(v(\xi_2) - v(\xi_1)) - \lambda (c(\xi_2) - c(\xi_1))  	\right\}.
									\end{eqnarray}
									Then $\Phi_\lambda (x,\xi_2) \le 
									b(\xi_2)^\top p^* + c(\xi_2)^\top y^*
									+\phi_\lambda^y(x,\xi_1,\xi_2)
									+ \phi_\lambda^p(x,\xi_1,\xi_2)$ because of two reasons: (1) the constraint in \eqref{eq:reformulation_xi2_y}  is more restrictive than \eqref{eq:reformulate-xi2_y-constr} as $\lambda(b_x(\xi_2) - b_x(\xi_1)) \le \lambda b_x(\xi_2) - Ay^*$, and (2) the constraint $\Delta p \ge 0$ of \eqref{eq:reformulation_xi2_p}is more restrictive than $\Delta p \ge - p^*$ in \eqref{eq:reformulate-xi2_p-constr}. 
									
									Now, consider
									\begin{subequations}\label{eq:lip-bound0}
										\begin{eqnarray}
											|\Phi_\lambda(x,\xi_2) - \Phi_\lambda (x,\xi_1)| &=& \Phi_\lambda(x,\xi_2) - \Phi_\lambda (x,\xi_1)\\
											\label{eq:lip-bound}	& \le& |\left(b(\xi_2) - b(\xi_1)\right)^\top p^*| + |\left(c(\xi_2) - c(\xi_1)\right)^\top y^*| + |\phi_\lambda^y(x,\xi_1,\xi_2)| + |\phi_\lambda^p(x,\xi_1,\xi_2)|.\hspace{10mm}
										\end{eqnarray}
									\end{subequations}
									We then show that each term in \eqref{eq:lip-bound} is bounded above by a product of a constant and $\|\xi_2 - \xi_1\|$ for Lipschitz continuity. Recall the linearity assumption \ref{asp:Linear} where 
									$c(\xi)=C\xi+{c_0}$ for $C\in\realR^{n\times k}$, {$c_0\in\realR^n$}, $v(\xi)=V\xi+{v_0}$ for $V\in\realR^{n\times k}$, {$v_0\in\realR^n$}, $b_x(\xi)=B_x\xi+{b_{x0}}$ for $B_x\in\realR^{m\times k}$, {$b_{x0}\in\realR^m$}.
									
									\begin{enumerate}
										\item $|\left(b(\xi_2) - b(\xi_1)\right)^\top p^*|$:  Similarly as in the proof of Lemma \ref{lem:SECSTAG-PESSBI-CLP-DUAL}, there are finitely many matrices $L_j^p\in\mathbb R^{n\times m}, \ j=1,\ldots, J_2$, one for every basis of the constraint matrix $A^\top$, so that $L_j^p \left(v(\xi_1) - \lambda c(\xi_1)\right)$ represents a basic feasible solution in $\left\{p\in\mathbb R^m: \ A^\top p +\lambda c(\xi_1) = v(\xi_1), \ p\ge 0 \right\}$. Let $L_{j^*}^p$ be the matrix associated with $p^*$. We then have 
										\begin{equation}\label{eq:lip-bound1}
											|\left(b(\xi_2) - b(\xi_1)\right)^\top p^*| \le \|B_x(\xi_2-\xi_1)\|\|L_{j^*}^p\left(v(\xi_1) - \lambda c(\xi_1)\right)\|\le L_1 \|\xi_2-\xi_1\|,
										\end{equation}
										where $L_1 := \|B_x\| \max_{1\le j\le J_2}\max_{\xi\in\mathcal S} \|L_j^p\left(v(\xi) - \lambda c(\xi)\right)\|$. The maximum is attainable due to the compactness of $\mathcal S$.
										\item  $|\left(c(\xi_2) - c(\xi_1)\right)^\top y^*|
										$: Similarly as in the proof of Lemma \ref{lem:SECSTAG-PESSBI-CLP-DUAL}, there are finitely many matrices $L_j^y\in\mathbb R^{m\times n}, \ j=1,\ldots, J_1$, one for each basis of the constraint matrix $A$ such that $\lambda L_j^y b_x(\xi_1)$ is a basic feasible solution in $\left\{y\in\mathbb R^n: \ Ay\le \lambda b_x(\xi_1)\right\}$. Let $L_{j^*}^y$ be the matrix  associated with $y^*$. So
										\begin{equation}\label{eq:lip-bound2}
											|\left(c(\xi_2) - c(\xi_1)\right)^\top y^*| \le \|C(\xi_2 - \xi_1)\| \|\lambda L_{j^*}^yb_x(\xi_1) \|\le L_2 \|\xi_2 - \xi_1\|,
										\end{equation}
										where $L_2:= |\lambda| \|C\|\max_{1\le j\le J_1}\max_{\xi\in\mathcal S} \|L_j^yb_x(\xi)\|$ and the maximum is attainable due to the compactness of $\mathcal S$.
										\item $|\phi_\lambda^y(x,\xi_1,\xi_2)|$: Following a similar argument as for the term $|\left(c(\xi_2) - c(\xi_1)\right)^\top y^*|
										$, let $L^y_{\bar{j}}$ be the matrix associated with an optimal $\Delta y$. Then
										\begin{equation}\label{eq:lip-bound3}
											|\phi_\lambda^y(x,\xi_1,\xi_2)| \le \|c(\xi_2)\|\| \lambda L_{\bar j}^yB_x(\xi_2-\xi_1)\|\le L_3\|\xi_2 -\xi_1\|,
										\end{equation}
										where $L_3:= |\lambda|\max_{\xi\in\mathcal S}\|c(\xi)\|\max_{1\le j\le J_1}\|L_j^yB_x\|$.
										\item $|\phi_\lambda^p(x,\xi_1,\xi_2)|$: Let $L_{\bar j}^p$ be the matrix associated with an optimal $\Delta p$. We have
										\begin{equation}\label{eq:lip-bound4}
											|\phi_\lambda^p(x,\xi_1,\xi_2)| \le \|b(\xi_2)\|\|L_{\bar j}^p\left(V-\lambda C\right)(\xi_2 -\xi_1)\|\le L_4 \|\xi_2 -\xi_1\|,
										\end{equation}
										where $L_4:= \max_{\xi\in\mathcal S}\|b(\xi)\|\max_{1\le j\le J_2}\|L_j^p\left(V-\lambda C\right)\|_{HS} $.
									\end{enumerate}
									Combining \eqref{eq:lip-bound0}-\eqref{eq:lip-bound4}, we obtain $|\Phi_{\lambda}(x,\xi_2)-\Phi_\lambda(x,\xi_1)|\le L\|\xi_2-\xi_1\|$, where $L:=L_1 + L_2 + L_3 + L_4$.
									
									Therefore, the proof is completed.
									\Halmos
								\end{proof}

								\section{Proof of Theorem \ref{thm:0-1SDP}: 0-1 SDP Approximation for $\mathcal D_M$}
								\label{sec:proof_thm_sdp}
								
								\begin{proof}{Proof of Theorem \ref{thm:0-1SDP}:}
									The constraint \eqref{eq:DR-LDR-constr1} is conservatively approximated by requiring the right-hand side to be a convex function. Denote $$\tilde{\mathcal Q} := 	\begin{bmatrix}
										Q-\frac{1}{2}\left(B_x^\top P+C^\top Y\right)-\frac{1}{2}\left(B_x^\top P+C^\top Y\right)^\top &\frac{1}{2}(q-B_x^\top p_0-P^\top b_{x0}-C^\top y_0-Y^\top c_0-W^\top\tau)\\
										\frac{1}{2}(q-B_x^\top p_0-P^\top b_{x0}-C^\top y_0-Y^\top c_0-W^\top\tau)^\top&r-b_{x0}^\top p_0-c_0^\top y_0+\tau^\top h
									\end{bmatrix}\in\mathbb S^{k+1},$$ where $\mathbb S^n$ denotes the space of $n\times n$ symmetric matrices.
									Applying the S-lemma \citep[e.g.,][]{polik2007survey}, given a nonnegative variable $\tau \ge 0$, the constraint \eqref{eq:DR-LDR-constr1} is implied by the following SDP constraint: $\tilde{\mathcal Q}\succeq0$,
									which is further equivalent to
									constraints \eqref{eq:DR-LDR-constr1-sdp}-\eqref{eq:DR-LDR-constr1-sdp2}  by replacing the bilinear terms with the auxiliary variables $\Gamma_i, \omega_i, \rho_i$.
									We complete the proof.		\Halmos
								\end{proof}
								
								\section{Proof of Proposition \ref{prop:extraneous}}
								\label{sec:proof-prop-cop}
								\begin{proof}{Proof of Proposition \ref{prop:extraneous}:}~
									\begin{enumerate}
										\item  A feasible solution of \eqref{model:CPP1} can be expressed as 
										\begin{align*}
											{\footnotesize
												\begin{bmatrix}
													1 & (z^*)^\top \\ z^* & Z^* 
												\end{bmatrix} = \sum_{i} \begin{bmatrix}
													{\nu^*_i}^2 & {\nu^*_i}{\bm{\beta}^*_i}^\top \\
													{\nu^*_i}{\bm{\beta}_i^*} & {\bm{\beta}_i^*}{\bm{\beta}_i^*}^\top
												\end{bmatrix} = \begin{bmatrix}
													\sum_i{\nu^*_i}^2 & & \sum_i{\nu^*_i}{\bm{\beta}_i^*}^\top \\
													\sum_i{\nu^*_i}{\bm{\beta}_i^*} & & \sum_i{\bm{\beta}_i^*}{\bm{\beta}_i^*}^\top
												\end{bmatrix},\quad \sum_i{\nu^*_i}{\bm{\beta}_{i,k+1}^*}=1\ \textrm{and}\ \sum_i\bm{\beta}^{*2}_{i,k+1}=1},
										\end{align*}
										where $\bm\beta^*_{i}\in\mathbb R^{k+1}, \ \forall i$ and $\bm\nu^*_i\in\mathbb R, \ \forall i.$ Let $\hat Z = \sum_i\bm{\beta}_i^*{\bm{\beta}_i^*}^\top$. It is easy to verify the feasibility of $Z^*$ to \eqref{model:CCP2}.
										\item A feasible solution of \eqref{model:CCP2} can be expressed as ${\footnotesize \bar Z =  \begin{bmatrix} \Psi^* & {\xi^*}^\top \\ \xi^* & 1 \end{bmatrix}=\sum_i\bm{\zeta}_i^* {\bm{\zeta}_i^*}^\top}$, where $\bm \zeta_i^*\in\mathbb R^{k+1}, \ \forall i$. Let $z^\prime  =\sum_i\zeta_{i,k+1}^*\zeta_i^*$ and $Z^\prime = \sum_i\zeta_i^*{\zeta_i^*}^\top$ and $\begin{bmatrix}
											1 & (z^\prime)^\top \\ z^\prime & Z^\prime 
										\end{bmatrix} $ is feasible to \eqref{model:CPP1}.
									\end{enumerate}
									The proof completes.
									\Halmos
								\end{proof}
								
								\section{Proof of Theorem \ref{thm:0-1COP}}
								\label{sec:proof-thm-cop}
								
								\begin{proof}{Proof of Theorem \ref{thm:0-1COP}:}
									Constraint  \eqref{eq:DR-LDR-constr1-max} is equivalent to
									\begin{equation}\label{eq:DR-LDR-constr1-max-cop}
										r \ge \min\left\{u: \ ue_{k+1}e_{k+1}^\top-\mathcal{Q}\in\mathcal{COP}(\hat \Xi) \right\}.
									\end{equation}
									Given that Problem \eqref{eq:tsnlp-ldr} is a minimization problem and that constraint \eqref{eq:DR-LDR-constr1-max-cop} is equivalent to \eqref{eq:DR-LDR-constr1}, for an optimal solution to \eqref{eq:tsnlp-ldr}, $r=u$. The proof concludes. 
									\Halmos
								\end{proof}

								\section{Proof of Theorem \ref{thm:discrete-worst-case}: Worst-case Distribution for $\mathcal D_\text{dis}$}
								\label{sec:discrete-distribution}

								\begin{proof}{Proof of Theorem \ref{thm:discrete-worst-case}:}
									
									Given a leader's solution $\hat{x}$, the dual of \eqref{eq:discrete-sub-primal} is derived as follows.
									\begin{subequations}\label{eq:discrete-sub-dual}
										\begin{eqnarray}
											\text{\textbf{SP}}_\text{dis}(\hat x): \quad &\max_{\sigma,\mu,\chi,\gamma} &\begin{array}{l}
												\sum_{s=1}^N{\chi^s}^\top v(\xi^s) \\
											\end{array}\\
											&\st & {\gamma_2\Sigma_0+\mu_0\mu_0^\top-\sum_{s=1}^N\gamma^s\xi^s{\xi^s}^\top+\sqrt{\gamma_1}\left(\mu_0 \mu^\top\Sigma_{0}^{1 / 2}+\Sigma_{0}^{1 / 2}\mu \mu_0^\top\right)\succeq0} \label{eq:CONSTR-b-BigMFree_SP}  \\
											& & \sum_{s=1}^N \gamma^s = 1,\quad \sqrt{\gamma_1}\Sigma_0^{1/2}\mu+\mu_0 - \sum_{s=1}^N \gamma^s\xi^s = 0,\quad \|\mu\|_2\leq1\label{eq:CONSTR-c-BigMFree_SP}\\
											& & -A\chi^s + \gamma^s b_{\hat x}(\xi^s)
											\geq 0,\ s=1,\ldots,N \label{eq:CONSTR-d-BigMFree_SP} \\
											&& \gamma^s c(\xi^s) + A^\top\sigma^s = 0,\ s=1,\ldots,N \label{eq:CONSTR-e-BigMFree_SP} \\ 
											&& \sum_{s=1}^N \left({\chi^s}^\top c(\xi^s)+ {\sigma^s}^\top b_{\hat x}(\xi^s)\right)
											\leq 0\label{eq:CONSTR-f-BigMFree_SP}\\
											& & \gamma^s\geq0,\ \sigma^s\geq0, \ s=1,\ldots,N.
										\end{eqnarray}
									\end{subequations}
									Denote the optimal solution $(\sigma^*,\mu^*,\chi^*,\gamma^*)$. We first show that the optimal dual solution $\gamma^*$ satisfy the three conditions required by the ambiguity set $\mathcal D_\text{dis}$, and  then show that the expected value: $\min_{\lambda\ge 0}\mathbb E\left[\Phi_\lambda(\hat x, \xi)\right]$ with respect to the distribution $\{\mathbb P(\xi = \xi^s) = \gamma^{s*}\}_{ s=1,\ldots,N}$ equals to the worst-case expectation.
									\begin{enumerate}
										\item $\sum_{s=1}^N \gamma^{s*}= 1, \ \gamma^{s*}\ge 0 , \ s=1,\ldots,N$: Because $\gamma^*$ is feasible to \eqref{eq:discrete-sub-dual}.
										\item ${\left(\sum_{s=1}^N \gamma^{s*} \xi^s-{\mu}_{0}\right)^{\top} {\Sigma}_{0}^{-1}\left(\sum_{s=1}^N \gamma^{s*} \xi^s-{\mu}_{0}\right) \leq \gamma_{1}}$: 
										\begin{equation*}
											\left(\sum_{s=1}^N \gamma^{s*} \xi^s-{\mu}_{0}\right)^{\top} {\Sigma}_{0}^{-1}\left(\sum_{s=1}^N \gamma^{s*} \xi^s-{\mu}_{0}\right) = \gamma_1{\mu^*}^\top\Sigma_0^{1/2}\Sigma_0^{-1}\Sigma_0^{1/2}\mu^* = \gamma_1  \|\mu^*\|^2 \le \gamma_1.
										\end{equation*} 
										The first  equality is due to the first constraint in \eqref{eq:CONSTR-c-BigMFree_SP}: $\sum_{s=1}^N\gamma^{s*}\xi^s-\mu_0 =\sqrt{\gamma_1}\Sigma_0^{1/2}\mu^*$. The last inequality is because of the last constraint in \eqref{eq:CONSTR-c-BigMFree_SP}. 
										\item ${\sum_{s=1}^N \gamma^{s*} \left[\left({\xi^s}-{\mu}_{0}\right)\left({\xi^s}-{\mu}_{0}\right)^{\top}\right] \preceq \gamma_{2} {\Sigma}_{0}}$:
										\begin{eqnarray*}
											\sum_{s=1}^N\gamma^{s*}\left(\xi^s-\mu_0\right)\left(\xi^s-\mu_0\right)^\top	&& = \sum_{s=1}^N\gamma^{s*}\xi^s{\xi^s}^\top-\mu_0\sum_{s=1}\gamma^{s*}{\xi^s}^\top - \sum_{s=1}^N\gamma^{s*}\xi^s\mu_0^\top+\mu_0\mu_0^\top\\
											&& = \sum_{s=1}^N\gamma^{s*}\xi^s{\xi^s}^\top - \mu_0(\sqrt{\gamma_1}\Sigma_0^{1/2}\mu^*+ \mu_0)^\top-(\sqrt{\gamma_1}\Sigma_0^{1/2}\mu^*+ \mu_0)\mu_0^\top+\mu_0\mu_0^\top\preceq \gamma_2\Sigma_0.
										\end{eqnarray*}
										The second equality holds because of the first constraint in \eqref{eq:CONSTR-c-BigMFree_SP}. The inequality is according to constraint \eqref{eq:CONSTR-b-BigMFree_SP}.
									\end{enumerate}
									Thus, the distribution $\{\mathbb P(\xi = \xi^s) = \gamma^{s*}\}_{ s=1,\ldots,N}\in \mathcal D_\text{dist}$ and the expected value with respect to $\gamma^*$: $$\min_{\lambda\ge 0}\mathbb E_{\gamma^*}\left[\Phi_\lambda (\hat x,\xi)\right]: = \min_{\lambda\ge 0}\sum_{s=1}^N \gamma^{s*} \Phi_\lambda (\hat x,\xi^s) \le \sup_{\mathcal F\in\mathcal D_\text{dis}}\min_{\lambda\ge 0}\mathbb E_F[\Phi_\lambda (\hat x,\xi)].$$
									It remains to show that the expected value $\min_{\lambda \ge 0}\mathbb E_{\gamma^*}\left[\Phi_\lambda (\hat x,\xi)\right]\ge \sup_{\mathcal F\in\mathcal D_\text{dis}}\min_{\lambda\ge 0}\mathbb E_F[\Phi_\lambda (\hat x,\xi)]$, where the equality is according to Theorem \ref{thm:sup-min-interchangeability}.
									%
									To see this, denote $\tilde \lambda, \tilde p^s, \tilde y^s, \ s=1,\ldots,N$ the optimal solution to %
									\begin{eqnarray*}
										\min_{\lambda\ge 0}	\mathbb E_{\gamma^*} [\Phi_\lambda (\hat x,\xi^s)] = \min_{\lambda\ge 0,p,y} && \sum_{s=1}^N \gamma^{s*}\left(b_{\hat x}(\xi^s)^\top p^s + c(\xi^s)^\top y^s\right)\\
										\label{eq:scenario-constr1}	\mbox{s.t.} && A^\top p^s + \lambda c(\xi^s) = v(\xi^s), \ s=1,\ldots,N\\
										\label{eq:scenario-constr2}	&& Ay^s \le \lambda b_{\hat x}(\xi^s), \ p^s \ge 0, \ s=1,\ldots,N.
									\end{eqnarray*}
									We then have $\min_{\lambda \ge 0} \mathbb E_{\gamma^{*}}[\Phi_{\lambda}(\hat x, \xi)] =$
									\begin{eqnarray*}
										\sum_{s=1}^N \gamma^{s*}\left(b_{\hat x}(\xi^s)^\top \tilde p^s+c(\xi^s)^\top \tilde y^s\right)& & \geq \sum_{s=1}^N \left({\chi^{s*}}^\top A^\top \tilde p^s-{\sigma^{s*}}^\top A\tilde y^s\right) \\
										& & \geq \sum_{s=1}^N \left[{\chi^{s*}}^\top \left(v(\xi^s)-\tilde\lambda c(\xi^s)\right) \tilde p^s-\tilde\lambda{\sigma^{s*}}^\top b_{\hat x}(\xi^s)\right]\\
										&  & \geq \sum_{s=1}^N {\chi^{*s}}^\top v(\xi^s) =\sup_{\mathcal F\in\mathcal D_\text{dis}}\min_{\lambda\ge 0}\mathbb E_F[\Phi_\lambda (\hat x,\xi)].
									\end{eqnarray*}
									The first inequality holds because of constraints \eqref{eq:CONSTR-d-BigMFree_SP}-\eqref{eq:CONSTR-e-BigMFree_SP} and the nonnegativity of $\tilde p^s$. The second inequality holds due to constraints \eqref{eq:scenario-constr1}-\eqref{eq:scenario-constr2} and the nonnegativity of $\sigma^{s*}$. The third inequality holds because of \eqref{eq:CONSTR-f-BigMFree_SP} and $\lambda\ge 0$. The equality holds since strong duality holds between \eqref{eq:discrete-sub-primal} and \eqref{eq:discrete-sub-dual} \citep{todd2001semidefinite}.
									
									It follows that $\min_{\lambda \ge 0}\mathbb E_{\gamma^*}\left[\Phi_\lambda (\hat x,\xi)\right]= \sup_{\mathcal F\in\mathcal D_\text{dis}}\min_{\lambda\ge 0}\mathbb E_F[\Phi_\lambda (\hat x,\xi)]$ and $\{\mathbb P(\xi = \xi^s) = \gamma^{s*}\}_{ s=1,\ldots,N}\in \mathcal D_\text{dist}$  characterize the worst-case distribution.
									\Halmos
								\end{proof}	
								\section{Subproblems for Cutting-Plane Algorithms}
								\label{sec:benders-app}
								In the cutting-plane algorithm, at each iteration, we solve the relaxed master problem MP \eqref{eq:RLX-MASTERPB} and obtain the optimal solution $(\hat x, \hat \nu)$. Then the leader's decision $\hat x$ is plugged into a subproblem, which admits attractable SDP formulations. In this section, we present the three subproblems for (1) the 0-1 SDP formulation \eqref{model:DR-PESSBI-POLY} obtained by conservatively approximating the the nonconvex quadratic constraint \eqref{eq:DR-LDR-constr1}, for (2) the 0-1 SDP approximation \eqref{model:DR-PESSBI-POLY-IA_COP} of the 0-1 COP exact reformulation of constraint \eqref{eq:DR-LDR-constr1}, and for the exact 0-1 SDP reformulation \eqref{eq:sdp-discrete} under the discrete ambiguity set.
								\subsection{Subproblems and Optimality Cuts for 0-1 SDPs}
								\subsubsection{Subproblem and optimality cut for  \eqref{model:DR-PESSBI-POLY}}
								\label{sec:sdp1-benders}
								The worst-case second-stage problem is approximated using the following SDP formulation:
								\begin{subequations}\label{eq:sdp-sub-primal}
									\begin{eqnarray}
										\min && r+t \\
										\label{cnstr:DR-PESSBI-POLY-5}	\mbox{s.t.} 	&& {AY+\Lambda W={\lambda B_x},\ \Lambda h-{Ay_0+\lambda b_{x0}}\ge0,\ \Lambda\ge0}			\\
										&& \tilde{\mathcal Q} \succeq 0\\
										&&			 \eqref{cnstr:DR-PESSBI-POLY-2}-\eqref{eq:DR-LDR-sign},\ \eqref{cnstr:DR-PESSBI-POLY-5b}, \ \eqref{eq:DR-LDR-constr1-sdp2}.\nonumber 
									\end{eqnarray}
								\end{subequations}
								The subproblem is the  dual problem of \eqref{eq:sdp-sub-primal}.
								\begin{subequations}\label{model:DR-PESSBI-POLY-FIXX-DUAL}
									\begin{eqnarray}
										\textbf{SP}_\text{SDP}(\hat{x}):&\underset{U,G,H,\zeta,\sigma,\eta,\mu,{\chi}}{\max}&\tr(VG)+{\chi^\top v_0}\\
										&\st&WUB_{\hat{x}}^\top+{W\eta b_{\hat{x}0}^\top}-WGA^\top-h\zeta^\top\ge0,\ \zeta\ge0\\
										&&\tr(B_{\hat{x}}H)\ge\tr(CG)+{\chi^\top c_0+\sigma^\top b_{\hat{x}0}}\\
										&&{B_{\hat{x}}\eta+b_{\hat{x}0}=A\chi+\zeta}\\
										& & {-U+\gamma_{2} \Sigma_{0}+\mu_{0} \mu_{0}^{\top}+\sqrt{\gamma_1}(\mu_0 \mu^\top\Sigma_{0}^{1 / 2}+\Sigma_{0}^{1 / 2}\mu \mu_0^\top)\succeq0}\label{cnstr:DR-PESSBI-POLY-DUAL1}\\
										&&\eta=\mu_0+\sqrt{\gamma_1}\Sigma_0^{1/2}\mu,\ W\eta\ge h,\ \|\mu\|_2\le 1\label{cnstr:DR-PESSBI-POLY-DUAL2}\\
										&&\begin{bmatrix}
											U    &\eta\\\eta^\top&1
										\end{bmatrix}\succeq0 \label{cnstr:DR-PESSBI-POLY-DUAL3}\\
										&&WH+h\sigma^\top\le 0,\ \sigma\ge0 \label{cnstr:DR-PESSBI-POLY-DUAL5}\\
										&&UC^\top+\eta c_0^\top=HA\quad 
										\label{cnstr:DR-PESSBI-POLY-DUAL6}\\
										&&A^\top\sigma+C\eta+c_0=0.\label{cnstr:DR-PESSBI-POLY-DUAL7}
									\end{eqnarray}
								\end{subequations}
								Note that the subproblem does not involve the big-M constants as in the 0-1 SDP formulation \eqref{model:DR-PESSBI-POLY} in Section \ref{sec:0-1SDP}. Let $(\widehat{U},\widehat{G},\widehat{H},\widehat{\zeta},\widehat\sigma,\widehat\eta,\widehat\mu,{\widehat\chi})$ be the optimal solution to the subproblem SP($\hat x$). When the optimal value is great than $\hat \nu$, we generate a valid cut into MP in the following form.
								\begin{proposition}\label{prop:optcut1}
									Let $[\cdot]^+$ ($[\cdot]^-$) denote the element-wise positive (negative) part of a matrix and $M$ be a sufficiently large big-M constant.
									The inequality
									\small
									\begin{align}\label{eq:OPTCUT}
										\begin{array}{lcl}
											&&\nu\ge\tr(V\widehat{G})+{\widehat\chi^\top v_0}\\&& -M\sum_{i=1}^d\left\{ \hat{x}_i\left(\tr[\mathbf 1_{k\times k}([\widehat U]^+ +[\widehat U]^-)]+[\tr(B_i\widehat H){-b_i^\top\widehat\sigma}]^{{+}} + \right. {([B_i\widehat\eta+b_i]^+ + [B_i\widehat\eta+b_i]^-)^\top\mathbf{1}_m +\left([\widehat\eta]^+ +[\widehat\eta]^-\right)^\top\mathbf{1}_k}\right)(1-x_i)
											\\&& \left. +(1-\hat{x}_i)\left(\tr[\mathbf 1_{k\times k}([\widehat U]^+ +[\widehat U]^-)]+[\tr(B_i\widehat H){-b_i^\top\widehat\sigma}]^- +{([B_i\widehat\eta+b_i]^+ + [B_i\widehat\eta+b_i]^-)^\top\mathbf{1}_m +\left([\widehat\eta]^+ +[\widehat\eta]^-\right)^\top\mathbf{1}_k}\right)x_i\right\}
										\end{array}
									\end{align}
									\normalsize
									is a specific optimality cut in the form of $\nu\ge u_l^\top x+a_l$ to the master problem  to solve 0-1 SDP \eqref{model:DR-PESSBI-POLY}.
								\end{proposition}
								\begin{proof}{Proof of Proposition \ref{prop:optcut1}:}
									To see that the inequality \eqref{eq:OPTCUT} is a valid cut, we first present an SDP formulation fo the worst-case second-stage problem with big-M constants. 
									\begin{equation}\label{eq:sub-mccormick-primal}
										\min \left\{r+t: \ \eqref{eq:DR-LDR-constr1-sdp}-\eqref{eq:DR-LDR-constr1-sdp2}, \ \eqref{cnstr:DR-PESSBI-POLY-2}-\eqref{eq:DR-LDR-sign}, \ \eqref{cnstr:DR-PESSBI-POLY-7}-\eqref{cnstr:DR-PESSBI-POLY-5b}	\right\}.
									\end{equation}
									The SDP \eqref{eq:sub-mccormick-primal} is equivalent to \eqref{eq:sdp-sub-primal} with binary $x=\hat{x}$. So their dual problems are equivalent as well. Specifically, the dual of \eqref{eq:sub-mccormick-primal} is 
									\small
									\begin{subequations}\label{eq:sub-mccormick-dual}
										\begin{eqnarray}
											{\displaystyle\max_{\substack{U,G,H,\zeta,\sigma,\eta,\mu,\Pi_{i}^j,\pi_{i}^j,\\ \chi,\phi_{i}^j,\psi_i^{j}}}}&&
											\tr(VG)+\chi^\top v_0 -M\sum_{i=1}^d\Big(\big[{\tr[\mathbf 1_{k\times k}(\Pi_{i}^1+\Pi_{i}^2)]}+\pi_{i}^2+(\phi_{i}^1+\phi_{i}^2)^\top\mathbf{1}_{m}+(\psi_{i}^1+\psi_{i}^2)^\top\mathbf{1}_{k}\big](1-\hat{x}_i)\nonumber \\
											&&+\big[{\mathbf 1_{k\times k}\tr(\Pi_{i}^3+\Pi_{i}^4)}+\pi_{i}^1+(\phi_{i}^3+\phi_{i}^4)^\top\mathbf{1}_{m}+(\psi_{i}^3+\psi_{i}^4)^\top\mathbf{1}_{k}\big]\hat{x}_i\Big)\\
											&&\sum_{i=1}^dW(\Pi_{i}^2-\Pi_{i}^1)B_{i}^\top+\sum_{i=1}^dW(\psi_{i}^2-\psi_{i}^1)b_{i}^\top+{WUB_0^\top}+W\eta b_0^\top-WGA^\top-h\zeta^\top\ge0,\nonumber\\
											&& \zeta\ge0 \label{cnstr:DR-PESSBI-POLY-DUAL4}\\
											&&{U}+\Pi_{i}^1-\Pi_{i}^2-\Pi_{i}^3+\Pi_{i}^4=0,\ i=1,\ldots,d \label{cnstr:DR-PESSBI-POLY-DUAL-PI3}\\
											&&B_0\eta+b_0+\sum_{i=1}^d (\phi_i^2-\phi_i^1) = A\chi+\zeta \label{cnstr:DR-PESSBI-POLY-DUAL8}\\
											&&-\tr(CG)-\chi^\top c_0+\tr(B_0H)-\sigma^\top b_0+\sum_{i=1}^d(\pi_{i}^2-\pi_{i}^4)\ge0\\ 
											&& \tr(B_iH)-b_i^\top\sigma+\pi_{i}^1-\pi_{i}^2-\pi^3_{i}+\pi^4_{i}=0,\ i=1,\ldots,d\label{cnstr:DR-PESSBI-POLY-DUAL-PI1}\\
											&&B_i\eta+b_i+\phi_i^1-\phi_i^2-\phi_i^3+\phi_i^4 = 0,\ i=1,\ldots,d\\
											&&\eta+\psi_i^1-\psi_i^2-\psi_i^3+\psi_i^4 = 0,\ i=1,\ldots,d\\
											&&\pi_{i}^j\ge0,\ \Pi_{i}^j\ge0,\ \phi_{i}^j\ge0,\ \psi_{i}^j\ge0,\ i=1,\ldots,d,\ j=1,\ldots,4 \label{cnstr:DR-PESSBI-POLY-DUAL-PI2}\\
											&&	\eqref{cnstr:DR-PESSBI-POLY-DUAL1}-\eqref{cnstr:DR-PESSBI-POLY-DUAL7}.\nonumber
										\end{eqnarray}
									\end{subequations}
									\normalsize
									Given an optimal solution $(\overline{U},\overline{G},\overline{H},\overline{\zeta},\overline{\sigma},\overline{\eta},\overline{\mu},\overline{\Pi}_{i}^j,\overline{\pi}_{i}^j,{\overline{\chi},\overline{\phi}_{i}^j,\overline{\psi}_{i}^j})$ to \eqref{eq:sub-mccormick-dual}, if the optimal value of \eqref{eq:sub-mccormick-dual} is greater than $\hat\nu$, following strong duality \citep{todd2001semidefinite},  an optimality cut $\nu\ge u_l^\top x+a_l$ to the master problem MP is specified as
									\begin{eqnarray}\label{cnstr:OPTCUT1}
										\nu&&\ge\tr(V\overline{G})+{\overline{\chi}^\top v_0 -M\sum_{i=1}^d\Big(\big[\tr[\mathbf 1_{k\times k}(\overline{\Pi}_{i}^1+\overline{\Pi}_{i}^2)]+\overline{\pi}_{i}^2+(\overline{\phi}_{i}^1+\overline{\phi}_{i}^2)^\top\mathbf{1}_{m}+(\overline{\psi}_{i}^1+\overline{\psi}_{i}^2)^\top\mathbf{1}_{k}\big](1-x_i)}\nonumber\\
										&&{+\big[\tr[\mathbf 1_{k\times k}(\overline{\Pi}_{i}^3+\overline{\Pi}_{i}^4)]+\overline{\pi}_{i}^1+(\overline{\phi}_{i}^3+\overline{\phi}_{i}^4)^\top\mathbf{1}_{m}+(\overline{\psi}_{i}^3+\overline{\psi}_{i}^4)^\top\mathbf{1}_{k}\big]x_i\Big)}.
									\end{eqnarray}
									To show that \eqref{eq:OPTCUT} is a valid cut,
									it remains to prove that the optimal $(\Pi_{i}^j,\pi_i^j,{\phi_i^j,\psi_i^j}),\ i=1,\ldots,d,\ j=1,\ldots,4$, of \eqref{eq:sub-mccormick-dual} can be expressed using $(\widehat{U},\widehat{H},{\widehat\eta,\widehat\sigma})$ which are optimal to SP($\hat x$). To this end, we present the following complementary slackness conditions and feasibility conditions associated with $(\Pi_{i}^j,\pi_{i}^j,{\phi_i^j,\psi_i^j}),\ i=1,\ldots,d,\ j=1,\ldots,4$. 
									\begin{subequations}\label{cnstr:CPLMTRY-SLACK-PI1}
										\begin{eqnarray}
											&&\Pi_{i}^1(B_i^\top TW+(1-\hat{x}_i)M\mathbf{1}_{k\times k}-\Gamma_i)=0,\ \Pi_{i}^2(B_i^\top TW-(1-\hat{x}_i)M\mathbf{1}_{k\times k}-\Gamma_i)=0\\
											&&\Pi_{i}^3(\hat{x}_iM\mathbf{1}_{k\times k}+\Gamma_i)=0,\ \Pi_{i}^4(\hat{x}_iM\mathbf{1}_{k\times k}-\Gamma_i)=0\\
											&&\pi_{i}^1(\theta_i-M\hat{x}_i)=0,\ \pi_{i}^2(\theta_i-\lambda+M(1-\hat{x}_i))=0\\
											&&\pi_{i}^3\theta_i=0,\ \pi_{i}^4(\theta_i-\lambda)=0,\ i=1,\ldots,d\\
											&&{(p_0+(1-\hat{x}_i)M\mathbf{1}_{m}-\omega_i)^\top\phi_{i}^1=0,\ (p_0-(1-\hat{x}_i)M\mathbf{1}_{m} -\omega_i)^\top\phi_{i}^2=0}\\
											&&{(\hat{x}_iM\mathbf{1}_{m}+\omega_i)^\top\phi_{i}^3=0,\ (\hat{x}_iM\mathbf{1}_{m}-\omega_i)^\top\phi_{i}^4 =0}\\
											&&{\big((TW)^\top b_i+(1-\hat{x}_i)M\mathbf{1}_{k}-\rho_i\big)^\top\psi_{i}^1=0,\ \big((TW)^\top b_i-(1-\hat{x}_i)M\mathbf{1}_{k}-\rho_i\big)^\top\psi_{i}^2=0}\\
											&&{(\hat{x}_iM\mathbf{1}_{k}+\rho_i)^\top\psi_{i}^3=0,\ (\hat{x}_iM\mathbf{1}_{k}-\rho_i)^\top\psi_{i}^4 =0}\\
											&&\eqref{cnstr:DR-PESSBI-POLY-DUAL4}-\eqref{cnstr:DR-PESSBI-POLY-DUAL-PI2}.\nonumber
										\end{eqnarray}
									\end{subequations}
									Denote $(\widehat{Q},\widehat q,\widehat r,\widehat t,\widehat Y,\widehat\lambda,\widehat\tau,\widehat T,\widehat\Lambda,\widehat{\Gamma}_i,\widehat{\theta}_i,{\widehat\omega_i,\widehat\rho_i,\widehat p_0,\widehat y_0})$ be one optimal solution to \eqref{eq:sub-mccormick-primal}. Let 
									\begin{eqnarray}\label{eq:OPT-PI3}
											&&	\widehat\Pi_{i}^1=\hat{x}_i [\widehat U]^{-},\ \widehat\Pi_{i}^2=\hat{x}_i [\widehat U]^{+},\ \widehat\Pi_{i}^3=(1-\hat{x}_i) [\widehat U]^{+},\ \widehat\Pi_{i}^4=(1-\hat{x}_i) [\widehat U]^{-}.\\
											&&	\widehat\pi_{i}^1=(1-\hat{x}_i)[\tr(B_i\widehat H){-b_i^\top\widehat\sigma}]^-,\ \widehat\pi_{i}^2=\hat{x}_i[\tr(B_i\widehat H){-b_i^\top\widehat\sigma}]^+,\ \widehat\pi_{i}^3=(1-\hat{x}_i)[\tr(B_i\widehat H){-b_i^\top\widehat\sigma}]^+,\hspace{8mm}\\ &&\widehat\pi_{i}^4=\hat{x}_i[\tr(B_i\widehat H){-b_i^\top\widehat\sigma}]^-,\\
											&&{\widehat\phi_{i}^1=\hat{x}_i[B_i\widehat\eta+b_i]^-,\ \widehat\phi_{i}^2=\hat{x}_i[B_i\widehat\eta+b_i]^+,\ \widehat\phi_{i}^3=(1-\hat{x}_i)[B_i\widehat\eta+b_i]^+,\ \widehat\phi_{i}^4=(1-\hat{x}_i)[B_i\widehat\eta+b_i]^-},\\
											\label{eq:OPT-PI7}	&&{\widehat\psi_{i}^1=\hat{x}_i[\widehat \eta]^-,\ \widehat\psi_{i}^2=\hat{x}_i[\widehat \eta]^+,\ \widehat\psi_{i}^3=(1-\hat{x}_i)[\widehat \eta]^+,\ \widehat\psi_{i}^4=(1-\hat{x}_i)[\widehat \eta]^-}.
									\end{eqnarray}
									It is easy to verify that \eqref{eq:OPT-PI3}-\eqref{eq:OPT-PI7} satisfy \eqref{cnstr:CPLMTRY-SLACK-PI1} and thus are optimal to \eqref{eq:sub-mccormick-dual}.  Substituting \eqref{eq:OPT-PI3}-\eqref{eq:OPT-PI7} in  \eqref{cnstr:OPTCUT1}, the optimality cut \eqref{cnstr:OPTCUT1}  is equivalent to 
									\eqref{eq:OPTCUT}.
									\Halmos
								\end{proof}

									\subsubsection{Subproblem and optimality cut for  \eqref{model:DR-PESSBI-POLY-IA_COP}}
									\label{sec:sdp2-benders}
									Denote	 $$\hat{\mathcal{Q}}:=
									\begin{bmatrix}
										\frac{1}{2}\left(B_x^\top P+C^\top Y\right)+\frac{1}{2}\left(P^\top B_x+Y^\top C\right)-Q & & \quad  \frac{1}{2}(B_x^\top p_0+P^\top b_{x0}+C^\top y_0+Y^\top c_0-q)\\
										\frac{1}{2}(B_x^\top p_0+P^\top b_{x0}+C^\top y_0+Y^\top c_0-q)^\top & & b_{x0}^\top p_0+c_0^\top y_0
									\end{bmatrix}\in\mathbb{S}^{(k+1)}.$$
									The worst-case second-stage problem is approximated using the following SDP formulation:
									\begin{subequations}\label{eq:cop-sub-primal}
										\begin{eqnarray}
											\min && r+t \\
											\mbox{s.t.} && re_{k+1}e_{k+1}^\top-\hat{\mathcal{Q}}= \mathcal H^\top U \mathcal H, \ U\ge 0, \ U\in\mathbb S^l
											\\
											&&			 \eqref{cnstr:DR-PESSBI-POLY-2}-\eqref{eq:DR-LDR-sign},\ \eqref{cnstr:DR-PESSBI-POLY-5b},\ \eqref{cnstr:DR-PESSBI-POLY-5}.\nonumber 
										\end{eqnarray}
									\end{subequations}
									The subproblem is the dual of \eqref{eq:cop-sub-primal}:
									\begin{subequations}\label{eq:cop-sub-dual}
										\begin{eqnarray}
											{\textbf{SP}_\text{COP}}(\hat{x}):\ \max_{U,G,H,E,\zeta,\sigma, \eta,\mu,\chi} && \tr(VG)+\chi^\top v_0 \\
											\st 
											&& -\tr(CG)-\chi^\top c_0+\tr(B_{\hat{x}}H)-\sigma^\top b_{\hat{x}0}\geq0 \\
											&& -WH-h\sigma^\top\geq0 \\
											&& -WGA^\top-h\zeta^\top+\frac{1}{2}W(U+U^\top)B_{\hat{x}}^\top+W\eta b_{\hat{x}0}^\top\geq0\\
											&& -\zeta-A\chi+B_{\hat{x}}\eta+b_{\hat{x}0}=0\\
											&& -HA + \frac{1}{2}(U+U^\top)C^\top+\eta c_0^\top=0,\ A^\top\sigma+C\eta+c_0=0\\
											&& -\frac{1}{2}(U+U^\top)+\gamma_{2} \Sigma_{0}+\mu_{0} \mu_{0}^{\top}+\sqrt{\gamma_1}(\mu_0 \mu^\top\Sigma_{0}^{1 / 2}+\Sigma_{0}^{1 / 2}\mu \mu_0^\top)\succeq0\\
											&& \mu_0+\sqrt{\gamma_1}\Sigma_{0}^{1 / 2}\mu -\eta=0\\
											&& \frac{1}{2}W(U+U^\top)W^\top-h\eta^\top W^\top-W\eta h^\top+hh^\top-\frac{1}{2}(E+E^\top)=0\\
											&& \zeta\geq0,\ \sigma\geq0,\ \|\mu\|_2\leq1,\ E\geq0.
										\end{eqnarray}
									\end{subequations}
									Given an optimal solution $(U^*,G^*,H^*,E^*,\zeta^*,\sigma^*, \eta^*,\mu^*,\chi^*)$ to problem \eqref{eq:cop-sub-dual}, if the optimal value of \eqref{eq:cop-sub-dual} is great than $\hat\nu$, we generate the following optimality cut into the relaxed master problem MP.
									\begin{proposition}
										The inequality
										\small
										\begin{align}
											\hspace{-0.5cm}
											\nu\geq& \tr(VG^*)+{\chi^*}^\top v_0 - M\sum_{i=1}^d(1-\hat{x}_i) \Bigg\{\left[\tr(B_iH^*)-{\sigma^*}^\top b_i\right]^- +\frac{1}{2}\tr\left[\mathbf 1_{k\times k}\left(\left[U^*+{U^*}^\top\right]^+ +\left[U^*+{U^*}^\top\right]^-\right)\right]\nonumber\\ 
											& +\mathbf 1_m^\top\left(\left[B_i\eta^*+b_i\right]^+ +\left[B_i\eta^*+b_i\right]^-\right)+\mathbf 1_k^\top\left(\left[\eta^*\right]^+ +\left[\eta^*\right]^-\right)\Bigg\}x_i -  M\sum_{i=1}^d \hat{x}_i\Bigg\{\left[\tr(B_iH^*)-{\sigma^*}^\top b_i\right]^+  \nonumber\\
											&  +\frac{1}{2}\tr\left[\mathbf 1_{k\times k}\left(\left[U^*+{U^*}^\top\right]^+ +\left[U^*+{U^*}^\top\right]^-\right)\right]+ \mathbf 1_m^\top\left(\left[B_i\eta^*+b_i\right]^+ +\left[B_i\eta^*+b_i\right]^-\right)+\mathbf 1_k^\top\left(\left[\eta^*\right]^+ +\left[\eta^*\right]^-\right)\Bigg\}(1-x_i)  \label{eq:OPTCUT-cop}
										\end{align}
										\normalsize
										is a specific optimality cut in the form of $\nu\ge u_l^\top x+a_l$ to the master problem  to solve 0-1 SDP \eqref{model:DR-PESSBI-POLY-IA_COP}.
									\end{proposition}
									The proof is similar to that of Proposition \ref{prop:optcut1} and the details are omitted for brevity.
									\subsection{Subproblem and optimality cut for \eqref{eq:sdp-discrete}}
									
									The subproblem is $\text{SP}_\text{dis}(\hat x)$ presented in Section \ref{sec:discrete-distribution} of the Appendices.
									Given an optimal solution $(\hat\sigma,\hat\mu,\hat\chi,\hat\gamma)$, the optimality cut is specified in Proposition \ref{prop:discrete-sub}.
									\begin{proposition}\label{prop:discrete-sub}
										The inequality
										\small
										\begin{align}
											\nu \geq& \sum_{s=1}^N{\hat{\chi^s}}^\top v(\xi^s)-M\sum_{i=1}^d\Bigg\{\hat x_i \left[\sum_{s=1}^N ([\hat{\gamma^s}(B_i\xi^s+b_i)]^- +[\hat{\gamma^s}(B_i\xi^s+b_i)]^+)^\top \mathbf{1}_m+ [-\sum_{s=1}^N  {\hat{\sigma^s}}^\top(B_i\xi^s+b_i)]^+\right](1-{x_i}) \nonumber \\ 
											&+ (1-\hat x_i)\left[\sum_{s=1}^N ([\hat{\gamma^s}(B_i\xi^s+b_i)]^+ +[\hat{\gamma^s}(B_i\xi^s+b_i)]^-)^\top \mathbf{1}_m+[-\sum_{s=1}^N {\hat{\sigma^s}}^\top(B_i\xi^s+b_i)]^-\right]{x_i}\Bigg\} \label{cnstr:OPTCUT_Discrete}
										\end{align}
										\normalsize
										is a specific optimality cut in the form of $\nu\ge u_l^\top x+a_l$ to the master problem  to solve 0-1 SDP \eqref{eq:sdp-discrete}.
									\end{proposition}
									The proof is similar to that of Proposition \ref{prop:optcut1} and thus is omitted for brevity.
									
									\section{Computational Performance with $(\gamma_1,\gamma_2)=(1,1)$}
									\label{sec:cpu_11}
									
									Table \ref{tab:gap-11}  reports the 25\%, 50\% and 75\% quantiles of the optimality gap for the three approximation approaches using linear decision rules when $(\gamma_1,\gamma_2)=(1,1)$.   Table \ref{tab:cpu_11} summarizes the computational performance, across the same test instances as those reported in Table \ref{tab:gap-11}.

									\begin{table}[htbp]
										\centering
										\caption{Quantiles of optimality gaps with $(\gamma_1,\gamma_2) = (1,1)$ }
										\resizebox{.95\textwidth}{!}{%
											\begin{tabular}{c|rrr|rrr|rrr|rrr|rrr|rrr}
												\toprule
												\multirow{3}[2]{*}{Setting} & \multicolumn{9}{c|}{$C=V=0$}                                             & \multicolumn{9}{c}{$C, V \neq 0$} \\
												& \multicolumn{3}{c|}{SDP} & \multicolumn{3}{c|}{IA-COP} & \multicolumn{3}{c|}{MILP-Cut} & \multicolumn{3}{c|}{SDP} & \multicolumn{3}{c|}{IA-COP} & \multicolumn{3}{c}{MILP-Cut} \\
												& \multicolumn{1}{c}{25\%-Q} & \multicolumn{1}{c}{50\%-Q} & \multicolumn{1}{c|}{75\%-Q} & \multicolumn{1}{c}{25\%-Q} & \multicolumn{1}{c}{50\%-Q} & \multicolumn{1}{c|}{75\%-Q} & \multicolumn{1}{c}{25\%-Q} & \multicolumn{1}{c}{50\%-Q} & \multicolumn{1}{c|}{75\%-Q} & \multicolumn{1}{c}{25\%-Q} & \multicolumn{1}{c}{50\%-Q} & \multicolumn{1}{c|}{75\%-Q} & \multicolumn{1}{c}{25\%-Q} & \multicolumn{1}{c}{50\%-Q} & \multicolumn{1}{c|}{75\%-Q} & \multicolumn{1}{c}{25\%-Q} & \multicolumn{1}{c}{50\%-Q} & \multicolumn{1}{c}{75\%-Q} \\
												\hline
												1     & 24.46 & 28.12 & 31.10 & 24.46 & 28.12 & 31.10 & 24.46 & 28.12 & 31.10 & 27.80 & 32.10 & 36.21 & 41.84 & 42.98 & 46.06 & 27.80 & 32.10 & 36.21 \\
												2     & 0.00  & 34.48 & 41.15 & 0.00  & 34.48 & 41.15 & 0.00  & 34.48 & 41.15 & 0.00  & 42.43 & 45.24 & 0.00  & 51.31 & 54.10 & 0.00  & 42.43 & 45.24 \\
												3     & 29.88 & 32.71 & 38.77 & 29.88 & 32.71 & 38.77 & 29.88 & 32.71 & 38.77 & 36.66 & 38.01 & 43.29 & 48.57 & 49.64 & 52.10 & 36.66 & 38.01 & 43.29 \\
												4     & 0.00  & 43.60 & 52.54 & 0.00  & 43.60 & 52.54 & 0.00  & 43.60 & 52.54 & 0.00  & 49.68 & 52.96 & 0.00  & 60.08 & 64.57 & 0.00  & 49.68 & 52.89 \\
												5     & 11.53 & 12.80 & 15.50 & 11.53 & 12.80 & 15.50 & 11.53 & 12.80 & 15.50 & 14.91 & 17.25 & 19.48 & 26.06 & 27.93 & 30.08 & 14.91 & 16.47 & 18.25 \\
												6     & 17.69 & 26.46 & 56.50 & 17.69 & 26.36 & 56.50 & 17.69 & 26.36 & 56.50 & 19.88 & 29.05 & 59.46 & 31.23 & 39.79 & 64.81 & 19.45 & 26.53 & 55.59 \\
												7     & 13.07 & 13.89 & 16.98 & 13.07 & 13.89 & 16.98 & 13.07 & 13.89 & 16.98 & 16.06 & 18.83 & 21.38 & 28.35 & 30.85 & 32.41 & 16.06 & 18.03 & 19.16 \\
												8     & 22.24 & 41.57 & 69.07 & 22.24 & 41.57 & 69.07 & 22.24 & 41.57 & 69.07 & 25.89 & 36.96 & 69.43 & 39.38 & 52.23 & 78.39 & 23.76 & 36.94 & 69.43 \\
												9     & 27.58 & 43.66 & 75.87 & 27.57 & 43.65 & 75.87 & 27.57 & 43.65 & 75.87 & 27.68 & 42.23 & 75.77 & 44.20 & 53.79 & 83.59 & 27.59 & 42.18 & 75.76 \\
												10    & 53.61 & 68.75 & 84.44 & 53.61 & 68.72 & 84.44 & 53.61 & 68.72 & 84.44 & 45.94 & 56.06 & 90.66 & 62.18 & 77.44 & 95.75 & 45.84 & 55.91 & 90.64 \\
												\bottomrule
											\end{tabular}%
										}
										\label{tab:gap-11}%
									\end{table}%
									\begin{table}[htbp]
										\centering
										\caption{Computational comparison of the three approximation approaches with $(\gamma_1,\gamma_2)=(1,1)$}
										\resizebox{\textwidth}{!}{%
											\begin{tabular}{c|rrr|rrr|rrr|rr|rrr|rrr|rrr|rr}
												\toprule
												\multirow{3}[2]{*}{Setting} & \multicolumn{11}{c|}{$C=V=0$}                                                           & \multicolumn{11}{c}{$C, V \neq 0$} \\
												& \multicolumn{3}{c|}{SDP} & \multicolumn{3}{c|}{IA-COP} & \multicolumn{3}{c|}{MILP-Cut} & \multicolumn{2}{c|}{Discrete} & \multicolumn{3}{c|}{SDP} & \multicolumn{3}{c|}{IA-COP} & \multicolumn{3}{c|}{MILP-Cut} & \multicolumn{2}{c}{Discrete} \\
												& \multicolumn{1}{c}{$t_\text{tot}$ (s)} & \multicolumn{1}{c}{\# It.} & \multicolumn{1}{c|}{Gap (\%)} & \multicolumn{1}{c}{$t_\text{tot}$ (s)} & \multicolumn{1}{c}{\# It.} & \multicolumn{1}{c|}{Gap (\%)} & \multicolumn{1}{c}{$t_\text{tot}$ (s)} & \multicolumn{1}{c}{\# It.} & \multicolumn{1}{c|}{Gap (\%)} & \multicolumn{1}{c}{$t_\text{tot}$ (s)} & \multicolumn{1}{c|}{\# It.} & \multicolumn{1}{c}{$t_\text{tot}$ (s)} & \multicolumn{1}{c}{\# It.} & \multicolumn{1}{c|}{Gap (\%)} & \multicolumn{1}{c}{$t_\text{tot}$ (s)} & \multicolumn{1}{c}{\# It.} & \multicolumn{1}{c|}{Gap (\%)} & \multicolumn{1}{c}{$t_\text{tot}$ (s)} & \multicolumn{1}{c}{\# It.} & \multicolumn{1}{c|}{Gap (\%)} & \multicolumn{1}{c}{$t_\text{tot}$ (s)} & \multicolumn{1}{c}{\# It.} \\
												\hline
												1     & \textbf{3.87} & 9.0   & \textbf{27.65} & 4.67  & 8.8   & \textbf{27.65} & 3.99  & 6.5   & 27.66 & 1.36  & 8.8   & 8.29  & 9.0   & \textbf{31.96} & 10.45 & 8.7   & 43.59 & \textbf{4.86} & 6.6   & \textbf{31.69} & 1.58  & 8.9 \\
												2     & 2.12  & 5.0   & \textbf{36.64} & 2.08  & 4.5   & \textbf{36.64} & \textbf{1.78} & 5.7   & \textbf{36.64} & 0.75  & 5.0   & 4.10  & 5.0   & \textbf{39.94} & \textbf{3.75} & 4.5   & 44.93 & 10.48 & 22.0  & \textbf{39.94} & 0.85  & 5.0 \\
												3     & \textbf{3.73} & 9.0   & \textbf{36.47} & 4.65  & 8.9   & \textbf{36.47} & 9.05  & 12.7  & \textbf{36.47} & 1.39  & 8.9   & 8.19  & 9.0   & \textbf{39.32} & 10.76 & 8.9   & 51.98 & \textbf{5.83} & 8.7   & \textbf{39.32} & 1.57  & 8.9 \\
												4     & 2.09  & 5.0   & \textbf{43.07} & 2.14  & 4.6   & \textbf{43.07} & \textbf{1.26} & 4.6   & \textbf{43.07} & 0.76  & 5.0   & 4.13  & 5.0   & 44.99 & \textbf{3.85} & 4.6   & 50.48 & 11.01 & 22.7  & \textbf{44.98} & 0.85  & 5.0 \\
												5     & 9.01  & 9.0   & 13.5  & 8.30  & 9.0   & \textbf{13.43} & \textbf{1.00} & 2.1   & \textbf{13.43} & 3.37  & 9.0   & 45.97 & 9.0   & 17.53 & 91.07 & 8.9   & 28.54 & \textbf{1.88} & 2.2   & \textbf{16.44} & 3.71  & 8.9 \\
												6     & 5.09  & 5.0   & 38.17 & \textbf{3.73} & 4.7   & 37.91 & 15.62 & 19.8  & \textbf{37.43} & 1.73  & 4.6   & \textbf{19.72} & 5.0   & 40.63 & 33.24 & 4.7   & 48.14 & 49.99 & 37.0  & \textbf{37.48} & 2.20  & 4.7 \\
												7     & 9.11  & 9.0   & \textbf{14.55} & 8.08  & 9.0   & \textbf{14.55} & \textbf{1.09} & 2.2   & \textbf{14.55} & 3.30  & 9.0   & 44.10 & 9.0   & 18.91 & 85.55 & 8.9   & 30.52 & \textbf{1.79} & 2.2   & \textbf{17.66} & 3.63  & 8.9 \\
												8     & 5.00  & 5.0   & 47.78 & \textbf{3.60} & 4.7   & 47.78 & 7.79  & 14.4  & \textbf{47.43} & 1.75  & 4.6   & 19.94 & 5.0   & 48.33 & 31.21 & 4.7   & 57.96 & \textbf{17.60} & 16.6  & \textbf{46.01} & 2.18  & 4.7 \\
												9     & 83.78 & 33.0  & 49.95 & \textbf{42.82} & 33.0  & 49.95 & 46.72 & 11.6  & 49.95 & 22.81 & 33.0  & 387.88 & 33.0  & 48.63 & \textbf{80.62} & 33.0  & 58.98 & 163.48 & 20.1  & \textbf{48.60} & 26.47 & 33.0 \\
												10    & 83.79 & 33.0  & 66.45 & \textbf{43.08} & 33.0  & \textbf{66.44} & 57.82 & 13.2  & \textbf{66.44} & 22.90 & 33.0  & 390.89 & 33.0  & 60.19 & \textbf{82.09} & 33.0  & 72.27 & 160.40 & 17.1  & \textbf{60.15} & 27.58 & 33.0 \\
												\bottomrule
											\end{tabular}%
										}
										\label{tab:cpu_11}
									\end{table}%

									\section{Out-of-Sample Performance on Misspecified Distributions}
									\label{sec:oos-misspecified}
									
									Denote $\mu$ and $\sigma^2$ the mean and variance associated with the in-sample uniform distribution over $[30,240]$.
									We generate two out-sample sets with $N^\prime = 5,000$ using the following two types of misspecified distributions, respectively.
									\begin{itemize}
										\item Misspecified moment information: uniform distribution with mean $\lambda^\mu \mu=124$ and variance $\lambda^\sigma \sigma^2=2940$ (i.e., uniform distribution on [30, 218]), where $\mu$ and $\sigma^2$ are the mean and the variance of the in-sample uniform distribution. 
										
										\item Misspecified distribution type:  truncated normal distributions\footnote{See \url{https://en.wikipedia.org/wiki/Truncated_normal_distribution}.} $\mathcal{N}(\mu,\sigma^2)$  over the interval $[30, 240]$.
									\end{itemize}

									\begin{figure}[h]
										\centering
										\begin{subfigure}{0.45\textwidth}
											\centering
											\includegraphics[ width=1.2\linewidth]{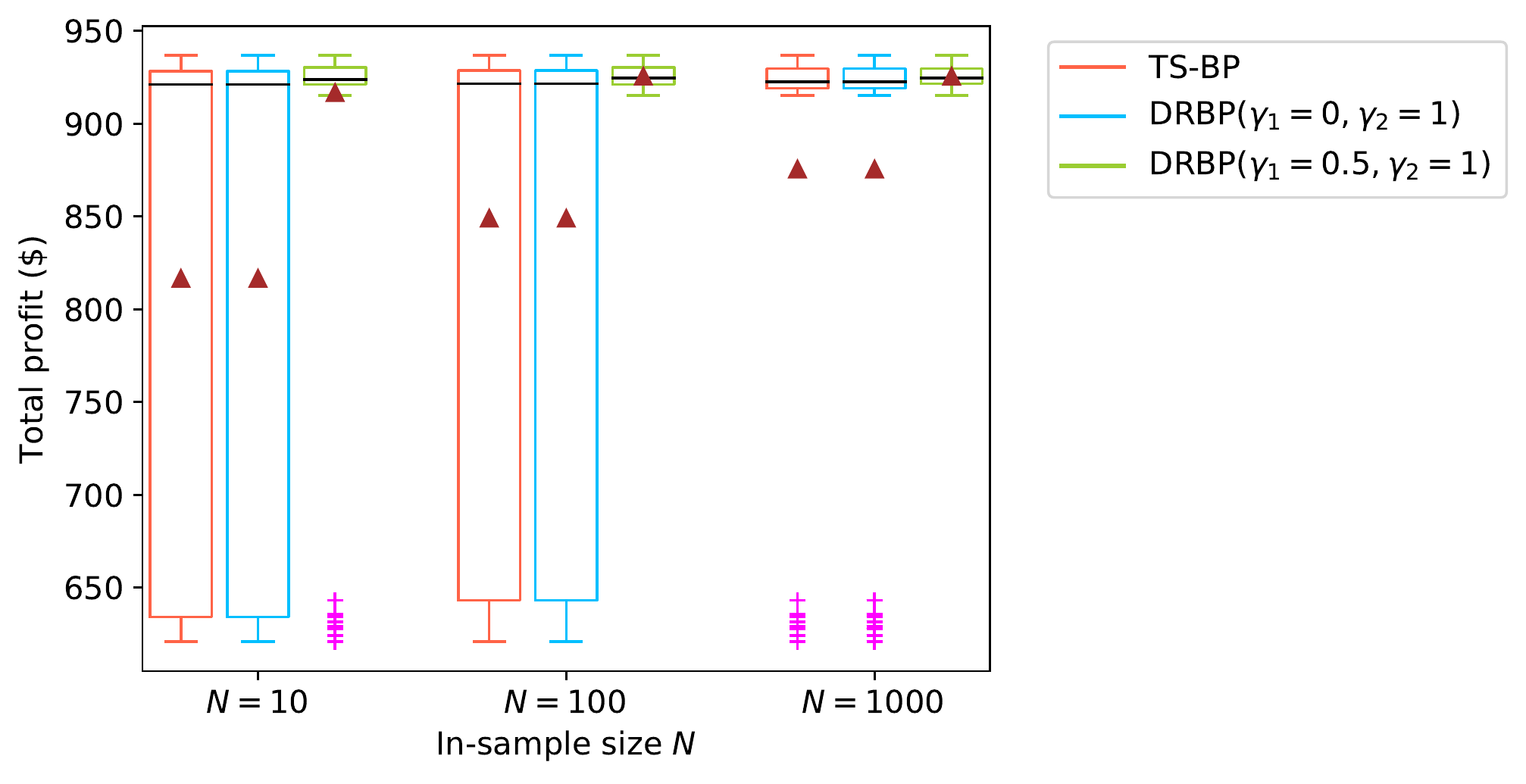}
											\caption{$C=V=0$}
											\label{fig:ofs-0-mis-u}
										\end{subfigure}
										\hfill
										\begin{subfigure}{0.45\textwidth} 
											\centering
											\includegraphics[ width=1.2\linewidth]{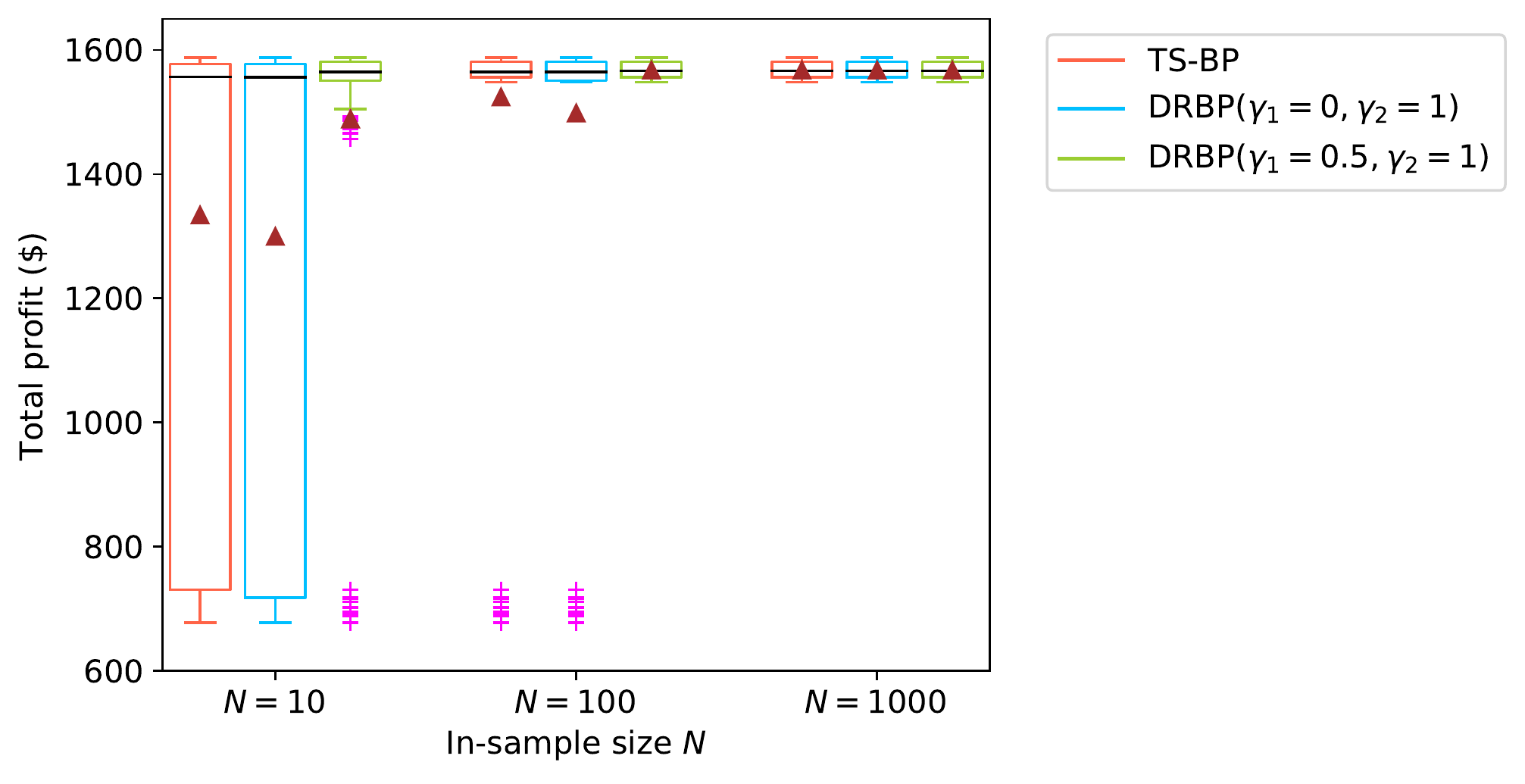}
											\caption{$C,V \neq 0$}
											\label{fig:ofs-n0-mis-u}
										\end{subfigure}
										\centering
										\caption{Out-of-sample performance using misspecified uniform distributions: expected profits}
										\label{fig:out-of-sample-mis-u}
									\end{figure}
									
									\begin{figure}[h]
										\centering
										\begin{subfigure}{0.45\textwidth}
											\centering
											\includegraphics[ width=1.2\linewidth]{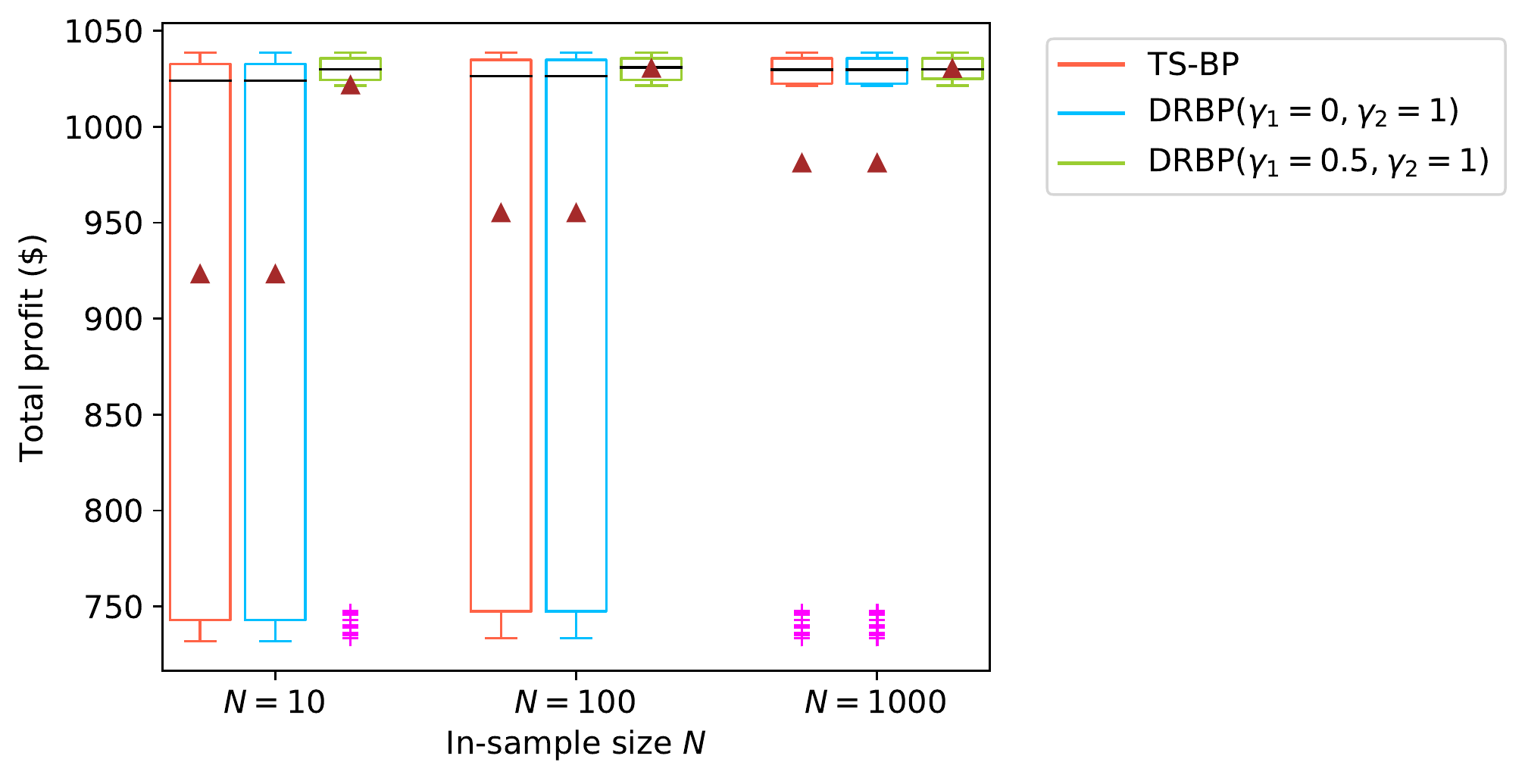}
											\caption{$C=V=0$}
											\label{fig:ofs-0-mis-n}
										\end{subfigure}
										\hfill
										\begin{subfigure}{0.45\textwidth} 
											\centering
											\includegraphics[ width=1.2\linewidth]{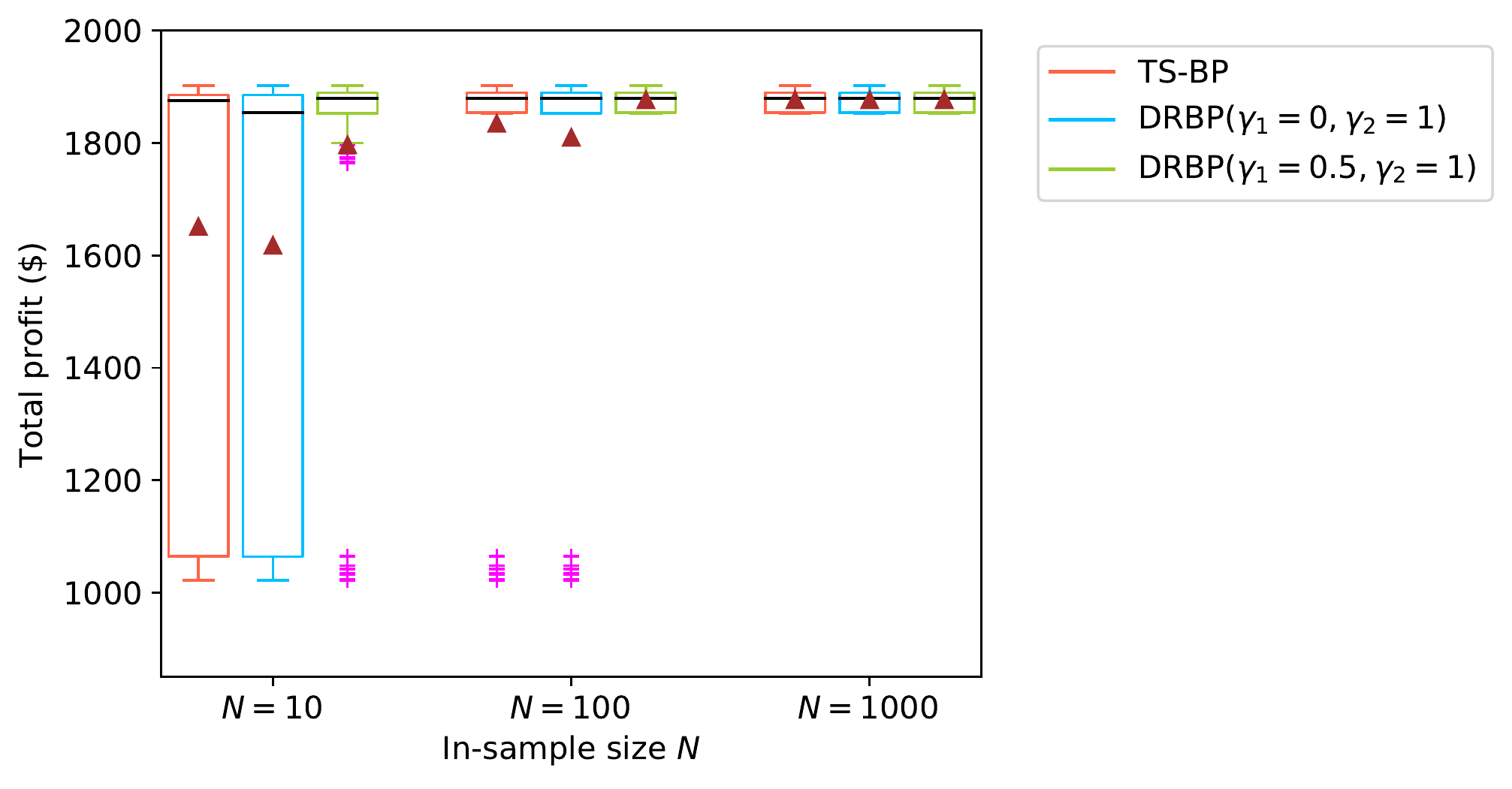}
											\caption{$C,V \neq 0$}
											\label{fig:ofs-n0-mis-n}
										\end{subfigure}
										\centering
										\caption{Out-of-sample performance using misspecified normal distributions: expected profits}
										\label{fig:out-of-sample-mis-n}
									\end{figure}

								\end{APPENDICES}
		\end{document}